\documentclass[a4paper,11pt,leqno]{article}%
\usepackage{graphicx}
\usepackage{amsmath}
\usepackage{amsfonts}
\usepackage{amssymb}%
\newtheorem{theorem}{Theorem}

\newtheorem{corollary}[theorem]{Corollary}

\newtheorem{definition}{Definition}
\newtheorem{example}{Example}

\newtheorem{lemma}[theorem]{Lemma}

\newtheorem{proposition}[theorem]{Proposition}
\newtheorem{remark}{Remark}

\newenvironment{proof}[1][Proof]{\textbf{#1.} }{\ \rule{0.5em}{0.5em}}
\begin{document}

\title{Regularity, Local and Microlocal Analysis in Theories of Generalized Functions}
\author{Jean-Andr\'{e} Marti\\Equipe Analyse Alg\'{e}brique Non Lin\'{e}aire-\textit{Laboratoire} \textit{G
T S I}\\Universit\'{e} Antilles-Guyane}
\maketitle

\begin{abstract}
We introduce a general context involving a presheaf $\mathcal{A}$ and a
subpresheaf $\ \mathcal{B}$ of $\mathcal{A}$. We show that all previously
considered cases of local analysis of generalized functions (defined from
duality or algebraic techniques) can be interpretated as the $\mathcal{B}%
$-local analysis of sections of $\mathcal{A}$.

But the microlocal analysis of the sections of sheaves or presheaves under
consideration is dissociated into a "frequential microlocal analysis " and
into a "microlocal asymptotic analysis". The frequential microlocal analysis
based on the Fourier transform leads to the study of propagation of
singularities under only linear (including pseudodifferential) operators in
the theories described here, but has been extended to some non linear cases in
classical theories involving Sobolev techniques. The microlocal asymptotic
analysis can inherit from the algebraic structure of $\mathcal{B}$ some good
properties with respect to nonlinear operations.

\end{abstract}
\tableofcontents

\section{Introduction}

The notion of regularity in algebras or spaces of generalized functions can be
formulated in a general way with the help of sheaf theory. In section 2, when
$\mathcal{A}$ is a presheaf of algebras or vector spaces on a topological
space $X$, and $\mathcal{B}$ a subpresheaf of $\mathcal{A}$, for each open set
$\Omega$ in $X$, we consider $\mathcal{B}\left(  \Omega\right)  $ as the space
or algebra of some regular elements of $\mathcal{A}\left(  \Omega\right)  $.
This leads to the notion of $\mathcal{B}$-singular support which refines the
notion of support of a section $u\in\mathcal{A}\left(  \Omega\right)  $
provided the localization principle $\left(  F_{1}\right)  $ holds: if $u$ and
$v$ are global sections of $\mathcal{A}$ which agree on each open set of a
family $\left(  \Omega_{i}\right)  _{i\in I}$ of open set in $X$, they agree
on the union $\underset{i\in I}{\cup}\Omega_{i}$.

We can give many examples of this situation in the framework of theories of
generalized functions: distributions \cite{Schwartz1} or Colombeau-type
algebras \cite{AraBia, JFC3, JFC4, JFC5, JFC6}. To illustrate this, let us
consider the following sequence of sheaf embeddings, defined for each $\Omega$
open set in $\mathbb{R}^{n}$ by
\[
\mathrm{C}^{\infty}\left(  \Omega\right)  \rightarrow\mathcal{L}\left(
\mathrm{C}_{c}^{\infty}(\Omega),\mathbb{C}\right)  =\mathcal{D}^{\prime
}\left(  \Omega\right)  \rightarrow\mathcal{G}\left(  \Omega\right)
\rightarrow\mathcal{L}\left(  \mathcal{G}_{c}(\Omega),\widetilde{\mathbb{C}%
}\right)
\]
$\mathcal{G}\left(  \Omega\right)  $ being the Colombeau algebra,
$\widetilde{\mathbb{C}}$ the ring of Colombeau's generalized numbers,
$\mathcal{G}_{c}(\Omega)$ the set of elements in $\mathcal{G}\left(
\Omega\right)  $ with compact support and $\mathcal{L}\left(  \mathcal{G}%
_{c}(\Omega),\widetilde{\mathbb{C}}\right)  $ defined in \cite{Gar} the space
of all continuous and $\widetilde{\mathbb{C}}$ linear functionals on
$\mathcal{G}_{c}(\Omega)$. Each term of the sequence can be considered as the
$\mathcal{B}\left(  \Omega\right)  $ regular space or algebra of the following
algebra or space $\mathcal{A}\left(  \Omega\right)  $. It is the basis for a
local analysis of the elements in $\mathcal{A}\left(  \Omega\right)  $. Some
results on propagation of singularities under $\mathcal{B}$-compatible
operators permit to explain and summarize the classical results involving
differential or pseudo-differential ones. But if we want to define a more
precise "microlocal" analysis which gives some informations not only on the
locus, but on the causis of the singularities described as a fibered space
above that locus, we have first to give a precise local characterization of
the singularities under consideration.

A review on the ideas, technics and results on microlocalization is given in
section 3. The first step was to follow the H\"{o}rmander ideas about the wave
front set $WF(u)$ of a distribution $u$, whose construction is deduced from
the classical Fourier characterization of smoothness of distributions with
compact support. For a general $v\in\mathcal{E}^{\prime}\left(  \mathbb{R}%
^{n}\right)  $ H\"{o}rmander introduces the cone $\Sigma\left(  v\right)  $ of
all $\eta\in\mathbb{R}^{n}\backslash0$ having no conic neighbourhood $V$ such
that the Fourier transform $\widehat{v}$ is rapidly decreasing in $V$. Lemma
8.1.1. in \cite{HorPDOT1} proves that \textit{if }$\Phi\in\mathcal{D}\left(
\mathbb{R}^{n}\right)  $ and $v\in\mathcal{E}^{\prime}\left(  \mathbb{R}%
^{n}\right)  $ then\textit{ }$\Sigma\left(  \Phi v\right)  \subset
\Sigma\left(  v\right)  $. It follows that if $\Omega$ is an open set in
$\mathbb{R}^{n}$ and $u\in\mathcal{D}^{\prime}\left(  \Omega\right)  $,
setting: $\Sigma_{x}(u)=\underset{\Phi}{\cap}\Sigma\left(  \Phi u\right)  ;$
$\Phi\in\mathcal{D}\left(  \Omega\right)  ,$ $\Phi\left(  x\right)  \neq0$,
one can define the wave front set of $u$ as $WF(u)=\left\{  \left(
x,\xi\right)  \in\Omega\times\left(  \mathbb{R}^{n}\backslash0\right)  ;\xi
\in\Sigma_{x}(u)\right\}  $.

This way was led in the sheaf $\mathcal{A}=\mathcal{G}$ of Colombeau
simplified algebras by Nedeljkov, Pilipovic and Scarpalezos \cite{NePiSc} by
taking $\mathcal{B}=\mathcal{G}^{\infty}$ as regular subsheaf of $\mathcal{G}%
$. This subsheaf, introduced by Oberguggenberger \cite{Ober1}, generalizes in
a natural way in $\mathcal{G}$ the regular properties of \textrm{C}$^{\infty}$
in $\mathcal{D}^{\prime}$. We can find in the literature a description of the
main properties of Fourier transform of compacted supported elements in
$\mathcal{G}^{\infty}(\Omega)$, which leads to a frequential microlocal
analysis similar to the H\"{o}rmander's one (see \cite{Scarpa1} for instance).
The crucial point was the conservation of the power of the lemma 8.1.1.
leading to the definition of the generalized wave front set of $u\in
\mathcal{G}(\Omega)$ denoted $WF_{g}(u)$.

Recently, A. Delcroix has extended in \cite{ADRReg} the $\mathcal{G}^{\infty}$
regularity to a so called $\mathcal{G}^{\mathcal{R}}$ regularity, which still
preserves the statements of the above quoted lemma, and gives a $\mathcal{G}%
^{\mathcal{R}}$ frequential microanalysis. We can chose $\mathcal{R}$ such
that $\mathcal{G}^{\mathcal{R}}$ countains an embedding of $\mathcal{D}%
^{\prime}$ into $\mathcal{G}$, which is not the case for $\mathcal{G}^{\infty
}$ ($\mathcal{G}^{\infty}\cap\mathcal{D}^{\prime}=$\textrm{C}$^{\infty}$ is a
result of \cite{Ober1}). Then it becomes possible to investigate the
frequential $\mathcal{D}^{\prime}$-singularities of $u\in\mathcal{G}(\Omega)$.

Inspired by the classical theories, many results on propagation of
singularities and pseudodifferential techniques have been obtained during the
last years by De Hoop, Garetto, H\"{o}rmann, Gramchev, Grosser, Kunzinger,
Steinbauer and others (see \cite{GaGrOb, GarHorm, GKOS, HorDeH, HorKun,
HorObPil}). For example, when $u\in\mathcal{G}(\Omega)$, H\"{o}rmann and
Garetto \cite{GarHorm} obtain characterisations of $WF_{g}(u)$ in terms of
intersections of some domains corresponding to pseudodifferential operators
similarly to H\"{o}rmander's characterizations of $WF(u)$ for $u\in
\mathcal{D}^{\prime}\left(  \Omega\right)  $ \cite{HorFIO}. Following the
ideas and technics of \cite{GarHorm} and making use of the theory of
pseudodifferential operators with generalized symbols (\cite{GaGrOb,
GarHorm}), Garetto \cite{Gar} has recently extended the definition of
$WF_{g}(u)$ when $u\in\mathcal{G}(\Omega)$ to the definitions of
$WF_{\mathcal{G}}(T)$ and $WF_{\mathcal{G}^{\infty}}(T)$ when $T\in
\mathcal{L}\left(  \mathcal{G}_{c}(\Omega),\widetilde{\mathbb{C}}\right)  $.
She can also give a Fourier-transform characterization of these wave front
sets when $T$ is a basic functional. Nevertheless, these very interesting and
deep results are still mainly limited to linear cases, at least in the
framework developed above. More precisely, even when $\mathcal{A}$ is a sheaf
of factor algebras, we don't know any study on the microlocal behaviour of
singularities under nonlinear operations by means of frequential methods based
on the Fourier transform.

However, such studies exist in a classical framework involving some spaces of
Sobolev type. In section 8 of \cite{HorLNLHDE}, H\"{o}rmander uses the results
on microlocal $H_{\left(  s\right)  }^{loc}$-regularity of nonlinear
operations for tempered distributions in $\mathcal{S}^{\prime}\left(
\mathbb{R}^{n}\right)  $ to discuss semi-linear equations, following Rauch
\cite{Rauch}. By means of paradifferential techniques some general results for
quasilinear equations are given in Bony \cite{Bony I}. Then, from a general
result on propagation of singularities for pseudo-differential operators and a
Bony's linearization theorem, H\"{o}rmander (\cite{HorLNLHDE}, section 11) can
discuss fully nonlinear equations and obtain precise propagation results for
hyperbolic second order semi-linear equations. Extensions of the previous
results can be found in works of Beals \cite{Beals I, Beals II} and
Bony\cite{Bony II}.

The Fourier transform is still the main tool involved in other generalized
cases, where the $\mathcal{G}^{\infty}$-regularity is subordinated to an
additional condition (such as an estimate on the growth of derivatives)
characterizing a special property such as to belong to an analytic, Gevrey or
\textrm{C}$^{L}$ class in H\"{o}rmander sense (\cite{HorPDOT1}, section 8.4).
It is the case of "analytic" algebra: $\mathcal{G}^{A}$ studied by Pilipovic,
Scarpalezos and Valmorin \cite{PSV}, of "\textrm{C}$^{L}$ class"algebra:
$\mathcal{G}^{L}$ introduced by Marti \cite{JAM3}, which are subalgebras of
$\mathcal{G}$, of "regular Gevrey ultradistributions" algebra: $\mathcal{G}%
^{\sigma,\infty}$ of Bouzar and Benmeriem \cite{BouBen} which is a subalgebra
of $\mathcal{G}^{\sigma}$, the "generalized Gevrey ultradistributions". In
these examples, the aim is always to perform the $\mathcal{B}$-frequential
microlocalization of generalized functions from the starting algebra
$\mathcal{A}\left(  \Omega\right)  $, when $\mathcal{B}$ is $\mathcal{G}^{A}$,
$\mathcal{G}^{L}$or $\mathcal{G}^{s,\infty}$ and $\mathcal{A}$ is
$\mathcal{G}$ or $\mathcal{G}^{\sigma}$. All these cases are special cases of
the more general one obtained when taking $\mathcal{A=G}^{r}$ which extends
the Colombeau algebra $\mathcal{G}$ and $\mathcal{B}=\mathcal{G}%
^{r,\mathcal{R}\text{,}L}$ which generalizes all the previous regularity
cases. The Fourier transform is still used to characterize the $\mathcal{B}%
$-regularity with the corresponding constraints. But this is not so easy or
natural. For instance, in \cite{JAM3} one starts by giving a characterization
of local $\mathcal{G}^{L}$-regularity by means of some sequence $u_{k}$ of
generalized functions with compact support whose the Fourier transform
$\widehat{u_{k}}$ verifies an estimate involving a special sequence $\left(
L_{k}\right)  _{k\in\mathbb{N}}$. Indeed $u_{k}$ is constructed as product of
$u\in\mathcal{G}(\Omega)$ and a suitable cutoff sequence $\mathcal{X}_{k}$
whose derivatives are controled up to the order $k$. This leads to define the
$\mathcal{G}^{L}$-wave front set of a generalized function: $WF_{g}^{L}\left(
u\right)  \subset\Omega\times(\mathbb{R}^{n}\setminus0)$ and prove, by
refining the cutoff sequence $\mathcal{X}_{k}$, that its projection on
$\Omega$ is the $\mathcal{G}^{L}$-singular support of $\,u$. Then, $WF_{g}%
^{L}\left(  u\right)  $ gives a spectral decomposition of \nolinebreak
sing\thinspace supp$^{L}\,u$. A generalization of these results to the case of
local $\mathcal{G}^{r,\mathcal{R},L}$-regularity of elements in $\mathcal{G}%
^{r}\left(  \Omega\right)  $ is given in \cite{ADJAM}. They lead to define the
$\mathcal{G}^{r,\mathcal{R},L}$-wave front set of $\mathcal{G}^{r}%
$-generalized functions.

These sophisticated constructions give a synthetic description of the
frequential microanalysis but the proofs seem to tell that the Fourier
transform is not really the good tool to perform this description in the above
cases. Perhaps the Fourier-Bros-Iagolnitzer transform would permit to give a
better approach of the problem in the future. For the analytic generalized
wave front set, we also can think of refering to boundary values techniques
which embed distributions into hyperfunctions. But we don't expect results
about nonlinear cases in these ways.

We recall here that generalized functions in the initial definition (sections
of the sheaf $\mathcal{G}$) are classes of families $\left(  u_{\varepsilon
}\right)  _{\varepsilon}$ of classical functions. But in the definitions of
generalized (frequential) wave front set considered above (when $\mathcal{B}$
is $\mathcal{G}^{\infty}$, $\mathcal{G}^{\mathcal{R}}$, $\mathcal{G}^{A}$ or
$\mathcal{G}^{L}$), the parameter $\varepsilon$ does not play a specific role.
It has only to ensure the correct use of some notions as regularity, rapid
decrease, analyticity, and so on in the definition of the different algebras
under consideration. For example the generalized wave front set of any
Dirac-delta function $\Delta$ in the generalized framework is exactly the same
as the classical wave front set of the distribution $\delta$. It is still the
same as the generalized wave front set of any power $\Delta^{m}$ of $\Delta$
without possibility to compare them. The main reason lies in the structure of
Fourier transform. A paradigmatic alternative can be found in the concept of
asymptotic analysis.

The idea of an "asymptotic" analysis \cite{ADJAMMO, JAM0, JAM1} of $u=\left[
u_{\varepsilon}\right]  \in\mathcal{A}\left(  \Omega\right)  =\mathcal{G}%
(\Omega)$ is the following. Let $\mathcal{F}$ be a subsheaf of vector spaces
(or algebras) of $\mathcal{G}$. One defines first the sheaf $\mathcal{B}$ such
that, for any open set $V$ in $\mathbb{R}^{n}$, $\mathcal{B}(V)$ is the space
of elements $u=\left[  u_{\varepsilon}\right]  \in\mathcal{A}(V)$ such that
$u_{\varepsilon}$ has a limit in $\mathcal{F}(V)$ when $\varepsilon$ tends to
$0$. Then $\mathcal{O}_{\mathcal{G}}^{\mathcal{F}}\left(  u\right)  $ is the
set of all $x\in\Omega$ such that $u$ agrees with a section of $\mathcal{B}$
above some neighbourhood of $x$. The $\mathcal{F}$-singular (or $\mathcal{B}%
$-singular) support of $u$ is $\Omega\backslash\mathcal{O}_{\mathcal{G}%
}^{\mathcal{F}}\left(  u\right)  $. For fixed $x$ and $u$, $N_{x}(u)$ is the
set of all $r\in\mathbb{R}_{+}$ such that $\varepsilon^{r}u_{\varepsilon}$
tends to a section of $\mathcal{F}$ above some neighbourhood of $x$. The
$\mathcal{F}$-singular spectrum of $u$ is the set of all $\left(  x,r\right)
\in\Omega\times\mathbb{R}_{+}$ such that $r\in\mathbb{R}_{+}\backslash
N_{x}(u)$. It gives a spectral decomposition of the $\mathcal{F}$-singular
support of $u$. As example, take $\delta_{\varepsilon}(x)=\dfrac
{1}{\varepsilon}\varphi\left(  \dfrac{x}{\varepsilon}\right)  $ where
$\varphi\in\mathcal{D}(\mathbb{R}),\varphi\geq0$ and\ ${\int}\varphi\left(
x\right)  dx=1$. Then, for $m\geq1$, $\Delta^{m}=\left[  \delta_{\varepsilon
}^{m}\right]  $ is a generalized function in $\mathcal{G}(\mathbb{R})$. Except
for $m=1$, for which $\Delta$ is associated with $\delta\in\mathcal{D}%
^{\prime}(\mathbb{R})$, $\Delta^{m}$ is not locally associated with an element
of $\mathcal{D}^{\prime}\left(  V\right)  $ in any neighbourhood $V$ of $0$.
But, for $r\geq m-1$, $\left[  \varepsilon^{r}\right]  \Delta^{m}$ is locally
associated with such an element. It follows that for all $m$, the
$\mathcal{D}^{\prime}$-singular support of $\Delta^{m}$ is $\left\{
0\right\}  $ but its $\mathcal{D}^{\prime}$-singular spectrum is the set
$\left\{  \left(  0,\emptyset\right)  \right\}  $ if $m=1$ or $\left\{
\left(  0,[0,m-1)\right)  \right\}  $ if $m>1$. It gives a more precise
description of the singularities of $\Delta^{m}$ that its $\mathcal{D}%
^{\prime}$-singular support and even that its frequential generalized wave
front set $\left\{  \left(  0,\mathbb{R}\backslash0\right)  \right\}  $ which
doesn't depend upon $m$.

This asymptotic analysis is extended to $(\mathcal{C},\mathcal{E}%
,\mathcal{P})$ algebras.\ This gives the general asymptotic framework, in
which the net $\left(  \varepsilon^{r}\right)  _{\varepsilon}$ is replaced by
any net $a$ satisfying some technical conditions, leading to the concept of
the $(a,\mathcal{F})$-singular parametric spectrum.\ The main advantage is
that this asymptotic analysis is compatible with the algebraic structure of
the presheaf $\mathcal{F}$ asymptotically associated to $(\mathcal{C}%
,\mathcal{E},\mathcal{P})$ algebras.\ Thus the $(a,\mathcal{F})$-singular
asymptotic spectrum inherits good properties with respect to nonlinear
operations when $\mathcal{F}$ is a presheaf of topological algebras. Moreover,
even when $\mathcal{F}$ is a presheaf (or sheaf) of vector spaces (for
instance $\mathcal{F}=\mathcal{D}^{\prime}$), some results on microlocal
analysis are still obtained for nonlinear operations (see paragraph 4.3.1) on
$(a,\mathcal{D}^{\prime})$-singular asymptotic spectrum of powers of $\delta$
functions. In \cite{ADJAMMO}, various examples of propagation of singularities
through nonlinear differential operators are given, connected to some results
of Oberguggenberger, Rauch, Reeds and Travers (\cite{OberBiaVolume, RauchReed,
Travers}).

The paper is organized as follows. In section 2 we introduce the local
analysis of generalized functions. Subsections of 2 give the basic
ingredients, some examples in algebraic or duality theories, and define
$\mathcal{G}^{r,\mathcal{R},L}$-local analysis. $\mathcal{G}$-local analysis
of functional sections of $\mathcal{L}\left(  \mathcal{G}_{c},\widetilde
{\mathbb{C}}\right)  $ and $\mathcal{F}$-local analysis for sections of some
$(\mathcal{C},\mathcal{E},\mathcal{P})$ algebra are recalled. Section 3 is
devoted to the frequential microlocal analysis, with characterization of
$\mathcal{G}^{r,\mathcal{R}}$ and $\mathcal{G}^{r,\mathcal{R},L}$-local
regularities and corresponding wave front sets. We also give the result proved
in \cite{Gar} on the Fourier transform characterisation of $WF_{\mathcal{G}%
}(T)$ when $T$ is a basic functional. The asymptotic microlocal analysis
studied in \cite{ADJAMMO} is detailed in section 4, with examples and
applications to nonlinear partial differential equations.

\section{The local analysis of generalized functions}

The purpose of this section is to localize the singularities of some
generalized functions. We refer the reader to \cite{Godement} for more details
on the sheaf theory involved in the sequel.

\subsection{The basis ingredients}

The basis ingredients of such an analysis are very simple and general ; even
in this subsection no algebraic condition is required.

$\bullet$ $\mathcal{A}$ is a given sheaf of sets (or presheaf with
localization principle $(F_{1})$ in addition) over a topological space $X$.

$\bullet$ $\mathcal{B}$ is a given subsheaf (or subpresheaf) of $\mathcal{A}$.

\begin{definition}
$\mathbf{:}$\textbf{ }$\mathcal{B}$-\textbf{global} \textbf{regularity\medskip
}

For any open set $\Omega$ in $X,$ the elements in $\mathcal{B}\left(
\Omega\right)  $ are considered as regular, and called $\mathcal{B}$-regular
elements of $\mathcal{A}\left(  \Omega\right)  $.
\end{definition}

\begin{definition}
$\mathbf{:}$\textbf{ }$\mathcal{B}$-\textbf{local regularity\medskip}

An element $u\in\mathcal{A}(\Omega)$, where $\Omega$ is any open set in $X$,
is called $\mathcal{B}$-regular at $x\in\Omega$ if there exists an open
neighbourhood $V$ of $x$ such that the restriction $u\mid_{V}$ is in
$\mathcal{B}\left(  V\right)  $.
\end{definition}

\smallskip

\begin{definition}
$\mathbf{:}$\textbf{ }$\mathcal{B}$-\textbf{regular} \textbf{open set\medskip}

We denote by $\mathcal{O}_{\mathcal{A}}^{\mathcal{B}}\left(  u\right)  $ the
set of all $x\in\Omega$ such that $u$ is $\mathcal{B}$-regular at $x$. We also
can write%
\[
\mathcal{O}_{\mathcal{A}}^{\mathcal{B}}\left(  u\right)  =\left\{  x\in
\Omega,\exists V\in\mathcal{V}_{x},u\mid_{V}\in\mathcal{B}\left(  V\right)
\right\}
\]
$\mathcal{V}_{x}$ being the family of all open neighbourhoods of $x$.
\end{definition}

This very simple framework suffices to state the following

\begin{definition}
$\mathbf{:}$\textbf{ }$\mathcal{B}$-\textbf{singular support\medskip}

For any section $u\in\mathcal{A}(\Omega)$, $\Omega$ any open set in $X$, the
$\mathcal{B}$-singular support of $u$ is%
\[
\mathcal{S}_{\mathcal{A}}^{\mathcal{B}}\left(  u\right)  =\Omega
\setminus\mathcal{O}_{\mathcal{A}}^{\mathcal{B}}\left(  u\right)  .
\]

\end{definition}

\begin{remark}
\ $\left(  i\right)  $~The gluing principle $(F_{2})$ is not needed to get the
notion of $\mathcal{B}$-singular support of a section $u\in\mathcal{A}%
(\Omega)$. More precisely, when $\left\{  b\right\}  $ is the constant
presheaf defined by a gobal section of $\mathcal{B}$, the localization
principle $(F_{1})$ is sufficient to prove the following: the set%
\[
\mathcal{O}_{\mathcal{A}}^{\{b\}}\left(  u\right)  =\left\{  x\in
\Omega\ \ \exists V\in\mathcal{V}_{x},\ u\left\vert _{V}\right.
=b\!\!\left\vert _{V}\right.  \right\}
\]
is exactly the union $\Omega_{\mathcal{A}}\left(  u\right)  $ of the open
subsets of $\Omega$ on which $u$ agrees with $b$.\newline Indeed, $(F_{1})$
allows to show that $u$ agrees with $b$ on an open subset $\mathcal{O}$ of
$\Omega$ if, and only if, it agrees with $b$ on an open neighborhood of every
point of $\mathcal{O}$. This leads immediately to the required
assertion.\newline Moreover, $\Omega_{\mathcal{A}}\left(  u\right)
=\mathcal{O}_{\mathcal{A}}^{\{b\}}\left(  u\right)  $ is the largest open set
on which $u$ agrees with $b$, and the $\mathcal{B}$-singular support of $u$ is
a closed subset of its $\{b\}$-singular support $\mathcal{S}_{\mathcal{A}%
}^{\{b\}}\left(  u\right)  =\Omega\setminus\mathcal{O}_{\mathcal{A}}%
^{\{b\}}\left(  u\right)  $.\smallskip

$\left(  ii\right)  $ When the embedding $\mathcal{B}\rightarrow\mathcal{A}$
is a sheaf morphism of abelian groups where $0$ denote the null global
section, $\mathcal{S}_{\mathcal{A}}^{\{0\}}\left(  u\right)  =\Omega
\setminus\mathcal{O}_{\mathcal{A}}^{\{0\}}\left(  u\right)  $ is exactly the
support of $u$ in its classical definition.

$\left(  iii\right)  $~In contrast to the situation described above for the
support or the $\{b\}$-singular support, we need the gluing principle
$(F_{2})$ if we want to prove that the restriction of $u$ to $\mathcal{O}%
_{\mathcal{A}}^{\mathcal{B}}\left(  u\right)  $ belongs to $\mathcal{B}%
(\mathcal{O}_{\mathcal{A}}^{\mathcal{B}}\left(  u\right)  )$. We make this
precise in the following
\end{remark}

\begin{proposition}
Let $u\in\mathcal{A}(\Omega)$. Set $\Omega_{\mathcal{A}}^{\mathcal{B}}\left(
u\right)  =\cup_{i\in I}\Omega_{i},\left(  \Omega_{i}\right)  _{i\in I}$
denoting the collection of all open subsets of $\Omega$ such that $u\left\vert
_{\Omega_{i}}\right.  \in\mathcal{B}\left(  \Omega_{i}\right)  $. Then, if
$\mathcal{B}$ is a sheaf (even if $\mathcal{A}$ is only a prehesaf),

$\left(  i\right)  $~$\Omega_{\mathcal{A}}^{\mathcal{B}}\left(  u\right)  $ is
the largest open subset $\mathcal{O}$ of $\Omega$ such that $u\left\vert
_{\mathcal{O}}\right.  $ is in $\mathcal{B}\left(  \mathcal{O}\right)  $;

$\left(  ii\right)  ~\Omega_{\mathcal{A}}^{\mathcal{B}}\left(  u\right)
=\mathcal{O}_{\mathcal{A}}^{\mathcal{B}}(u)$ and $\mathcal{S}_{\mathcal{A}%
}^{\mathcal{B}}\left(  u\right)  =\Omega\setminus\Omega_{\mathcal{A}%
}^{\mathcal{B}}\left(  u\right)  $.
\end{proposition}

\begin{proof}
$\left(  i\right)  $ For $i\in I$, set $u\left\vert _{\Omega_{i}}\right.
=f_{i}\in\mathcal{B}\left(  \Omega_{i}\right)  $.\ The family $\left(
f_{i}\right)  _{i\in I}$ is coherent by assumption: From $(F_{2})$, there
exists $f\in\mathcal{B}(\Omega_{\mathcal{A}}^{\mathcal{B}}\left(  u\right)  )$
such that $f\left\vert _{\Omega_{i}}\right.  =f_{i}$. But from $(F_{1})$, we
have $f=u$ on $\cup_{i\in I}\Omega_{i}=\Omega_{\mathcal{A}}^{\mathcal{B}%
}\left(  u\right)  $. Thus $u\,|_{\Omega_{\mathcal{A}}^{\mathcal{B}}\left(
u\right)  }\in\mathcal{B}(\Omega_{\mathcal{A}}^{\mathcal{B}}\left(  u\right)
)$, and $\Omega_{\mathcal{A}}^{\mathcal{B}}\left(  u\right)  $ is clearly the
largest open subset of $\Omega$ having this property.

$\left(  ii\right)  $ First, $\mathcal{O}_{\mathcal{A}}^{\mathcal{B}}\left(
u\right)  $ is clearly an open subset of $\Omega$. For $x\in\mathcal{O}%
_{\mathcal{A}}^{\mathcal{B}}\left(  u\right)  $, set $u\left\vert _{V_{x}%
}\right.  =f_{x}\in\mathcal{B}\left(  V_{x}\right)  $ for some suitable
neighborhood $V_{x}$.\ The open set $\mathcal{O}_{\mathcal{A}}^{\mathcal{B}%
}\left(  u\right)  $ can be covered by the family $\left(  V_{x}\right)
_{x\in\mathcal{O}_{\mathcal{A}}^{\mathcal{F}}\left(  u\right)  }$.\ As the
family $\left(  f_{x}\right)  $ is coherent, we get from $(F_{2})$ that there
exists $f\in\mathcal{B}\left(  \cup_{x\in\mathcal{O}_{\mathcal{A}%
}^{\mathcal{B}}\left(  u\right)  }V_{x}\right)  $ such that $f\left\vert
_{V_{x}}\right.  =f_{x}$. From $(F_{1})$, we have $u=f$ on $\cup
_{x\in\mathcal{O}_{\mathcal{A}}^{\mathcal{B}}\left(  u\right)  }V_{x}$ and,
therefore, $u\,|_{\mathcal{O}_{\mathcal{A}}^{\mathcal{B}}\left(  u\right)
}\in\mathcal{B}(\mathcal{O}_{\mathcal{A}}^{\mathcal{B}}\left(  u\right)  ).$
Thus $\mathcal{O}_{\mathcal{A}}^{\mathcal{B}}\left(  u\right)  $ is contained
in $\Omega_{\mathcal{A}}^{\mathcal{B}}\left(  u\right)  $. Conversely, if
$x\in\Omega_{\mathcal{A}}^{\mathcal{B}}\left(  u\right)  $, there exists an
open neighborhood $V_{x}$ of $x$ such that $u\left\vert _{V_{x}}\right.
\in\mathcal{B}\left(  V_{x}\right)  $. Thus $x\in\mathcal{O}_{\mathcal{A}%
}^{\mathcal{B}}\left(  u\right)  $ and the assertion $\left(  ii\right)  $ holds.
\end{proof}

\begin{remark}
When $\mathcal{A}$ is a sheaf and $\mathcal{B}$ a subpresheaf of $\mathcal{A}%
$, we can associate to $\mathcal{B}$ a subsheaf of $\mathcal{A}$ as follows.
When $\Omega$ is an open set of $X$, we note $u\in\mathcal{B}\left(  x\right)
$ if $u\in\mathcal{A}\left(  \Omega\right)  $ is $\mathcal{B}$-regular at $x$
according to Definition 2. Set%
\[
\mathcal{B}_{\ast}\left(  \Omega\right)  =\left\{  u\in\mathcal{A}\left(
\Omega\right)  \mid\forall x\in\Omega\text{ }u\in\mathcal{B}\left(  x\right)
\right\}  \text{.}%
\]

Let $\mathcal{B}_{\ast}$ be the functor $\Omega\rightarrow\mathcal{B}_{\ast
}\left(  \Omega\right)  $. We intend to prove that $\mathcal{B}_{\ast}$ is a
subsheaf of $\mathcal{A}$. The presheaf structure of $\mathcal{A}$ induces
immediately the same one for $\mathcal{B}_{\ast}$ and principle $(F_{1})$ is
also fullfiled. To prove that the gluing principle $(F_{2})$ holds, we
consider a collection $\left(  \Omega_{i}\right)  _{i\in I}$ of open subset
$\Omega_{i}$ of $\Omega$ such that $\Omega=\cup_{i\in I}\Omega_{i}$ and a
coherent family $\left(  u_{i}\right)  _{i\in I}$ of elements $u_{i}%
\in\mathcal{B}_{\ast}\left(  \Omega_{i}\right)  $. First, we can glue the
$u_{i}$ into $u\in\mathcal{A}\left(  \Omega\right)  $. Now, we have to prove
that $u\in\mathcal{B}_{\ast}\left(  \Omega\right)  $. For any $x\in\Omega$
choose $i$ such that $x\in$ $\Omega_{i}$. Then, we have%
\[
\mathcal{B}_{\ast}\left(  \Omega_{i}\right)  \ni u_{i}=u\left\vert
_{\Omega_{i}}\right.  \in\mathcal{A}\left(  \Omega_{i}\right)  \text{.}%
\]
Therefore there exists $V_{i}\subset\Omega_{i}$, $V_{i}\in\mathcal{V}_{x}$
such that $u_{i}\left\vert _{V_{i}}\right.  \in\mathcal{B}\left(
V_{i}\right)  $, from what we deduce%
\[
u\left\vert _{V_{i}}\right.  =\left(  u\left\vert _{\Omega_{i}}\right.
\right)  \left\vert _{V_{i}}\right.  =u_{i}\left\vert _{V_{i}}\right.
\in\mathcal{B}\left(  V_{i}\right)
\]
which proves that $u\in\mathcal{B}\left(  x\right)  $ for each $x\in\Omega$.
Then $\mathcal{B}_{\ast}$ is a sheaf. In fact it is the sheaf associated to
$\mathcal{B}$ in the sense of \cite{Godement}. Its construction is simplified
by using the sheaf structure of $\mathcal{A}$. Roughly speaking,
$\mathcal{B}_{\ast}$ is constructed thanks to a local procedure which adds
many sections to the $\mathcal{B}$ ones. Then, one wishes to compare the
corresponding singular supports of the same $u\in\mathcal{A}\left(
\Omega\right)  $. The answer is given by the following proposition which shows
that the singularities of sections of $\mathcal{A}$ don't decrease when
replacing $\mathcal{B}$ by $\mathcal{B}_{\ast}$.
\end{remark}

\begin{proposition}
Suppose that $\mathcal{A}$ is a sheaf and $\mathcal{B}$ a subpresheaf of
$\mathcal{A}$. Let $\mathcal{B}_{\ast}$ the subsheaf of $\mathcal{A}$
associated to $\mathcal{B}$. Then, for any section $u$ of $\mathcal{A}$ over
the open set $\Omega$ of $X$ we have%
\[
\mathcal{S}_{\mathcal{A}}^{\mathcal{B}}\left(  u\right)  =\mathcal{S}%
_{\mathcal{A}}^{\mathcal{B}_{\ast}}\left(  u\right)  \text{.}%
\]

\end{proposition}

\begin{proof}
For $u\in\mathcal{A}\left(  \Omega\right)  $, the presheaf mapping:
$\mathcal{B}\rightarrow\mathcal{B}_{\ast}$ leads immediately to the set
inclusion: $\mathcal{S}_{\mathcal{A}}^{\mathcal{B}}\left(  u\right)
\supset\mathcal{S}_{\mathcal{A}}^{\mathcal{B}_{\ast}}\left(  u\right)  $.
Conversely, let be $x\in\mathcal{O}_{\mathcal{A}}^{\mathcal{B}_{\ast}}\left(
u\right)  $. There exists $V\in\mathcal{V}_{x}$ such that $u\left\vert
_{V}\right.  \in\mathcal{B}_{\ast}\left(  V\right)  $. But as $x\in V$, there
exists $W\in\mathcal{V}_{x}\cap V$ such that%
\[
u\left\vert _{W}\right.  \in\mathcal{B}\left(  W\right)  .
\]
Then $x\in\mathcal{O}_{\mathcal{A}}^{\mathcal{B}}\left(  u\right)  $. We have
proved the inclusion $\mathcal{O}_{\mathcal{A}}^{\mathcal{B}_{\ast}}\left(
u\right)  \subset\mathcal{O}_{\mathcal{A}}^{\mathcal{B}}\left(  u\right)  $
which gives the converse one for the respective singular supports and leads to
the required equality.
\end{proof}

\subsection{Some properties of $\mathcal{B}$-singular support}

\subsubsection{Elementary algebraic properties}

\begin{proposition}
\label{algebraic}We suppose that $\mathcal{B}$ and $\mathcal{A}$ are
presheaves of $\mathbb{K}$-vector spaces, (resp. algebras). Let $\left(
u_{j}\right)  _{1\leq j\leq p}$ be any finite family of elements in
$\mathcal{A}(\Omega)$ and $\left(  \lambda_{j}\right)  _{1\leq j\leq p}$ any
finite family of elements in $\mathbb{K}$. We have:%
\[
\mathcal{S}_{\mathcal{A}}^{\mathcal{B}}(%
{\textstyle\sum\limits_{1\leq j\leq p}}
\lambda_{j}u_{j})\subset%
{\textstyle\bigcup\limits_{1\leq j\leq p}}
\mathcal{S}_{\mathcal{A}}^{B}(u_{j}).
\]
In the resp. case, we have in addition:%
\[
\mathcal{S}_{\mathcal{A}}^{\mathcal{B}}(%
{\textstyle\prod\limits_{1\leq j\leq p}}
u_{j})\subset%
{\textstyle\bigcup\limits_{1\leq j\leq p}}
\mathcal{S}_{\mathcal{A}}^{B}(u_{j})\text{.}%
\]
In particular, if $u_{j}=u$ for $1\leq j\leq p$, we have $\mathcal{S}%
_{\mathcal{A}}^{\mathcal{B}}(u^{p})\subset\mathcal{S}_{\mathcal{A}}^{B}(u).$
\end{proposition}

\begin{proof}
If $x\in\Omega$ is in $\underset{1\leq j\leq p}{\cap}$ $\mathcal{O}%
_{\mathcal{A}}^{B}(u_{j})$, there exists $V_{j}$ in $\mathcal{V}_{x}$ such
that $u_{j}\left\vert _{V_{j}}\right.  \in\mathcal{B}(V_{j})$. Thus $\left(
{\textstyle\sum\limits_{1\leq j\leq p}}
\lambda_{j}u_{j}\right)  \left\vert _{\underset{1\leq j\leq p}{\cap}V_{j}%
}\right.  \in\mathcal{B}(\underset{1\leq j\leq p}{\cap}V_{j})$ (resp. $\left(
%
{\textstyle\prod\limits_{1\leq j\leq p}}
u_{j}\right)  \left\vert _{\underset{1\leq j\leq p}{\cap}V_{j}}\right.
\in\mathcal{B}(\underset{1\leq j\leq p}{\cap}V_{j})$, which implies%
\[
\underset{1\leq j\leq p}{\cap}\mathcal{O}_{\mathcal{A}}^{B}(u_{j}%
)\subset\mathcal{O}_{\mathcal{A}}^{\mathcal{B}}(%
{\textstyle\sum\limits_{1\leq j\leq p}}
\lambda_{j}u_{j})\ \ \ \ \ \text{(resp.\ }\underset{1\leq j\leq p}{\cap
}\mathcal{O}_{\mathcal{A}}^{B}(u_{j})\subset\mathcal{O}_{\mathcal{A}%
}^{\mathcal{B}}\left(
{\textstyle\prod\limits_{1\leq j\leq p}}
u_{j}\right)  \text{.}%
\]
The result follows by taking the complementary sets in $\Omega$.\medskip
\end{proof}

\subsubsection{$\mathcal{B}$-compatible operators and propagation of
singularities}

We begin by a general result which doesn't need algebraic assumptions. Let
$\Omega$ be a given open subset of $X$. A presheaf operator $A$ in
$\mathcal{A}\left(  \Omega\right)  $ is defined as a presheaf morphism
$\mathcal{A}\left(  \Omega\right)  \rightarrow\mathcal{A}\left(
\Omega\right)  $ compatible with restrictions. More precisely if
$\mathcal{O}_{\Omega}$ denote the category of all open sets in $\Omega$, $A$
may be given by a collection $\left(  A_{V}\right)  _{V\in\mathcal{O}_{\Omega
}}$ of mappings $A_{V}:\mathcal{A}\left(  V\right)  \rightarrow\mathcal{A}%
\left(  V\right)  $ such that for each $V\in\mathcal{O}_{\Omega}$ and
$u\in\mathcal{A}\left(  \Omega\right)  $ we have: $A_{\Omega}\left(  u\right)
\left\vert _{V}\right.  =A_{V}\left(  u\left\vert _{V}\right.  \right)  $.
Thus we can simplify the notations and write $A$ instead $A_{V}$ when acting
on sections over $V$.

\begin{definition}
Let $A$ be a presheaf operator in $\mathcal{A}\left(  \Omega\right)  $. We say
that $A$ is locally $\mathcal{B}$-compatible if for each triple $\left(
x,V,v\right)  \in\Omega\times\mathcal{V}_{x}\times\mathcal{B}\left(  V\right)
$ there exists $W\in\mathcal{V}_{x}$, $W\subset V$, such that $A\left(
v\right)  \left\vert _{W}\right.  \in\mathcal{B}\left(  W\right)  $.
\end{definition}

\begin{proposition}
\label{operator}Suppose that the assumption given in subsection 2.1\ are
fulfilled. Let $A$ be a presheaf operator in $\mathcal{A}$ $\left(
\Omega\right)  $ locally $\mathcal{B}$-compatible. Then we have%
\[
\mathcal{S}_{\mathcal{A}}^{\mathcal{B}}(A\left(  u\right)  )\subset
\mathcal{S}_{\mathcal{A}}^{B}(u).
\]

\end{proposition}

\begin{proof}
If $x\in\Omega$ belongs to $\mathcal{O}_{\mathcal{A}}^{B}(u)$ there exists $V$
in $\mathcal{V}_{x}$ such that $u\left\vert _{V}\right.  \in\mathcal{B}(V)$.
Then, there exists $W\in\mathcal{V}_{x}$, $W\subset V$, such that $A\left(
u\left\vert _{V}\right.  \right)  \left\vert _{W}\right.  \in\mathcal{B}%
\left(  W\right)  $. We have%
\[
A\left(  u\left\vert _{V}\right.  \right)  \left\vert _{W}\right.  =\left(
A\left(  u\right)  \left\vert _{V}\right.  \right)  \left\vert _{W}\right.
=A\left(  u\right)  \left\vert _{V\cap W}\right.  =A\left(  u\right)
\left\vert _{W}\right.  \text{.}%
\]

Then $x$ belongs to $\mathcal{O}_{\mathcal{A}}^{B}(A\left(  u\right)  )$, and
we have proved that $\mathcal{O}_{\mathcal{A}}^{B}(u)\subset\mathcal{O}%
_{\mathcal{A}}^{B}(A\left(  u\right)  )$. The result follows by taking the
complementary sets in $\Omega$.\medskip
\end{proof}

The following weakened form of locally $\mathcal{B}$-compatibility may be more
practical for applications

\begin{definition}
\label{B-compatible}A presheaf operator $A$ in $\mathcal{A}\left(
\Omega\right)  $ is said $\mathcal{B}$-compatible if for each open set $V$ of
$\Omega$ it maps $\mathcal{B}\left(  V\right)  $ into itself.
\end{definition}

It is easy to see that a $\mathcal{B}$-compatible operator is locally
$\mathcal{B}$-compatible but with the above definition we can get some useful
results. The simplest one concerns the composition product, with an obvious proof.

\begin{proposition}
If a presheaf operator $A$ in $\mathcal{A}\left(  \Omega\right)  $ is
$\mathcal{B}$-compatible, then for any $p\in N$, the composition product
$A^{p}=\overset{p}{\overbrace{A\circ A\circ...A}}$ is $\mathcal{B}$-compatible.
\end{proposition}

Adding some algebraic hypothesis leads to the following definitions and
results as corollaries of propositions \ref{algebraic} and \ref{operator}.
When $\mathcal{B}$ is a presheaf of algebras and $\mathcal{A}$ a presheaf of
vector spaces and a $\mathcal{B}$-module we recall that the external sheaf
product $\mathcal{B\times A\rightarrow A}$ extending the usual algebra product
$\mathcal{B\times B\rightarrow B}$ is defined for any $\Omega\in
\mathcal{O}_{X}$ and $\left(  b,u\right)  \in\mathcal{B}\left(  \Omega\right)
\mathcal{\times A}\left(  \Omega\right)  $ by $\left(  b,u\right)  \mapsto bu$
$\in\mathcal{A}\left(  \Omega\right)  $ with $bu\left\vert _{V}\right.
=b\left\vert _{V}\right.  u\left\vert _{V}\right.  $ for each open set $V$ in
$\Omega$.

\begin{definition}
We suppose that $\mathcal{B}$ is a presheaf of algebras and $\mathcal{A}$ a
presheaf of vector space and a $\mathcal{B}$-module. Let $b$ be a given
element in $\mathcal{B}\left(  \Omega\right)  $. We define the $B$-operator of
multiplication in $\mathcal{A}\left(  \Omega\right)  $ by the map
$\mathcal{A}\left(  \Omega\right)  \rightarrow\mathcal{A}\left(
\Omega\right)  $ such that $B\left(  u\right)  =bu$.
\end{definition}

\begin{proposition}
$B$ is a presheaf operator $\mathcal{B}$-compatible.
\end{proposition}

\begin{proof}
First, we have for each open set $V$ in $\Omega$%
\[
B\left(  u\right)  \left\vert _{V}\right.  =bu\left\vert _{V}\right.
=b\left\vert _{V}\right.  u\left\vert _{V}\right.  =B\left(  u\left\vert
_{V}\right.  \right)  .
\]

When $v$ is in $\mathcal{B}\left(  V\right)  $, the external product
$b\left\vert _{V}\right.  v$ agrees with the standard product in the algebra
$\mathcal{B}\left(  V\right)  $ and lies in it. Then, for each pair $\left(
V,v\right)  \in\mathcal{O}_{\Omega}\times\mathcal{B}\left(  V\right)  $,
$B\left(  v\right)  $ is in $\mathcal{B}\left(  V\right)  $.
\end{proof}

\begin{corollary}
\label{sum}We keep the same assumption as above and consider a family $\left(
b_{\alpha}\right)  _{\alpha\in\mathfrak{A}}$ of elements in $\mathcal{B}%
\left(  \Omega\right)  $ and another family of $\mathcal{B}$-compatible
operators $\left(  A_{\alpha}\right)  _{\alpha\in\mathfrak{A}}$ in
$\mathcal{A}\left(  \Omega\right)  $ where $\mathfrak{A}$ is a set of indices.
Then, for any finite$\mathfrak{\ }$part $\mathfrak{A}_{0}$ of $\mathfrak{A}$,
$\underset{\alpha\in\mathfrak{A}_{0}}{%
{\textstyle\sum}
}b_{\alpha}A_{\alpha}$ is a $\mathcal{B}$-compatible operator in
$\mathcal{A}\left(  \Omega\right)  $.
\end{corollary}

\begin{proof}
For each $\alpha\in\mathfrak{A}_{0}$ and $V$ open set in $\Omega$ we have
\[
b_{\alpha}A_{\alpha}\left(  u\right)  \left\vert _{V}\right.  =b_{\alpha
}\left\vert _{V}\right.  A_{\alpha}\left(  u\right)  \left\vert _{V}\right.
=b_{\alpha}\left\vert _{V}\right.  A_{\alpha}\left(  u\left\vert _{V}\right.
\right)  =b_{\alpha}A_{\alpha}\left(  u\left\vert _{V}\right.  \right)  .
\]

When $v$ belongs to $\mathcal{B}\left(  V\right)  $, $A_{\alpha}\left(
v\right)  $ belongs to $\mathcal{B}\left(  V\right)  $ from the hypothesis on
$A_{\alpha}$. Then, the external product $b_{\alpha}\left\vert _{V}\right.  $
$A_{\alpha}\left(  v\right)  $ agrees with the standard product in the algebra
$\mathcal{B}\left(  V\right)  $ and lies in it. Then, for each pair $\left(
V,v\right)  \in\mathcal{O}_{\Omega}\times\mathcal{B}\left(  V\right)  $,
$b_{\alpha}A_{\alpha}\left(  v\right)  $ is in $\mathcal{B}\left(  V\right)  $
and it is the same for the finite sum$\underset{\alpha\in\mathfrak{A}_{0}%
}{\text{ }%
{\textstyle\sum}
}b_{\alpha}A_{\alpha}\left(  v\right)  $.
\end{proof}

\begin{corollary}
Let $P$ be the polynomial in $\mathcal{A}\left(  \Omega\right)  $ defined for
each $u\in\mathcal{A}\left(  \Omega\right)  $ by $P\left(  u\right)  =%
{\textstyle\sum\limits_{1\leq j\leq p}}
b_{j}u^{j}$ where $b_{j}\in\mathcal{B}\left(  \Omega\right)  $. We suppose
that $\mathcal{B}$ and $\mathcal{A}$ are presheaves of algebras. Then $P$ is
$\mathcal{B}$-compatible.
\end{corollary}

\begin{proof}
It suffices to remark that for each $j\in\mathbb{N}$ the map $A_{j}:u\mapsto
u^{j}$ is a $\mathcal{B}$-compatible presheaf operator of $\mathcal{A}\left(
\Omega\right)  $. Putting $\alpha=j$ in the above corollary gives the result.
\end{proof}

Collecting all these informations we can summarize the previous results in the following:

\begin{proposition}
\label{product}We suppose that $\mathcal{B}$ is a presheaf of algebras and
$\mathcal{A}$ a presheaf of vector spaces and a $\mathcal{B}$-module.
$\mathfrak{A}$ being a set of indices, let $\left(  b_{\alpha}\right)
_{\alpha\in\mathfrak{A}}$ be a family of elements in $\mathcal{B}\left(
\Omega\right)  $, $\left(  A_{\alpha}\right)  _{\alpha\in\mathfrak{A}}$ a
family of $\mathcal{B}$-compatible operators in $\mathcal{A}\left(
\Omega\right)  $ and $\left(  p_{\alpha}\right)  _{\alpha\in\mathfrak{A}}$ a
family of positive integers. Then,

$\left(  i\right)  $ For any finite$\mathfrak{\ }$part $\mathfrak{A}_{0}$ of
$\mathfrak{A}$, $\underset{\alpha\in\mathfrak{A}_{0}}{%
{\textstyle\sum}
}b_{\alpha}A_{\alpha}^{p_{\alpha}}$ is a $\mathcal{B}$-compatible operator in
$\mathcal{A}\left(  \Omega\right)  $.

$\left(  ii\right)  $ If $\mathcal{A}$ is a presheaf of algebras,
$u\mapsto\underset{\alpha\in\mathfrak{A}_{0}}{%
{\textstyle\sum}
}b_{\alpha}\left(  A_{\alpha}\left(  u\right)  \right)  ^{p_{\alpha}}$ is a
$\mathcal{B}$-compatible operator in $\mathcal{A}\left(  \Omega\right)  $.
\end{proposition}

\subsection{Examples}

\begin{example}
$:\mathrm{C}^{\infty}$\textbf{-}$\mathbf{local}$ $\mathbf{analysis}$
$\mathbf{in}$ $\mathcal{D}^{\prime}\medskip\left(  \Omega\right)  $

Let $\mathcal{A}=\mathcal{D}^{\prime}$, $\mathcal{B}=\mathrm{C}^{\infty}$.
Then, for any distribution $u\in\mathcal{D}^{\prime}\left(  \Omega\right)  $
where $\Omega$ is an open set of $\mathbb{R}^{n}$%
\[
\mathcal{S}_{\mathcal{D}^{\prime}}^{\mathrm{C}^{\infty}}\left(  u\right)
=sing\,supp\left(  u\right)
\]
where sing\thinspace supp\thinspace$u$ is, in the H\"{o}rmander sense, the
closet subset aff all $x\in\Omega$ having no neighbourhood in which the
distribution $u$ is smooth.
\end{example}

\begin{example}
$:\mathcal{G}^{\infty}$\textbf{-}$\mathbf{local}$ $\mathbf{analysis}$
$\mathbf{in}$\textbf{ }$\mathcal{G}$\textbf{\medskip}$\left(  \Omega\right)  $

Let $\mathcal{A}=\mathcal{G}$ and $\mathcal{B}=\mathcal{G}^{\infty}$ the
"regular" subsheaf of $\mathcal{G}$, the sheaf of Colombeau's generalized
functions. $\mathcal{G}^{\infty}\left(  \Omega\right)  $ is defined as the
sections of $\mathcal{G}\left(  \Omega\right)  $ having a representative
verifying%
\[
\forall K\Subset\Omega\text{ }\exists p\geq0\text{ }\forall\alpha\in
\mathbb{N}^{n}\text{ }\underset{x\in K}{\text{sup}}\left\vert \partial
^{\alpha}u_{\varepsilon}\left(  x\right)  \right\vert =\mathcal{O}\left(
\varepsilon^{-p}\right)  \text{ as }\varepsilon\rightarrow0.
\]
Then%
\[
\mathcal{S}_{\mathcal{G}}^{\mathcal{G}^{\infty}}\left(  u\right)
=sing\,supp_{g}\left(  u\right)
\]
where sing\thinspace supp$_{g}$\thinspace$u$ is the generalized singular
support of $u$ defined in the literature as the set of all $x\in\Omega$ having
no neighbourhood $V$ such that $u\left\vert _{V}\right.  \in\mathcal{G}%
^{\infty}\left(  V\right)  $.
\end{example}

\begin{example}
$:\mathcal{G}^{\mathcal{R}}$-$\mathbf{local}$ $\mathbf{analysis}$
$\mathbf{in}$\textbf{ \medskip}$\mathcal{G}\left(  \Omega\right)  $

In \cite{ADRReg} the $\mathcal{G}^{\infty}$-regularity is extended into a
$\mathcal{G}^{\mathcal{R}}$ one. Starting from a set $\mathcal{R}$ of
sequences of positive numbers, $\mathcal{G}^{\mathcal{R}}$ is defined by the
sections $u\in\mathcal{G}\left(  \Omega\right)  $ having a representative
verifying%
\[
\forall K\Subset\Omega\text{ }\exists\left(  N_{l}\right)  _{l\geq0}%
\in\mathcal{R}\text{ }\underset{x\in K}{\text{sup}}\left\vert \partial
^{\alpha}u_{\varepsilon}\left(  x\right)  \right\vert =\mathcal{O}\left(
\varepsilon^{-N_{\left\vert \alpha\right\vert }}\right)  \text{ as
}\varepsilon\rightarrow0.
\]

Under certain stability conditions on the set $\mathcal{R}$ exposed in
\cite{ADRReg}, $\mathcal{G}^{\mathcal{R}}$ is a subsheaf of differential
algebras of $\mathcal{G}$ and when $\mathcal{R}$ consists of the set of all
bounded sequences, then $\mathcal{G}^{\mathcal{R}}=\mathcal{G}^{\infty}$.

\noindent Then%
\[
\mathcal{S}_{\mathcal{G}}^{\mathcal{G}^{\mathcal{R}}}\left(  u\right)
=sing\,supp^{\mathcal{R}}\left(  u\right)
\]
where $sing\,supp^{\mathcal{R}}\,u$ is the set of all $x\in\Omega$ having no
neighbourhood $V$ such that $u\!\mid_{V}\in\mathcal{G}^{\mathcal{R}}\left(
V\right)  $.
\end{example}

\begin{example}
$:\mathbf{Local}$ $\mathbf{analysis}$ $\mathbf{of}$\textbf{ }$\mathcal{G}^{L}%
$-$\mathbf{type}\medskip$

In \cite{JAM3} one constructs $\mathcal{B}=\mathcal{G}^{L}$ as a special
regular sub(pre)sheaf of $\mathcal{A}=$ $\mathcal{G}$, extending in a
generalized sense the $\mathrm{C}^{L}$ classes of H\"{o}rmander
\cite{HorPDOT1} containing analytic and Gevrey classes and constructed from an
increasing sequence $L_{k}$ of positive numbers such that $L_{0}$=$1$ and%
\[
k\leq L_{k},\;L_{k+1}\leq CL_{k}%
\]
for some constant $C$. When taking $L_{k}=k+1$, we obtain the analytic case
$\mathcal{G}^{A}$ studied in \cite{PSV} involving special properties of
holomorphic generalized functions which give to $\mathcal{G}^{A}$ a sheaf property.
\end{example}

\begin{example}
$:\mathcal{G}^{\sigma,\infty}$-$\mathbf{local}$ $\mathbf{analysis}$
$\mathbf{in}$\textbf{ \medskip}$\mathcal{G}^{\sigma}\left(  \Omega\right)  $

In \cite{BouBen} Bouzar and Benmeriem introduce a sheaf of algebra
$\mathcal{G}^{\sigma}\neq\mathcal{G}$ of Gevrey ultradistributions with
another asymptotic scale than the Colombeau one by replacing the estimate
$O\left(  \varepsilon^{-m}\right)  $ by $O\left(  e^{\varepsilon^{-\frac
{1}{2\sigma-1}}}\right)  ^{m}$ (resp. $O\left(  \varepsilon^{p}\right)  $ by
$O\left(  e^{-\varepsilon^{-\frac{1}{2\sigma-1}}}\right)  ^{p}$)in the
definition of moderate (resp. null) elements When taking $L_{k}=\left(
k+1\right)  ^{\sigma}$ they can construct a subpresheaf $\mathcal{G}%
^{\sigma,\infty}$ of $\mathcal{G}^{\sigma}$and give a study of $\mathcal{G}%
^{s,\infty}$-singularity.
\end{example}

By choosing $\mathcal{R}$ as a regular subset of $\mathbb{R}_{+}^{\mathbb{N}}%
$, Delcroix has extended the $\mathcal{G}^{\infty}$-regularity into the
$\mathcal{G}^{\mathcal{R}}$ one, and in a work in progress Bouzar replaces the
classical regularity by the $\mathcal{R}$-regularity to extend the regular
generalized Gevrey ultradistributions. We follow this way in view of
constructing a general model containing all the previous examples but we have
to add two other parameters: an asymptotic scale $r=\left(  r_{\lambda
}\right)  _{\lambda}\in\left(  \mathbb{R}_{+}^{\ast}\right)  ^{\Lambda}$ and a
sequence $L_{k}$ of positive numbers such that $L_{0}$=$1$ and $k\leq
L_{k},\;L_{k+1}\leq CL_{k}$ for some constant $C$.

\subsection{ $\mathcal{G}^{r,\mathcal{R}}$ and $\mathcal{G}^{r,\mathcal{R}%
\text{,}L}$-local analysis in $\mathcal{G}^{r}\left(  \Omega\right)  $}

\subsubsection{The $\mathcal{G}^{r}$ sheaf of algebras}

\bigskip Let us consider

$\bullet$ $\mathcal{E}=\mathrm{C}^{\infty}$ as starting sheaf of algebras, for
each open set $\Omega$ in $\mathbb{R}^{n\text{ }}$, $\mathrm{C}^{\infty
}\left(  \Omega\right)  $ is endowed by the usual family $p_{K,\alpha}$ of seminorms

$\bullet$ $\Lambda$ a set of indices left-filtering for the given (partial)
order relation $\prec$.

$\bullet$ an asymptotic scale $r=\left(  r_{\lambda}\right)  _{\lambda}%
\in\left(  \mathbb{R}_{+}^{\ast}\right)  ^{\Lambda}$ such that $\lim
\limits_{\Lambda}$ $r_{\lambda}=0$, (or $r_{\lambda}\rightarrow0$), that is to
say: for each $\mathbb{R}$\textit{-}neighbourhood $W$\ of $0$, there exists
$\lambda_{0}\in\Lambda$\ such that%
\[
\lambda\prec\lambda_{0}\Longrightarrow r_{\lambda}\in W.
\]

Define the functors $\mathcal{X}^{r}$ (resp. $\mathcal{N}^{r}$)$:\Omega
\mapsto\mathcal{X}^{r}\left(  \Omega\right)  $ (resp. $\mathcal{N}^{r}\left(
\Omega\right)  $) by%
\begin{align*}
\mathcal{X}^{r}\left(  \Omega\right)   &  =\ \left\{  \left(  u_{\lambda
}\right)  _{\lambda}\in\left[  \mathrm{C}^{\infty}\left(  \Omega\right)
\right]  ^{\Lambda},\forall K\Subset\Omega,\forall\alpha\in\mathbb{N}%
^{n},\exists N\in\mathbb{N},p_{K,\alpha}\left(  u_{\lambda}\right)  =O\left(
r_{\lambda}^{-N}\right)  \text{for }r_{\lambda}\rightarrow0\right\} \\
\mathcal{N}^{r}\left(  \Omega\right)   &  =\ \left\{  \left(  u_{\lambda
}\right)  _{\lambda}\in\left[  \mathrm{C}^{\infty}\left(  \Omega\right)
\right]  ^{\Lambda},\forall K\Subset\Omega,\forall\alpha\in\mathbb{N}%
^{n},\forall m\in\mathbb{N},p_{K,\alpha}\left(  u_{\lambda}\right)  =O\left(
r_{\lambda}^{m}\right)  \text{for }r_{\lambda}\rightarrow0\right\}
\end{align*}
it is not difficult to prove with the same techniques as Colombeau ones
\cite{JFC3} that $\mathcal{X}^{r}$ and $\mathcal{N}^{r}$ are respectively a
sheaf of differential algebras and a subsheaf of ideals of $\mathcal{X}^{r}%
$over the ring%
\[
\mathcal{X}^{r}\left(  \mathbb{C}\right)  =\left\{  \left(  s_{\lambda
}\right)  _{\lambda}\in\mathbb{C}^{\Lambda},\exists N\in\mathbb{N},\left\vert
s_{\lambda}\right\vert =O\left(  r_{\lambda}^{-N}\right)  \text{for
}r_{\lambda}\rightarrow0\right\}  \text{.}%
\]
Then $\mathcal{X}^{r}/\mathcal{N}^{r}=\mathcal{G}^{r}$ is a priori a factor
presheaf of Colombeau type. It is well known that $\mathcal{G}$ is a sheaf and
even a fine sheaf. The first assumption (a result from Aragona and Biagioni
\cite{AraBia}) is based on the existence of a \textrm{C}$^{\infty}$ partition
of unity associated to any open covering of $\Omega$ (due to the fact that
$\mathbb{R}^{d}$ is a locally compact Hausdorff space). On the other hand, we
can notice that \textrm{C}$^{\infty}$ is a fine sheaf because multiplication
by a smooth function defines a sheaf homomorphism in a natural way. Hence the
usual topology and \textrm{C}$^{\infty}$ partition of unity define the
required sheaf partition of unity according to the definition in sheaf theory.
This leads very easily to the second assumption, from the well known result
that any sheaf of modules on a fine sheaf is itself a fine sheaf: it is
precisely the case of $\mathcal{G}$ which is a sheaf of \textrm{C}$^{\infty}$
modules. And it is the same for $\mathcal{G}^{r}$ which is a fine sheaf of
\textrm{C}$^{\infty}$ modules and also a sheaf of differential algebras over
the ring $\mathcal{X}^{r}\left(  \mathbb{C}\right)  /\mathcal{N}^{r}\left(
\mathbb{C}\right)  $ with%
\[
\mathcal{N}^{r}\left(  \mathbb{C}\right)  =\left\{  \left(  s_{\lambda
}\right)  _{\lambda}\in\mathbb{C}^{\Lambda},\forall m\in\mathbb{N},\left\vert
s_{\lambda}\right\vert =O\left(  r_{\lambda}^{m}\right)  \text{for }%
r_{\lambda}\rightarrow0\right\}  \text{.}%
\]

\subsubsection{The $\mathcal{G}^{r,\mathcal{R}}$ subsheaf of $\mathcal{G}^{r}%
$}

The $\mathcal{G}^{r,\mathcal{R}}$ regularity of $\mathcal{G}^{r}$ generalizes
the $\mathcal{G}^{\infty}$ regularity of $\mathcal{G}$. We begin by defining a
regular subspace $\mathcal{R}$ of $\mathbb{R}_{+}^{n}$ in the Delcroix sense
\cite{ADRReg}:

\begin{definition}
\label{Regular}A subspace $\mathcal{R}$ of $\mathbb{R}_{+}^{\mathbb{N}}$ is
\emph{regular} if $\mathcal{R}$ is non empty and\newline

$\left(  i\right)  $~$\mathcal{R}$ is \textquotedblleft
overstable\textquotedblright\ by translation and by maximum%
\[
\forall N\in\mathcal{R},\ \forall\left(  k,k^{\prime}\right)  \in
\mathbb{N}^{2},\ \exists N^{\prime}\in\mathcal{R},\ \ \forall n\in
\mathbb{N},\ \ \ N\left(  n+k\right)  +k^{\prime}\leq N^{\prime}\left(
n\right)  ,\vspace{-0.04in}%
\]%
\[
\forall N_{1}\in\mathcal{R},\ \forall N_{2}\in\mathcal{R},\ \exists
N\in\mathcal{R},\ \forall n\in\mathbb{N},\ \ \ \max\left(  N_{1}\left(
n\right)  ,N_{2}\left(  n\right)  \right)  \leq N\left(  n\right)  .
\]

$\left(  ii\right)  $~For all $N_{1}\ $and $N_{2}$ in $\mathcal{R}$, there
exists $N\in\mathcal{R}$ such that%
\[
\forall\left(  l_{1},l_{2}\right)  \in\mathbb{N}^{2},\ \ \ N_{1}\left(
l_{1}\right)  +N_{2}\left(  l_{2}\right)  \leq N\left(  l_{1}+l_{2}\right)  .
\]

\end{definition}

Then, for any regular subset $\mathcal{R}$ of $\mathbb{R}_{+}^{\mathbb{N}}$ we
can set%
\begin{align*}
\mathcal{X}^{r,\mathcal{R}}\left(  \Omega\right)   &  =\ \left\{  \left(
u_{\lambda}\right)  _{\lambda}\in\left[  \mathrm{C}^{\infty}\left(
\Omega\right)  \right]  ^{\Lambda},\forall K\Subset\Omega,\exists
N\in\mathcal{R},\forall\alpha\in\mathbb{N}^{n},p_{K,\alpha}\left(  u_{\lambda
}\right)  =O\left(  r_{\lambda}^{-N\left(  \left\vert \alpha\right\vert
\right)  }\right)  \text{for }r_{\lambda}\rightarrow0\right\} \\
\mathcal{N}^{r,\mathcal{R}}\left(  \Omega\right)   &  =\ \left\{  \left(
u_{\lambda}\right)  _{\lambda}\in\left[  \mathrm{C}^{\infty}\left(
\Omega\right)  \right]  ^{\Lambda},\forall K\Subset\Omega,\forall
m\in\mathcal{R},\forall\alpha\in\mathbb{N}^{n},p_{K,\alpha}\left(  u_{\lambda
}\right)  =O\left(  r_{\lambda}^{m\left(  \left\vert \alpha\right\vert
\right)  }\right)  \text{for }r_{\lambda}\rightarrow0\right\}  \text{.}%
\end{align*}

\begin{proposition}
$\medskip$

$\left(  i\right)  $~For any regular subspace $\mathcal{R}$ of $\mathbb{R}%
_{+}^{\mathbb{N}}$, the functor $\Omega\rightarrow\mathcal{X}^{r,\mathcal{R}%
}\left(  \Omega\right)  $ defines a sheaf of differential algebras over the
ring $\mathcal{X}^{r}\left(  \mathbb{C}\right)  $.

$\left(  ii\right)  $~The set $\mathcal{N}^{r,\mathcal{R}}\left(
\Omega\right)  $ is equal to $\mathcal{N}^{r}\left(  \Omega\right)  $. Thus,
the functor $\mathcal{N}^{\mathcal{R}}:\Omega\rightarrow\mathcal{N}%
^{\mathcal{R}}\left(  \Omega\right)  $ defines a sheaf of ideals of the sheaf
$\mathcal{X}^{\mathcal{R}}\left(  \cdot\right)  $.

$\left(  iii\right)  $~For any regular subspaces $\mathcal{R}_{1}$ and
$\mathcal{R}_{2}$ of $\mathbb{R}_{+}^{\mathbb{N}}$, with $\mathcal{R}%
_{1}\subset\mathcal{R}_{2}$, the sheaf $\mathcal{X}^{\mathcal{R}_{1}}\left(
\Omega\right)  $ is a subsheaf of the sheaf $\mathcal{X}^{\mathcal{R}_{2}%
}\left(  \Omega\right)  $.
\end{proposition}

\begin{proof}
The proof follows the same lines as in the case of $\mathcal{G}^{\mathcal{R}}$
algebras (see \cite{ADRReg}, Proposition 1.). We have to verify that our
asymptotic scale $\left(  r_{\lambda}\right)  _{\lambda}$ involving a more
general parametrization doesn't modify the results. We deduce assertion
$\left(  i\right)  $ from the assertion $\left(  i\right)  $ in the definition
of $\mathcal{R}$. For the equality $\mathcal{N}^{r,\mathcal{R}}\left(
\Omega\right)  =\mathcal{N}^{r}\left(  \Omega\right)  $, take first $\left(
u_{\lambda}\right)  _{\lambda}\in\mathcal{N}^{r,\mathcal{R}}\left(
\Omega\right)  $. For any $K\Subset\Omega$, $\alpha\in\mathbb{N}^{n}$ and
$m\in\mathbb{N}$, choose $N\in\mathcal{R}$. From $\left(  i\right)  $ in
definition \ref{Regular} there exists $N^{\prime}\in\mathcal{R}$ such that
$N+m\leq N^{\prime}$. Thus $p_{K,\alpha}\left(  u_{\lambda}\right)  =O\left(
r_{\lambda}^{N^{\prime}\left(  \left\vert \alpha\right\vert \right)  }\right)
=O\left(  r_{\lambda}^{m}\right)  $ and $\left(  u_{\lambda}\right)
_{\lambda}\in\mathcal{N}^{r}\left(  \Omega\right)  $. Conversely for given
$\left(  u_{\lambda}\right)  _{\lambda}\in\mathcal{N}^{r}\left(
\Omega\right)  $ and $N\in\mathcal{R}$ we have $p_{K,\alpha}\left(
u_{\lambda}\right)  =O\left(  r_{\lambda}^{N\left(  \left\vert \alpha
\right\vert \right)  }\right)  $ since this estimate holds for all
$m\in\mathbb{N}$. For the sheaf properties we have to replace Colombeau's
estimates by $\mathcal{X}^{r,\mathcal{R}}$ estimates and consider only a
finite number of terms by compactness. Thus, from $\left(  ii\right)  $ in
definition \ref{Regular} we have the results. The inclusion $\mathcal{X}%
^{\mathcal{R}_{1}}\left(  \Omega\right)  \subset\mathcal{X}^{\mathcal{R}_{2}%
}\left(  \Omega\right)  $ prove $\left(  iii\right)  $.
\end{proof}

According to same arguments as those used for $\mathcal{G}^{r}$ the presheaf
$\mathcal{G}^{r,\mathcal{R}}\,=\mathcal{X}^{r,\mathcal{R}}/\mathcal{N}%
^{r,\mathcal{R}}=\mathcal{X}^{r,\mathcal{R}}/\mathcal{N}^{r}$ turns to be a
sheaf of differentiable algebras on the ring $\mathcal{X}^{r,\mathcal{R}%
}\left(  \mathbb{C}\right)  /\mathcal{N}^{r}\left(  \mathbb{C}\right)  $ with
\[
\mathcal{N}^{r}\left(  \mathbb{K}\right)  =\left\{  \left(  s_{\lambda
}\right)  _{\lambda}\in\mathbb{C}^{\Lambda},\forall m\in\mathbb{N},\left\vert
s_{\lambda}\right\vert =O\left(  r_{\lambda}^{m}\right)  \text{ for
}r_{\lambda}\rightarrow0\right\}  ,\ \ \mathbb{K}=\mathbb{R}\text{ or
}\mathbb{K}=\mathbb{C}\text{.}%
\]

Moreover, from $\left(  ii\right)  $ in the above proposition, $\mathcal{G}%
^{r,\mathcal{R}}$ is a subsheaf of $\mathcal{G}^{r}$.

\begin{definition}
For any regular subset $\mathcal{R}$ of $\mathbb{R}_{+}^{\mathbb{N}}$, the
sheaf of algebras (subsheaf of $\mathcal{G}^{r}$)%
\[
\mathcal{G}^{r,\mathcal{R}}\,=\mathcal{X}^{r,\mathcal{R}}\,/\mathcal{N}%
^{r,\mathcal{R}}\,
\]
is called the sheaf of $\left(  r,\mathcal{R}\right)  $-regular algebras of
(nonlinear) generalized functions.
\end{definition}

\begin{example}
Taking $\lambda=\varepsilon\in\left]  0,1\right]  $, $r_{\varepsilon
}=\varepsilon$ and $\mathcal{R}=\mathbb{R}_{+}^{\mathbb{N}}$, we recover the
sheaf $\mathcal{G}$ of Colombeau simplified algebras.

Taking $\lambda=\varepsilon\in\left]  0,1\right]  $, $r_{\varepsilon
}=\varepsilon$ and $\mathcal{R}=\mathcal{B}o$ (the set of bounded sequences),
we obtain the sheaf of $\mathcal{G}^{\infty}$-generalized functions
\cite{ADRReg}.

Taking $\lambda=\varepsilon\in\left]  0,1\right]  $, $r_{\varepsilon
}=e^{\varepsilon^{-\frac{1}{2\sigma-1}}}$ and $\mathcal{R}=\mathbb{R}%
_{+}^{\mathbb{N}}$, we obtain the sheaf of so called $\mathcal{G}^{s}%
$-generalized functions in \cite{BouBen}.
\end{example}

\subsubsection{The $\mathcal{G}^{r,\mathcal{R},L}$ subpresheaf of
$\mathcal{G}^{r}$ and $\left(  r\text{,}\mathcal{R}\text{,}L\right)
$-analysis}

Let $L_{k}$ be an increasing sequence of positive numbers such that $L_{0}$=1
and%
\[
k\leq L_{k},\;L_{k+1}\leq CL_{k}%
\]
for some constant $C.$ According to H\"{o}rmander definition given in
subsection 8.4 of \cite{HorPDOT1}, we shall denote by $\mathrm{C}^{L}$ the
sheaf of $\mathbb{K}$-algebras on $\mathbb{R}^{n}$ ($\mathbb{K=R}$ or
$\mathbb{C}$) such that, for any open set $\Omega\subset\mathbb{R}^{n}$%

\[
\mathrm{C}^{L}\left(  \Omega\right)  =\left\{  u=\mathrm{C}^{\infty}\left(
\Omega\right)  \mid\forall K\Subset\Omega,\exists c>0,\forall\alpha
\in\mathbb{N}^{n},\underset{x\in K}{\sup}\left\vert \mathcal{D}^{\alpha
}u\left(  x\right)  \right\vert \leq c\left(  cL_{\left\vert \alpha\right\vert
}\right)  ^{\left\vert \alpha\right\vert }\right\}
\]

When $L_{k}=k+1$, $\mathrm{C}^{L}$ is the sheaf $\mathrm{A}$ of analytical
functions. If $L_{k}=\left(  k+1\right)  ^{a}$, $a>1$, $\mathrm{C}^{L}$ is the
sheaf $\mathrm{G}_{a}$ of the Gevrey class of order $a$.

It is possible to enlarge the above definition into a generalized one
involving three parameters $r$, $\mathcal{R}$, $L$ corresponding to a choice
of some asymptotic scale $r=\left(  r_{\lambda}\right)  _{\lambda}\in\left(
\mathbb{R}_{+}^{\ast}\right)  ^{\Lambda}$, a regular subset $\mathcal{R}$ of
$\mathbb{R}_{+}^{\mathbb{N}}$, and a sequence $L=\left(  L_{k}\right)  _{k}$.

\begin{definition}
Let us define the functors $\mathcal{X}^{r,\mathcal{R},L}$ (resp.
$\mathcal{N}^{r,\mathcal{R},L}$)$:\Omega\mapsto\mathcal{X}^{r,\mathcal{R}%
,L}\left(  \Omega\right)  $ (resp. $\mathcal{N}^{r,\mathcal{R},L}\left(
\Omega\right)  $) by

$\left\{
\begin{array}
[c]{l}%
\mathcal{X}^{r,\mathcal{R},L}\left(  \Omega\right)  =\left\{  \left(
u_{\lambda}\right)  _{\lambda}\in\left[  \mathrm{C}^{\infty}\left(
\Omega\right)  \right]  ^{\Lambda},\forall K\Subset\Omega,\exists
N\in\mathcal{R},\exists c>0,\exists\lambda_{0}\in\Lambda\right. \\
\multicolumn{1}{r}{\;\;\;\;\;\;\;\;\;\;\;\;\;\;\;\;\;\;\;\;\;\;\;\;\;\;\ \;\;\;\;\;\left.
\forall\alpha\in\mathbb{N}^{n},\forall\lambda\prec\lambda_{0}:\sup
\limits_{x\in K}\left\vert D^{\alpha}u_{\lambda}\left(  x\right)  \right\vert
\leq cr_{\lambda}^{-N\left(  \left\vert \alpha\right\vert \right)  }\left(
cL_{\left\vert \alpha\right\vert }\right)  ^{\left\vert \alpha\right\vert
}\right\}  \text{,}}\\
\mathcal{N}^{r,\mathcal{R},L}\left(  \Omega\right)  =\mathcal{X}%
^{r,\mathcal{R},L}\left(  \Omega\right)  \cap\mathcal{N}^{r}\left(
\Omega\right)  \text{.}%
\end{array}
\right.  $
\end{definition}

\begin{lemma}
$\mathcal{X}^{r,\mathcal{R},L}$ is a subsheaf of $\mathcal{X}^{r,\mathcal{R}}%
$, $\mathcal{N}^{r,\mathcal{R},L}$ is a sheaf of ideals of $\mathcal{X}%
^{r,\mathcal{R},L}$
\end{lemma}

\begin{proof}
For each $\Omega$, $\mathcal{X}^{r,\mathcal{R},L}\left(  \Omega\right)  $ is a
subalgebra of $\mathcal{X}^{r,\mathcal{R}}\left(  \Omega\right)  $, and the
restriction and localization processes are obvious. Let us try to glue
together the bits, giving some family $\left(  \Omega_{i}\right)  _{i\in I}$,
with $\Omega=\underset{_{1\leq i\leq L}}{\cup}\Omega_{i}$, and $U_{i}=\left(
u_{i,\lambda}\right)  _{\varepsilon}\in\mathcal{X}^{r,\mathcal{R},L}\left(
\Omega_{i}\right)  $ with $U_{i}=U_{j}$ on $\Omega_{i}\cap\Omega j$. We begin
to define $U\left(  x\right)  $ as $U_{i}\left(  x\right)  $ when $x\in\Omega$
lies in $\Omega_{i}$. Clearly $U$ belongs to $\left[  \mathrm{C}^{\infty
}\left(  \Omega\right)  \right]  ^{\Lambda}$. Let $K\Subset\Omega$. We can
cover $K$ by a finite number of $\Omega_{i}:\Omega_{1}$,...$\Omega_{p}$ such
that $K=\underset{_{1\leq i\leq L}}{\cup}K_{i}$, with $K_{i}=K\cap
\overline{\Omega_{i}^{\prime}}\subset\Omega_{i}$. This\ is possible by
choosing $\Omega_{i}^{\prime}\subsetneqq\Omega_{i}$,$d\left(  \Omega
_{i}^{\prime},\Omega_{i}\right)  \leq d(K,\Omega)$ for $1\leq i\leq p$. Then%
\[
\exists N_{i}\in\mathcal{R},\exists c_{i}>0,\exists\lambda_{0,i}\in
\Lambda,\forall\alpha\in\mathbb{N}^{n},\forall\lambda\prec\lambda_{0,i}%
:\sup\limits_{x\in K_{i}}\left\vert D^{\alpha}u_{i,\lambda}\left(  x\right)
\right\vert \leq c_{i}r_{\lambda}^{-N_{i}\left(  \left\vert \alpha\right\vert
\right)  }\left(  c_{i}L_{\left\vert \alpha\right\vert }\right)  ^{\left\vert
\alpha\right\vert }%
\]

From the assumption there exits some $\lambda_{0}\in\Lambda$ such that
$\lambda_{0}\prec\lambda_{0,i}$ and $N\in\mathcal{R}$ such that $N\geq N_{i}$
for $1\leq i\leq p$. Set $c=\underset{_{1\leq i\leq L}}{\max}c_{i}$. Then for
each $\alpha\in\mathbb{N}^{n}$ and $\lambda\prec\lambda_{0}$ we have%
\[
\sup\limits_{x\in K}\left\vert D^{\alpha}u_{\lambda}\left(  x\right)
\right\vert =\underset{1\leq i\leq p}{%
{\textstyle\sum}
}cr_{\lambda}^{-N_{i}\left(  \left\vert \alpha\right\vert \right)  }\left(
cL_{\left\vert \alpha\right\vert }\right)  ^{\left\vert \alpha\right\vert
}\leq pcr_{\lambda}^{-N\left(  \left\vert \alpha\right\vert \right)  }\left(
cL_{\left\vert \alpha\right\vert }\right)  ^{\left\vert \alpha\right\vert
}\leq c^{\prime}r_{\lambda}^{-N\left(  \left\vert \alpha\right\vert \right)
}\left(  c^{\prime}L_{\left\vert \alpha\right\vert }\right)  ^{\left\vert
\alpha\right\vert }.
\]

Thus $U$ belongs to $\mathcal{X}^{L}\left(  \Omega\right)  $. It is easy to
see that for each $\Omega$, $\mathcal{N}_{\ast}^{L}\left(  \Omega\right)  $ is
an ideal of $\mathcal{X}^{L}\left(  \Omega\right)  $, and the same proof as
above leads to the sheaf structure of $\mathcal{N}_{\ast}^{L}$. Then we can
define a new factor presheaf of $C^{L}$-type algebras which is a subpresheaf
of $\mathcal{G}^{r,\mathcal{R}}$ and $\mathcal{G}^{r}$, according to the
definition of the ideal $\mathcal{N}^{r,\mathcal{R},L}\left(  \Omega\right)
=\mathcal{X}^{r,\mathcal{R},L}\left(  \Omega\right)  \cap\mathcal{N}%
^{r}\left(  \Omega\right)  $.
\end{proof}

\begin{definition}
For any asymptotic scale $r=\left(  r_{\lambda}\right)  _{\lambda}\in\left(
\mathbb{R}_{+}^{\ast}\right)  ^{\Lambda}$, any regular subset $\mathcal{R}$ of
$\mathbb{R}_{+}^{\mathbb{N}}$, and any sequence $L=\left(  L_{k}\right)  _{k}%
$, the presheaf of algebras (subpresheaf of $\mathcal{G}^{r}$)%
\[
\mathcal{G}^{r,\mathcal{R},L}\,=\mathcal{X}^{r,\mathcal{R},L}\,/\mathcal{N}%
^{r,\mathcal{R},L}\,
\]
is called the sheaf of $\left(  r,\mathcal{R},L\right)  $-regular algebras of
(nonlinear) generalized functions.
\end{definition}

\begin{example}
Taking $\lambda=\varepsilon\in\left]  0,1\right]  $, $r_{\varepsilon
}=\varepsilon$ , $\mathcal{R}=\mathcal{B}o$ (the set of bounded sequences),
and some $L=\left(  L_{k}\right)  _{k}$, we obtain the presheaf of
$\mathcal{G}^{L}$-generalized functions (subpresheaf of $\mathcal{G}$)
\cite{JAM3}.

Taking $\lambda=\varepsilon\in\left]  0,1\right]  $, $r_{\varepsilon
}=\varepsilon$ , $\mathcal{R}=\mathcal{B}o$ and $L_{k}=k+1$, we obtain the
sheaf of $\mathcal{G}^{A}$-generalized functions (subsheaf of $\mathcal{G}$)
\cite{PSV}.

Taking $\lambda=\varepsilon\in\left]  0,1\right]  $, $r_{\varepsilon
}=e^{\varepsilon^{-\frac{1}{2\sigma-1}}}$ , $\mathcal{R}=\mathcal{B}o$ and
$L_{k}=\left(  k+1\right)  ^{\sigma}$, we obtain the presheaf of
$\mathcal{G}^{\sigma,\infty}$-generalized functions (subpresheaf of
$\mathcal{G}^{\sigma}$) \cite{BouBen}.
\end{example}

\begin{remark}
We proved that $\mathcal{G}^{r,\mathcal{R},L}$is a presheaf ; the localization
principle $(F_{1})$ is not difficult to prove. However, we lack some suitable
partition of unity or other argument which preserve the $\mathcal{X}%
^{r,\mathcal{R},L}$ estimates and permit to glue together the bits and get
$(F_{1})$. But instead of introducing the sheaf associated to $\mathcal{G}%
^{r,\mathcal{R},L}$, one can keep its presheaf structure in the following
localization processes.

Taking $\mathcal{A}=\mathcal{G}^{r}$ and $\mathcal{B}=\mathcal{G}%
^{r,\mathcal{R},L}$, the $\mathcal{G}^{r,\mathcal{R},L}$-singularities of
$u\in\mathcal{G}^{r}(\Omega)$ are localized in%
\[
\mathcal{S}_{\mathcal{G}^{r}}^{\mathcal{G}^{r,\mathcal{R},L}}\left(  u\right)
=sing\,supp^{^{(r,\mathcal{R},L)}}\left(  u\right)
\]
where $sing\,supp^{^{(r,\mathcal{R},L)}}(u)$ is the set of all $x\in\Omega$
having no neighbourhood $V$ such that $u\left\vert _{V}\right.  \in
\mathcal{G}^{r,\mathcal{R},L}\left(  V\right)  $. But we cannot prove that
$\Omega\backslash sing\,supp^{^{(r,\mathcal{R},L)}}\left(  u\right)  $ is the
largest open set $\mathcal{O}$ such that $u\in\mathcal{G}^{r,\mathcal{R}%
,L}\left(  \mathcal{O}\right)  $. Due to the lack of the $\left(
F_{2}\right)  $ principe in lemma 2, we don't know how to prove the existence
of such an open set.
\end{remark}

\begin{remark}
Following the previous definition 2 of $\mathcal{B}$-local regularity of an
element $u\in\mathcal{A}\left(  \Omega\right)  $, we can naturally set that
$u\in\mathcal{G}^{r}(\Omega)$ is locally $\mathcal{G}^{r,\mathcal{R},L}$ at
$x\in\Omega$ if for some neighbourhood $V$ of $x\in\Omega$ the restriction
$u\!\mid_{V}$ belongs to $\mathcal{G}^{r,\mathcal{R},L}\left(  V\right)  $,
that is to say if for some neighbourhood $V$ of $x$ there exists a
representative $\left(  u_{\lambda}\right)  _{\lambda}$ of $u$ such that
$\left(  u_{\lambda}\!\mid_{V}\right)  _{\lambda}$ belongs to $\mathcal{X}%
^{r,\mathcal{R},L}\left(  V\right)  $. Let $\mathcal{H}^{r,\mathcal{R}%
,L}\left(  \Omega\right)  $ be the set of all $u\in\mathcal{G}^{r}(\Omega)$
which are locally $\mathcal{G}^{r,\mathcal{R},L}$ at $x\in\Omega$. It is not
difficult to prove that the sheaf associated to $\mathcal{G}^{r,\mathcal{R}%
,L}$ is the functor $\Omega\mapsto\mathcal{H}^{r,\mathcal{R},L}\left(
\Omega\right)  $. It is a subsheaf of algebras of $\mathcal{G}^{r}$. But in
the general case, we cannot prove that $u\in\mathcal{H}^{r,\mathcal{R}%
,L}\left(  \Omega\right)  $ has a global representative in $\mathcal{X}%
^{r,\mathcal{R},L}\left(  \Omega\right)  $. However this is fulfilled when
taking $r_{\varepsilon}=\varepsilon$, $\mathcal{R=B}o$, $L_{k}=k+1$,
corresponding to the analytic case studied in \cite{PSV} involving special
properties of holomorphic generalized functions. And then $\mathcal{H}%
^{r,\mathcal{R},L}=\mathcal{G}^{A}$, the subsheaf of generalized analytic
functions of the sheaf $\mathcal{G}$. But when $\mathcal{A}=\mathcal{G}$ and
$\mathcal{B}=\mathcal{G}^{A}$, the $\mathcal{G}^{A}$-singularities of
$u\in\mathcal{G}(\Omega)$ are always localized in%
\[
\mathcal{S}_{\mathcal{G}}^{\mathcal{G}^{A}}\left(  u\right)  =sing\,supp^{A}%
\left(  u\right)
\]
where $sing\,supp^{A}\,u$ is the set of all $x\in\Omega$ having no
neighbourhood $V$ such that $u\left\vert _{V}\right.  \in\mathcal{G}%
^{A}\left(  V\right)  $. Here, the sheaf structure of $\mathcal{G}^{A}$
provides the following precision: from lemma 2 involving the $\left(
F_{2}\right)  $ principe, we can prove that $\Omega\backslash sing\,supp^{A}%
$\thinspace$u$ is also the largest open set $\mathcal{O}$ such that
$u\in\mathcal{G}^{A}\left(  \mathcal{O}\right)  $.
\end{remark}

\subsubsection{Canonical embeddings}

\begin{lemma}
We have the following commutative diagram in which the arrows are canonical
embeddings%
\[%
\begin{array}
[c]{ccc}%
\mathrm{C}^{L} & \rightarrow & \mathrm{C}^{\infty}\\
\downarrow &  & \downarrow\\
\mathcal{G}^{r,\mathcal{R},L} & \rightarrow & \mathcal{G}^{r,\mathcal{R}}%
\end{array}
\]

\end{lemma}

\begin{proof}
The canonical (pre)sheaf embedding of $\mathrm{C}^{L}$ into $\mathcal{G}%
^{r,\mathcal{R},L}$ (resp. of $\mathrm{C}^{\infty}$ into $\mathcal{G}%
^{r,\mathcal{R}}$) is defined for each open set $\Omega\subset\mathbb{R}^{n}$
in by the canonical map%
\[
\mathrm{C}^{L}\left(  \Omega\right)  \rightarrow\mathcal{G}^{r,\mathcal{R}%
,L}\text{ }\left(  \Omega\right)  \text{\ (resp. }\mathrm{C}^{\infty
}\rightarrow\mathcal{G}^{r,\mathcal{R}}\text{)}:\text{\ }u\mapsto\left[
u_{\lambda}\right]  \text{, with }u_{\lambda}=u\text{ for }\lambda\in\Lambda
\]
which is \ an injective homomorphism of algebras, $\left[  u_{\lambda}\right]
$ being the class of $u$ $\in\mathrm{C}^{L}\left(  \Omega\right)  $ (resp.
$\mathrm{C}^{\infty}\left(  \Omega\right)  $) in the factor algebra
$\mathcal{G}^{r,\mathcal{R},L}$ $\left(  \Omega\right)  $ (resp.$\mathcal{G}%
^{r,\mathcal{R}}\left(  \Omega\right)  $). In order to construct the sheaf
embedding : $\mathcal{G}^{r,\mathcal{R},L}\rightarrow\mathcal{G}%
^{^{r,\mathcal{R}}}$, we recall that
\[
\mathcal{X}^{r,\mathcal{R},L}\left(  \Omega\right)  \subset\mathcal{X}%
^{r,\mathcal{R}}\left(  \Omega\right)  ;\mathcal{N}^{r,\mathcal{R},L}\left(
\Omega\right)  =\mathcal{X}^{r,\mathcal{R},L}\left(  \Omega\right)
\cap\mathcal{N}^{r}\left(  \Omega\right)
\]
is a necessary and sufficient condition to embed $\mathcal{G}^{r,\mathcal{R}%
,L}$ $\left(  \Omega\right)  $ into $\mathcal{G}^{^{r,\mathcal{R}}}\left(
\Omega\right)  $, from which we deduce the required sheaf embedding.
\end{proof}

\begin{remark}
When $\lambda=\varepsilon$ and $r_{\varepsilon}=\varepsilon$, we can suppress
the symbol $r$ in the previous formulation. When taking $\mathcal{R}%
=\mathcal{B}o$, the symbol $\mathcal{R}$ becomes $\infty$. For example, when
we do that simultaneously, we have%
\[
\mathcal{G}^{r}=\mathcal{G}\text{ };\mathcal{G}^{r,\mathcal{R}}=\mathcal{G}%
^{\infty}\text{ };\mathcal{G}^{r,\mathcal{R},L}=\mathcal{G}^{L}.
\]

\end{remark}

\subsection{$\mathcal{G}$-local analysis in\textbf{ }$\mathcal{L}\left(
\mathcal{G}_{c}(\Omega),\widetilde{\mathbb{C}}\right)  $}

\subsubsection{Duality in the Colombeau context}

Starting from the usual family of semi norms $\left(  p_{i}\right)  _{i\in I}$
defining the topology of \textrm{C}$^{\infty}\left(  \Omega\right)  $ by
\[
p_{i}\left(  f\right)  =p_{K,l}\left(  f\right)  =\underset{x\in K,\left\vert
\alpha\right\vert \leq l}{\sup}\left\vert \partial^{\alpha}f\left(  x\right)
\right\vert
\]
the so-called sharp topology of $\mathcal{G}(\Omega)$ \cite{NePiSc}is defined
by the family of ultra-pseudo-seminorms $\left(  \mathcal{P}_{i}\right)
_{i\in I}$ such that $\mathcal{P}_{i}\left(  u\right)  =e^{-v_{p_{i}}\left(
u\right)  }$ where $v_{p_{i}}$ is the valuation defined for $u=\left[
u_{\varepsilon}\right]  $ by%
\[
v_{p_{i}}\left(  u\right)  =v_{p_{i}}\left(  \left(  u_{\varepsilon}\right)
_{\varepsilon}\right)  =\sup\left\{  b\in\mathbb{R}:p_{i}\left(
u_{\varepsilon}\right)  =O\left(  \varepsilon^{b}\right)  \text{ as
}\varepsilon\rightarrow0\right\}  .
\]
The valuation on the ring $\widetilde{\mathbb{C}}$ of generalized numbers
given for each $r=\left[  r_{\varepsilon}\right]  $ by
\[
v\left(  r\right)  =v\left(  \left(  r_{\varepsilon}\right)  _{\varepsilon
}\right)  =\sup\left\{  b\in\mathbb{R}:\left\vert u_{\varepsilon}\right\vert
=O\left(  \varepsilon^{b}\right)  \text{ as }\varepsilon\rightarrow0\right\}
\]
leads to the ultra-pseudo-norm on $\widetilde{\mathbb{C}}:\left\vert
r\right\vert _{e}=e^{-v\left(  r\right)  }$.

It is proved in \cite{Gar0} that a $\widetilde{\mathbb{C}}$-linear map
$T:\mathcal{G}_{c}(\Omega)\rightarrow\widetilde{\mathbb{C}\text{ }}$ is
continuous for the above topologies if and only if there exists a finite
subset $I_{0}\subset I$ \ and a constant $C>0$ such that, for all
$u\in\mathcal{G}_{c}(\Omega)$%
\[
\left\vert \left\langle T,u\right\rangle \right\vert _{e}\leq C\underset{i\in
I_{0}}{\max}\mathcal{P}_{i}\left(  u\right)  .
\]

In \cite{Gar}, the topological dual $\mathcal{L}\left(  \mathcal{G}_{c}%
(\Omega),\widetilde{\mathbb{C}}\right)  $ is endowed with the topology of
uniform convergence on bounded subsets which is defined by the
ultra-pseudo-seminorms%
\[
\mathcal{P}_{\mathcal{B}}\left(  T\right)  =\underset{u\in\mathcal{B}}{\sup
}\left\vert \left\langle T,u\right\rangle \right\vert _{e}%
\]
with $\mathcal{B}$ varying in the family of all bounded subsets of
$\mathcal{G}_{c}(\Omega)$, i.e., for each $i\in I$, $\underset{u\in
\mathcal{B}}{\sup}\mathcal{P}_{i}\left(  u\right)  <\infty$.

\subsubsection{Localization of $\mathcal{G}$-singularities}

The sheaf embedding $\mathcal{G}\rightarrow\mathcal{L}\left(  \mathcal{G}%
_{c},\widetilde{\mathbb{C}}\right)  $ is defined, for each open set $\Omega$
of $\mathbb{R}^{n}$, by the continuous (as recalled in \cite{Gar}) map%
\[
\mathcal{G}\left(  \Omega\right)  \ni u\mapsto T_{u}\in\mathcal{L}\left(
\mathcal{G}_{c}\left(  \Omega\right)  ,\widetilde{\mathbb{C}}\right)
\]
where $T_{u}$ is defined, for $u=\left[  u_{\varepsilon}\right]
\in\mathcal{G}\left(  \Omega\right)  $ and each $v=\left[  v_{\varepsilon
}\right]  \in\mathcal{G}_{c}\left(  \Omega\right)  $, by
\[
\left\langle T_{u},v\right\rangle =\left[
{\textstyle\int\nolimits_{K}}
u_{\varepsilon}\left(  x\right)  v_{\varepsilon}\left(  x\right)  dx\right]
\in\widetilde{\mathbb{C}}%
\]
where $K$ is an arbitrary compact set containing supp$v$ in its interior.

From Definition $2.9$ in \cite{Gar}, the $\mathcal{G}$-singular support of
$T\in\mathcal{L}\left(  \mathcal{G}_{c}(\Omega),\widetilde{\mathbb{C}}\right)
$ denoted by ($singsupp_{\mathcal{G}}\left(  T\right)  $) is the complement of
the set of all points $x\in\Omega$ such that the restriction of $T$ to some
neighborhood $V$ of $x$ belongs to $\mathcal{G}\left(  V\right)  $. Then we
still have with our standard notations ($\mathcal{A}=\mathcal{L}\left(
\mathcal{G}_{c},\widetilde{\mathbb{C}}\right)  $, $\mathcal{B=G}$)%
\[
singsupp_{\mathcal{G}}\left(  T\right)  =\mathcal{S}_{\mathcal{L}\left(
\mathcal{G}_{c},\widetilde{\mathbb{C}}\right)  }^{\mathcal{G}}\left(
T\right)  .
\]

\subsection{$\mathcal{F}$-local analysis in $(\mathcal{C},\mathcal{E}%
,\mathcal{P})$-algebras $\mathcal{A}\left(  \Omega\right)  $}

\subsubsection{ The algebraic structure of a $(\mathcal{C},\mathcal{E}%
,\mathcal{P})$ algebra}

We summarize the construction of the so-called $(\mathcal{C},\mathcal{E}%
,\mathcal{P})$ algebras \cite{ADJAMMO, JAM1} which generalize many cases met
in the literature. $\mathbb{K}$ is the real or complex field and $\Lambda$ a
set of indices. $\mathcal{C}$ is the factor ring $A/I$ where $I$ is an ideal
of $A$, a given subring of $\mathbb{K}^{\Lambda}$. $(\mathcal{E},\mathcal{P})$
is a sheaf of topological $\mathbb{K}$-algebras on a topological space $X$. A
presheaf of $(\mathcal{C},\mathcal{E},\mathcal{P})$\textbf{ }algebras on $X$
is a presheaf $\mathcal{A}=\mathcal{H}/\mathcal{J}$ of factor algebras where
$\mathcal{J}$ is an ideal of $\mathcal{H}$, a subsheaf of $\mathcal{E}%
^{\Lambda}$. The sections of $\mathcal{H}$ (resp. $\mathcal{J}$) of $X$ have
to verify some estimates given by means of $\mathcal{P}$ and $A$ (resp. $I$).

The above construction needs some technical conditions given in \cite{ADJAMMO}
on the structure of $\mathcal{C}$ and we suppose that for any open set
$\Omega$ in $X$, the algebra $\mathcal{E}(\Omega)$ is endowed with a family
$\mathcal{P}(\Omega)=(p_{i})_{i\in I(\Omega)}$ of semi-norms. Then, we set%
\begin{gather*}
\mathcal{H}\left(  \Omega\right)  =\mathcal{H}_{(A,\mathcal{E},\mathcal{P}%
)}(\Omega)=\{(u_{\lambda})_{\lambda}\in\left[  \mathcal{E}(\Omega)\right]
^{\Lambda}\mid\forall i\in I(\Omega),\left(  (p_{i}(u_{\lambda})\right)
_{\lambda}\in\left\vert A\right\vert \}\text{,}\\
\mathcal{J}\left(  \Omega\right)  =\mathcal{J}_{(I_{A},\mathcal{E}%
,\mathcal{P})}(\Omega)=\left\{  (u_{\lambda})_{\lambda}\in\left[
\mathcal{E}(\Omega)\right]  ^{\Lambda}\mid\forall i\in I(\Omega),\left(
p_{i}(u_{\lambda})\right)  _{\lambda}\in\left\vert I_{A}\right\vert \right\}
\text{.}%
\end{gather*}

The factor $\mathcal{H}_{(A,\mathcal{E},\mathcal{P})}/\mathcal{J}%
_{(I_{A},\mathcal{E},\mathcal{P})}$ is a presheaf verifying the localization
principle $\left(  F_{1}\right)  $ but generally not the gluing principle
$\left(  F_{2}\right)  $.\textit{ The element in }$\mathcal{A}(\Omega)$
\textit{defined by }$\left(  u_{\lambda}\right)  _{\lambda\in\Lambda}%
\in\mathcal{H}_{(A,\mathcal{E},\mathcal{P})}(\Omega)$ \textit{is denoted by
}$\left[  u_{\lambda}\right]  $. For $u\in\mathcal{A}\left(  \Omega\right)  $,
the notation $\left(  u_{\lambda}\right)  _{\lambda\in\Lambda}\in u$ means
that $\left(  u_{\lambda}\right)  _{\lambda\in\Lambda}$ is a representative of
$u$.\medskip

\subsubsection{ Association process}

We assume further that $A$ is unitary and $\Lambda$ is left-filtering for the
given (partial) order relation $\prec$. Let us denote by:\smallskip

$\bullet$ $\mathcal{F}$ a given sheaf of topological $\mathbb{K}$-vector
spaces (resp. $\mathbb{K}$-algebras) over $X$ containing $\mathcal{E}$ as a
subsheaf,\vspace{-0.04in}

$\bullet$ $a$ a map from $\mathbb{R}_{+}$ to $A_{+}$ such that $a(0)=1$ (for
$r\in\mathbb{R}_{+}$, we denote $a\left(  r\right)  $ by $\left(  a_{\lambda
}\left(  r\right)  \right)  _{\lambda}$).\smallskip

\noindent For $\left(  u_{\lambda}\right)  _{\lambda}\in\mathcal{H}%
_{(A,\mathcal{E},\mathcal{P})}\left(  \Omega\right)  $, we denote by
$\lim\limits_{\Lambda}\left.  _{\mathcal{F}(\Omega)}\right.  u_{\lambda}$ the
limit of $\left(  u_{\lambda}\right)  _{\lambda}$ for the $\mathcal{F}%
$\textit{-}topology when it exists. We recall that $\lim\limits_{\Lambda
}\left.  _{\mathcal{F}(V)}\right.  u_{\lambda}=f\in\mathcal{F}(\Omega)$ iff,
for each $\mathcal{F}$\textit{-}neighbourhood $W$\ of $f$, there exists
$\lambda_{0}\in\Lambda$\ such that%
\[
\lambda\prec\lambda_{0}\Longrightarrow u_{\lambda}\in W.
\]
We suppose also that we have, for each open subset $V\subset\Omega$,%
\[
\mathcal{J}_{(I_{A},\mathcal{E},\mathcal{P})}(V)\subset\left\{  \left(
v_{\lambda}\right)  _{\lambda}\in\mathcal{H}_{(A,\mathcal{E},\mathcal{P}%
)}(V):\lim\limits_{\Lambda}\left.  _{\mathcal{F}(V)}\right.  v_{\lambda
}=0\right\}  .
\]
Now, consider $u=\left[  u_{\lambda}\right]  \in\mathcal{A}(\Omega)$,
$r\in\mathbb{R}_{+}$, $V$ an open subset of $\Omega$ and $f\in\mathcal{F}%
(V)$.\ We say that $u$ is $a\left(  r\right)  $\emph{-associated to }$f$\emph{
in }$V$:%
\[
u\overset{a(r)}{\underset{\mathcal{F}\left(  V\right)  }{\sim}}f
\]
if $\lim\limits_{\Lambda}\left.  _{\mathcal{F}(V)}\right.  \left(  a_{\lambda
}\left(  r\right)  \,u_{\lambda}\left\vert _{V}\right.  \right)
=f.$\smallskip\newline In particular, if $r=0$, $u$ and $f$ are said
\emph{associated }in $V$. To ensure the independence of the definition with
respect to the representative of $u$, we must have, for any $(\eta_{\lambda
})_{\lambda}\in\mathcal{J}_{(I_{A},\mathcal{E},\mathcal{P})}(\Omega)$,
$\lim\limits_{\Lambda}\left.  _{\mathcal{F}(V)}\right.  a_{\lambda}\left(
r\right)  \,\eta_{\lambda}\left\vert _{V}\right.  =0$. As $\mathcal{J}%
_{(I_{A},\mathcal{E},\mathcal{P})}(V)$ is an ideal over $A$, $(a_{\lambda
}\left(  r\right)  \,\eta_{\lambda}\left\vert _{V}\right.  )_{\lambda}$ is in
$\mathcal{J}_{(I_{A},\mathcal{E},\mathcal{P})}(V)$.\ Thus, our claim follows
from the above assumption.

When taking $X=\mathbb{R}^{d}$, $\mathcal{F}=\mathcal{D}^{\prime}$,
$\Lambda=]0,1]$, $\mathcal{A=G}$, $V=\Omega$, $r=0$, the usual association
between $u=\left[  u_{\varepsilon}\right]  \in\mathcal{G}\left(
\Omega\right)  $ and $T\in\mathcal{D}^{\prime}\left(  \Omega\right)  $ is
defined by%
\[
u\sim T\Longleftrightarrow u\overset{a(0)}{\underset{\mathcal{D}^{\prime
}\left(  \Omega\right)  }{\sim}}T\Longleftrightarrow\lim_{\varepsilon
\rightarrow0}\left.  _{\mathcal{D}^{\prime}\left(  \Omega\right)  }\right.
u_{\varepsilon}=T.
\]
Using the previous notations and according to the previous assumption we have,
for any open set $\Omega$ in $X$%
\[
\mathcal{J}_{(I_{A},\mathcal{E},\mathcal{P})}(\Omega)\subset\mathcal{N}%
_{\mathcal{E}}^{\mathcal{F}}\left(  \Omega\right)  =\left\{  \left(
u_{\lambda}\right)  _{\lambda}\in\mathcal{H}_{(A,\mathcal{E},\mathcal{P}%
)}(\Omega)\ :\ \lim\limits_{\Lambda}\left.  _{\mathcal{F}(\Omega)}\right.
u_{\lambda}=0\right\}  .
\]
Set%
\[
\mathcal{F}_{\mathcal{A}}(\Omega)=\left\{  u\in\mathcal{A}(\Omega
)\,\mid\,\exists\left(  u_{\lambda}\right)  _{\lambda}\in u,\ \exists
f\in\mathcal{F}(\Omega)\ :\ \lim\limits_{\Lambda}\left.  _{\mathcal{F}%
(\Omega)}\right.  u_{\lambda}=f\right\}  .
\]
$\mathcal{F}_{\mathcal{A}}(\Omega)$ is well defined because if $\left(
\eta_{\lambda}\right)  _{\lambda}$ belongs to $\mathcal{J}_{(I_{A}%
,\mathcal{E},\mathcal{P})}(\Omega)$, we have $\lim\limits_{\Lambda}\left.
_{\mathcal{F}(\Omega)}\right.  \eta_{\lambda}=0$. \smallskip Moreover,
$\mathcal{F}_{\mathcal{A}}$ is a subpresheaf of vector spaces (resp. algebras)
of $\mathcal{A}$. Roughly speaking, it is the presheaf whose sections above
some open set $\Omega$ are the generalized functions in $\mathcal{A}\left(
\Omega\right)  $ associated to an element of $\mathcal{F}\left(
\Omega\right)  $.

\subsubsection{\textbf{ }Localization of $\mathcal{F}$-singularities}

We refer to definition and results given in section 2 and take here
$\mathcal{B}=\mathcal{F}_{\mathcal{A}}$. When $u$ belongs to $\mathcal{A}%
\left(  \Omega\right)  $, we can consider the set $\mathcal{O}_{\mathcal{A}%
}^{\mathcal{F}}\left(  u\right)  $ ($=\mathcal{O}_{\mathcal{A}}^{\mathcal{B}%
}\left(  u\right)  $) of all $x\in\Omega$ having a neighborhood $V$ on which
$u$ is associated to $f\in\mathcal{F}\left(  V\right)  $, that is:%
\[
\mathcal{O}_{\mathcal{A}}^{\mathcal{F}}\left(  u\right)  =\left\{  x\in
\Omega\ \left\vert \ \exists V\in\mathcal{V}_{x}:u\left\vert _{V}\right.
\in\mathcal{F}_{\mathcal{A}}(V)\right.  \right\}  ,
\]
$\mathcal{V}_{x}$ being the set of all the neighborhoods of $x$. \smallskip
This leads to the following definition: The $\mathcal{F}$-singular support of
$u\in\mathcal{A}(\Omega)$\ is denoted $\mathcal{S}_{\mathcal{A}}^{\mathcal{F}%
}\left(  u\right)  $ and defined as%
\[
\mathcal{S}_{\mathcal{A}}^{\mathcal{F}}\left(  u\right)  =\Omega
\backslash\mathcal{O}_{\mathcal{A}}^{\mathcal{F}}\left(  u\right)  .
\]
Since the support of $u\in\mathcal{A}(\Omega)$ is defined by%
\[
\mathrm{supp}(u)=\Omega\backslash\mathcal{O}_{\mathcal{A}}\left(  u\right)
\ \ \text{with}\ \ \ \mathcal{O}_{\mathcal{A}}\left(  u\right)  =\left\{
x\in\Omega\left\vert \ \exists V\in\mathcal{V}_{x}:u\left\vert _{V}\right.
=0\right.  \right\}  ,
\]
it is clear that $\mathcal{S}_{\mathcal{A}}^{\mathcal{F}}\left(  u\right)  $
is a closed subset containing $\mathrm{supp}(u)$. \medskip

\subsubsection{\textbf{ }Some results}

We can directly deduce the algebraic properties of $\mathcal{S}_{\mathcal{A}%
}^{\mathcal{F}}\left(  u\right)  $ (see \cite{ADJAMMO}) from subsection 2.2.
For the differential ones we suppose that $\mathcal{F}$ is a sheaf of
topological differential vector spaces, with continuous differentiation,
admitting $\mathcal{E}$ as a subsheaf of topological differential algebras.
Then the presheaf $\mathcal{A}$ is also a presheaf of differential algebras
with, for any $\alpha\in\mathbb{N}^{n}$ and $u\in\mathcal{A}\left(
\Omega\right)  $,
\[
\partial^{\alpha}u=\left[  \partial^{\alpha}u_{\lambda}\right]  \text{, where
}\left(  u_{\lambda}\right)  _{\lambda}\text{ is any representative of
}u\text{.}%
\]
The independence of $\partial^{\alpha}u$ on the choice of representative
follows directly from the definition of $\mathcal{J}_{(I_{A},\mathcal{E}%
,\mathcal{P})}$.) The behaviour of $\mathcal{S}_{\mathcal{A}}^{\mathcal{F}%
}\left(  u\right)  $ under differential operations is linked to the following

\begin{proposition}
Under the above hypothesis $\partial^{\alpha}$ is a $\mathcal{F}_{\mathcal{A}%
}$-compatible presheaf operator of $\mathcal{A}\left(  \Omega\right)  $.
\end{proposition}

\begin{proof}
Let $V$ be any open set of $\Omega$. We have%
\[
\mathcal{F}_{\mathcal{A}}(V)=\left\{  v\in\mathcal{A}(V)\,\mid\,\exists\left(
v_{\lambda}\right)  _{\lambda}\in v,\ \exists f\in\mathcal{F}(V)\ :\ \lim
\limits_{\Lambda}\left.  _{\mathcal{F}(V)}\right.  v_{\lambda}=f\right\}  .
\]
Le be $v\in\mathcal{F}_{\mathcal{A}}(V)$ Then $\partial^{\alpha}v$ has a
representative $\partial^{\alpha}v_{\lambda}$ verifying%
\[
\lim\limits_{\Lambda}\left.  _{\mathcal{F}(V)}\right.  \partial^{\alpha
}v_{\lambda}=\partial^{\alpha}f\in\mathcal{F}(V)
\]
and $\partial^{\alpha}v$ is in $\mathcal{F}_{\mathcal{A}}(V)$.
\end{proof}

This result permits to obtain in the following subsection all the expected
results on the propagation of local $\mathcal{F}_{\mathcal{A}}$-singularities
under differential operations.

\begin{example}
Taking $\mathcal{E}=\mathrm{C}^{\infty}$; $\mathcal{F}=\mathcal{D}^{\prime}$;
$\mathcal{A}=\mathcal{G}$ leads to the $\mathcal{D}^{\prime}$-singular support
of an element of the Colombeau algebra. this notion is complementary to the
usual concept of local association in the Colombeau sense. We refer the reader
to \cite{JAM0, JAM1} for more details.
\end{example}

\begin{example}
Those following examples are considered for $X=\mathbb{R}^{d}$, $\mathcal{E}%
=\mathrm{C}^{\infty}$ and $\mathcal{A=G}$. \smallskip\newline$\left(
i\right)  $~Take $u\in\sigma_{\Omega}\left(  \mathrm{C}^{\infty}\left(
\Omega\right)  \right)  $, where $\sigma_{\Omega}:\mathrm{C}^{\infty}\left(
\Omega\right)  \rightarrow\mathcal{G}\left(  \Omega\right)  $ is the well
known canonical embedding. Then $\mathcal{S}_{\mathcal{G}}^{\mathrm{C}^{p}%
}(u)=\emptyset$, for all $p\in\overline{\mathbb{N}}$.\smallskip\newline%
$\left(  ii\right)  $~Take $\varphi\in\mathcal{D}\left(  \mathbb{R}\right)  $,
with $\int\varphi\left(  x\right)  \,\mathrm{d}x=1$, and set $\varphi
_{\varepsilon}\left(  x\right)  =\varepsilon^{-1}\varphi\left(  x/\varepsilon
\right)  $. As $\varphi_{\varepsilon}\underset{\mathcal{D}^{\prime}%
(\mathbb{R})}{\overset{\varepsilon\rightarrow0}{\longrightarrow}}\delta$, we
have: $\mathcal{S}_{\mathcal{G}}^{\mathcal{D}^{\prime}}\left(  \left[
\varphi_{\varepsilon}\right]  \right)  =\left\{  0\right\}  $. We note also
that $\mathcal{S}_{\mathcal{G}}^{\mathrm{C}^{p}}\left(  \left[  \varphi
_{\varepsilon}\right]  \right)  =\left\{  0\right\}  $. (Indeed, for any
$K\Subset\mathbb{R}^{\ast}=\mathbb{R}\backslash\left\{  0\right\}  $ and
$\varepsilon$ small enough, $\varphi_{\varepsilon}$ is null on $K$ and,
therefore, $\varphi_{\varepsilon}\underset{\mathrm{C}^{\infty}(\mathbb{R}%
^{\ast})}{\overset{\varepsilon\rightarrow0}{\longrightarrow}}=0$%
.)\smallskip\newline$\left(  iii\right)  $~Take $u=\left[  u_{\varepsilon
}\right]  $ with $u_{\varepsilon}(x)=\varepsilon\sin(x/\varepsilon)$. We have,
for all $K\Subset\mathbb{R}$, $\lim p_{K,0}(u_{\varepsilon})=0$, for all
$K\Subset\mathbb{R}$, whereas $\lim p_{K,1}(u_{\varepsilon})$ does not exists,
for $l\geq1$.\ Therefore
\[
\mathcal{S}_{\mathcal{G}}^{\mathrm{C}^{0}}(u)=\mathbb{\emptyset\,}%
,\text{\ \ \ }\mathcal{S}_{\mathcal{G}}^{\mathrm{C}^{1}}(u)=\mathbb{R}\text{.}%
\]
Remark that we have, for any $\left(  p,q\right)  \in\overline{\mathbb{N}}%
^{2}$, with $p\leq q$, and $u\in\mathcal{G}$, $\mathcal{S}_{\mathcal{G}%
}^{\mathrm{C}^{p}}\left(  u\right)  \subset\mathcal{S}_{\mathcal{G}%
}^{\mathrm{C}^{q}}(u)$.
\end{example}

\begin{example}
: $\mathcal{D}_{3\sigma-1}^{\prime}$-$\mathbf{local}$ $\mathbf{analysis}$
$\mathbf{in}$\textbf{ }$\mathcal{G}^{\sigma}\left(  \Omega\right)  \medskip$

In \cite{Gram} Gramchev proves the embedding of some spaces of
ultradistributions in $\mathcal{G}\left(  \Omega\right)  $. In \cite{BouBen}
Benmeriem and Bouzar prove the imbedding of $E_{3\sigma-1}^{\prime}\left(
\Omega\right)  $ (the Gevrey ultradistibutions with compact support and
$3\sigma-1$ indice) into $\mathcal{G}^{\sigma}\left(  \Omega\right)  $ (which
is in fact in relation with the index $2\sigma-1$). The imbeding of
$\mathcal{D}_{3\sigma-1}^{\prime}(\Omega)$ (the Gevrey ultradistibution of
$3\sigma-1$ indice) into $\mathcal{G}^{\sigma}\left(  \Omega\right)  $ is also
proved. However if $u\in\mathcal{G}^{\sigma}\left(  \Omega\right)  $ it is
possible to define an association with an ultradistribution\ (for example of
$\mathcal{D}_{3\sigma-1}^{\prime}(\Omega)$) in the following way: for
$T\in\mathcal{D}_{3\sigma-1}^{\prime}(\Omega)$ and $\left[  u_{\varepsilon
}\right]  =u\in\mathcal{G}^{\sigma}\left(  \Omega\right)  $ we set%
\[
u\thicksim T\Longleftrightarrow\underset{\varepsilon\underset{\mathcal{D}%
_{3s-1}^{\prime}}{\rightarrow}0}{\lim}u_{\varepsilon}=T.
\]
It suffices to verify that $\mathcal{N}_{\mathcal{G}^{\sigma}}^{\mathcal{D}%
_{3\sigma-1}^{\prime}}(\Omega)=\left\{  \left(  u_{\varepsilon}\right)
_{\varepsilon}\in\mathcal{E}_{m}^{\sigma}\left(  \Omega\right)  ,\underset
{\varepsilon\underset{\mathcal{D}_{3\sigma-1}^{\prime}}{\rightarrow}0}{\lim
}u_{\varepsilon}=0\right\}  \ $contains $\mathcal{N}^{\sigma}(\Omega)$ to
ensure that the previous definition don't depends upon the representative
$\left(  u_{\varepsilon}\right)  _{\varepsilon}$ of $u\in\mathcal{G}^{\sigma
}\left(  \Omega\right)  $. The subspace of $\mathcal{G}^{s}\left(
\Omega\right)  $, the Gevrey generalized ultradistributions associated to
ultradistributions $\mathcal{D}_{3\sigma-1}^{\prime}(\Omega)$ is
\[
\mathcal{D}_{3\sigma-1,\mathcal{G}^{\sigma}}^{\prime}(\Omega)=\left\{
u=\left[  u_{\varepsilon}\right]  \in\mathcal{G}^{\sigma}\left(
\Omega\right)  ,\exists T\in\mathcal{D}_{3\sigma-1}^{\prime}(\Omega
)\ \underset{\varepsilon\underset{\mathcal{D}_{3\sigma-1}^{\prime}%
}{\rightarrow}0}{\lim}u_{\varepsilon}=T\right\}
\]
$\mathcal{D}_{3\sigma-1,\mathcal{G}^{\sigma}}^{\prime}(\Omega)$ is well
defined because the limit doesn't depend on the representative of $u$.

One can consider $\mathcal{O}_{\mathcal{G}^{\sigma}}^{\mathcal{D}_{3\sigma
-1}^{\prime}}(u)$ \textquotedblleft the set of all $x$ in the neighbourhood of
which $u$ is associated to an ultradistribution\textquotedblright\ :
\[
\mathcal{O}_{\mathcal{G}^{\sigma}}^{\mathcal{D}_{3\sigma-1}^{\prime}%
}(u)=\left\{  x\in\Omega/\exists V\in\mathcal{V}(x):u\mid_{V}\in
\mathcal{D}_{3\sigma-1}^{\prime}(\Omega)\right\}
\]
where $\mathcal{V}(x)$ is the set of all the open neighbourhoods of $x$. The
$\mathcal{D}_{3s-1}^{\prime}$-asymptotic singular support of\textit{ }%
$u\in\mathcal{G}^{s}\left(  \Omega\right)  $ is obtained by taking
$\mathcal{A}=\mathcal{G}^{\sigma}$ and $\mathcal{B}=\mathcal{D}_{3\sigma
-1}^{\prime}$%
\[
S_{\mathcal{G}^{\sigma}}^{\mathcal{D}_{3\sigma-1}^{\prime}}(u)=\Omega
\backslash\mathcal{O}_{\mathcal{G}^{\sigma}}^{\mathcal{D}_{3\sigma-1}^{\prime
}}(u).
\]

\end{example}

\subsection{$\mathcal{B}$-compatibility of differential or pseudo-differential
operators}

Here is a list of particular cases of subpresheaf $\mathcal{B}$ of presheaf
$\mathcal{A}$ of interest to us:%

\[%
\begin{array}
[c]{lllllllllll}%
\mathcal{B}: & \mathrm{C}^{\infty} & \mathrm{C}^{L} & \mathcal{G}^{L} &
\mathcal{G}^{\infty} & \mathcal{G}^{\mathcal{R}} & \mathcal{G}^{\sigma
,\mathcal{\infty}} & \mathcal{G}^{r,\mathcal{R}} & \mathcal{G}^{r,\mathcal{R}%
,L} & \mathcal{G} & \mathcal{F}_{\mathcal{A}}\\
\mathcal{A}: & \mathcal{D}^{\prime} & \mathcal{D}^{\prime} & \mathcal{G} &
\mathcal{G} & \mathcal{G} & \mathcal{G}^{\sigma} & \mathcal{G}^{r} &
\mathcal{G}^{r} & \mathcal{L}\left(  \mathcal{G}_{c}(\Omega),\widetilde
{\mathbb{C}}\right)  & \mathcal{A}%
\end{array}
\]

$\mathcal{B}$ is always a presheaf (or a sheaf) of differential algebras and
$\mathcal{A}$ is a differential $\mathcal{B}$-module with a differentiation
$\partial^{\alpha}$ ($\alpha\in\mathfrak{A}=\mathbb{N}^{n}$) extending the
$\mathcal{B}$-one. Then in each case and each open set $V$ in $\Omega$ (open
set of $X=\mathbb{R}^{n}$) it is easy to prove that $\partial^{\alpha}v$ maps
$\mathcal{B}\left(  V\right)  $ into itself. Thus $\partial^{\alpha}$ is a
presheaf operator $\mathcal{B}$-compatible in $\mathcal{A}\left(
\Omega\right)  $ according to Definition \ref{B-compatible}. If we give now a
family $\left(  b_{\alpha}\right)  _{\alpha\in\mathbb{N}^{n}}$ of elements in
$\mathcal{B}\left(  \Omega\right)  $, then, $P(\partial)$=$\underset
{\left\vert \alpha\right\vert \leq m}{%
{\textstyle\sum}
}b_{\alpha}\partial^{\alpha}$ is a $\mathcal{B}$-compatible operator in
$\mathcal{A}\left(  \Omega\right)  $ from Corollary \ref{sum}.

Moreover at least in some cases, when $\mathcal{A}$ is $\mathcal{D}^{\prime}$,
(resp.$\mathcal{G}$, $\mathcal{L}\left(  \mathcal{G}_{c}(\Omega),\widetilde
{\mathbb{C}}\right)  $), a \textit{pseudo-differential} \textit{operator }$A$
can be defined in $\mathcal{A}\left(  \Omega\right)  $. Setting $b(x,\xi
)=\underset{\left\vert \alpha\right\vert \leq m}{\left(  2\pi\right)  ^{-n}%
{\textstyle\sum}
}b_{\alpha}\left(  x\right)  \left(  i\xi\right)  ^{\alpha}$ when $b_{\alpha}$
belongs to \textrm{C}$^{\infty}\left(  \Omega\right)  $, the differential
operator $P(\partial)$=$\underset{\left\vert \alpha\right\vert \leq m}{%
{\textstyle\sum}
}b_{\alpha}\partial^{\alpha}$ verifying the formula
\[
P(\partial)u\left(  x\right)  =%
{\textstyle\int\nolimits_{\mathbb{R}^{n}}}
e^{ix\xi}b(x,\xi)\widehat{u}(\xi)d\xi=%
{\textstyle\int}
{\textstyle\int\nolimits_{\Omega\times\mathbb{R}^{n}}}
e^{i\left\langle \left(  x-y\right)  ,\xi\right\rangle }b(x,\xi)u(y)dyd\xi
\]
maps $\mathcal{D}\left(  \Omega\right)  $ into $\mathrm{C}^{\infty}\left(
\Omega\right)  $. Via the theory of oscillatory integrals, the above formula
can be extended into%
\[
Au\left(  x\right)  =%
{\textstyle\int}
{\textstyle\int\nolimits_{\mathbb{R}^{n}}}
e^{i\left\langle \left(  x-y\right)  ,\xi\right\rangle }a(x,y,\xi)u(y)dyd\xi
\]
which defines a pseudo-differential operator $A$ \cite{HorFIO} when
$a=a(x,y,\xi)$ in $S_{\rho,\delta}^{m}\left(  \Omega\times\Omega
\times\mathbb{R}^{n}\right)  $ of H\"{o}rmander symbols of order $m$ and type
$\left(  \rho,\delta\right)  $. $A$ extends continuously to a map
$\mathcal{E}^{\prime}\left(  \Omega\right)  \rightarrow\mathcal{D}^{\prime
}\left(  \Omega\right)  $.

We can find (for example in \cite{Dapic}) a definition of generalized
pseudo-differential operators with generalized symbols extending the classical
one. The set $\widetilde{\mathcal{S}}_{\rho,\delta}^{m}\left(  \Omega
\times\mathbb{R}^{n}\right)  $ of generalized symbols can be described
\cite{Gar} as the algebra $\mathcal{G}_{\mathcal{S}_{\rho,\delta}^{m}}\left(
\Omega\times\mathbb{R}^{n}\right)  $ based on $\mathcal{S}_{\rho,\delta}%
^{m}\left(  \Omega\times\mathbb{R}^{n}\right)  $ and obtained as a
$\mathcal{G}_{E}$-module by choosing $E=\mathcal{S}_{\rho,\delta}^{m}\left(
\Omega\times\mathbb{R}^{n}\right)  $. Then, the pseudo differential operator
with generalized symbols $b\in\widetilde{\mathcal{S}}_{\rho,\delta}^{m}\left(
\Omega\times\mathbb{R}^{n}\right)  $ is the map $\mathcal{G}_{c}\left(
\Omega\right)  \rightarrow\mathcal{G}\left(  \Omega\right)  $ given by%
\[
Au:=%
{\textstyle\int\nolimits_{\mathbb{R}^{n}}}
e^{ix\xi}b(x,\xi)\widehat{u}(\xi)d\xi:=\left[
{\textstyle\int\nolimits_{\mathbb{R}^{n}}}
e^{ix\xi}b_{\varepsilon}(x,\xi)\widehat{u}_{\varepsilon}(\xi)d\xi\right]
\]

One can define more generally pseudo differential operators by means of
symbols in $\widetilde{\mathcal{S}}_{\rho,\delta}^{m}\left(  \Omega
\times\Omega\times\mathbb{R}^{n}\right)  $ and generalized oscillatory
integrals (see \cite{GaGrOb}). Such an operator $A$ is given by%
\[
Au:=%
{\textstyle\int}
{\textstyle\int\nolimits_{\Omega\times\mathbb{R}^{n}}}
e^{i\left\langle \left(  x-y\right)  ,\xi\right\rangle }a(x,y,\xi
)u(y)dyd\xi:=\left[
{\textstyle\int}
{\textstyle\int\nolimits_{\Omega\times\mathbb{R}^{n}}}
e^{i\left\langle \left(  x-y\right)  ,\xi\right\rangle }a_{\varepsilon
}(x,y,\xi)u_{\varepsilon}(y)dyd\xi\right]
\]
which defines a generalized function in $\mathcal{G}\left(  \Omega\right)  $
when $u$ is in $\mathcal{G}_{c}\left(  \Omega\right)  $.

We can find in (\cite{Gar}, def. 2.5) an extension of the action of $A$ to the
dual $\mathcal{L}\left(  \mathcal{G}_{c}(\Omega),\widetilde{\mathbb{C}%
}\right)  $ namely%
\[
AT\left(  u\right)  =T\left(  ^{t}Au\right)  ,u\in\mathcal{G}_{c}(\Omega)
\]
where $^{t}A$ (the transpose of $A)$ is the pseudo-differential operator
defined by%
\[
^{t}Au:=%
{\textstyle\int}
{\textstyle\int\nolimits_{\Omega\times\mathbb{R}^{n}}}
e^{i\left\langle \left(  x-y\right)  ,\xi\right\rangle }a(x,y,-\xi
)u(y)dyd\xi\text{.}%
\]

When $\mathcal{A}$ is $\mathcal{D}^{\prime}$, (resp.$\mathcal{G}$,
$\mathcal{L}\left(  \mathcal{G}_{c}(\Omega),\widetilde{\mathbb{C}}\right)  $),
it is proved in each case that if $A$ \textit{is a properly supported
pseudo-differential operator, }it maps $\mathcal{A}\left(  \Omega\right)  $
into itself. Moreover, when $\mathcal{B}$ is $\mathrm{C}^{\infty}$,
(resp.$\mathcal{G}^{\infty}$, $\mathcal{G}$), for each open set $V$ in
$\Omega$, $A$ maps $\mathcal{B}\left(  V\right)  $ into itself. In other words
$A$ is a $\mathcal{B}$-compatible operator in $\mathcal{A}\left(
\Omega\right)  $.

Therefore proposition \ref{operator} allows to deduce the classical inclusions%
\begin{align*}
\mathcal{S}_{\mathcal{A}}^{\mathcal{B}}(P(\partial)u)  &  \subset
\mathcal{S}_{\mathcal{A}}^{B}(u)\\
\text{or }\mathcal{S}_{\mathcal{A}}^{\mathcal{B}}(Au)  &  \subset
\mathcal{S}_{\mathcal{A}}^{B}(u)
\end{align*}
from the presheaf property of an operator $\mathcal{B}$-compatible.

Through Proposition \ref{product}, we can even obtain some non linear results,
when $\mathcal{A}$ is a presheaf of algebras as
\begin{align*}
\mathcal{S}_{\mathcal{A}}^{\mathcal{B}}(\underset{\left\vert \alpha\right\vert
\leq m}{%
{\textstyle\sum}
}b_{\alpha}\left(  \partial^{\alpha}u\right)  ^{p_{\alpha}})  &
\subset\mathcal{S}_{\mathcal{A}}^{B}(u)\\
\text{or }\mathcal{S}_{\mathcal{A}}^{\mathcal{B}}(Au)^{p}  &  \subset
\mathcal{S}_{\mathcal{A}}^{B}(u)
\end{align*}
where $p$ is any positive integer and $\left(  p_{\alpha}\right)  _{\alpha
\in\mathbb{N}^{n}}$ any given family of positive integers.

\begin{remark}
In Definition 1.3 of \cite{Gar}, a functional $T\in\mathcal{L}\left(
\mathcal{G}_{c}\left(  \Omega\right)  ,\widetilde{\mathbb{C}}\right)  $ is
said to be basic if it is of the form $\left\langle T,u\right\rangle =\left[
\left\langle T_{\varepsilon},u_{\varepsilon}\right\rangle \right]  $ where
$\left(  T_{\varepsilon}\right)  _{\varepsilon}$ is a net of distributions in
$\mathcal{D}^{\prime}\left(  \Omega\right)  $ satisfying the following
condition%
\[
\forall K\Subset\Omega\text{ }\exists j\in\mathbb{N}\text{ }\exists\eta
\in\left(  0,1\right]  \text{ }\forall u\in\mathcal{D}_{K}\left(
\Omega\right)  \text{ }\forall\varepsilon\in\left(  0,\eta\right]  \left\vert
T_{\varepsilon}\left(  u\right)  \right\vert \leq\underset{x\in K,\text{
}\left\vert \alpha\right\vert \leq j}{\sup}\left\vert \partial^{\alpha
}u\left(  x\right)  \right\vert .
\]
When $T$ is a basic functional in $\mathcal{L}\left(  \mathcal{G}_{c}%
(\Omega),\widetilde{\mathbb{C}}\right)  $, Theorem 2.10\ proves that the
inclusion%
\[
singsupp_{\mathcal{G}}\left(  AT\right)  \subset singsupp_{\mathcal{G}}\left(
T\right)
\]
is valid for any pseudo differential operator $A$ with amplitude in
$\widetilde{\mathcal{S}}_{\rho,\delta}^{m}\left(  \Omega\times\Omega
\times\mathbb{R}^{n}\right)  $ and a fortiori when $A$ is properly supported.
In this case our above remark on the presheaf property of an operator
$\mathcal{B}$-compatible shows that this result remains valid for any
functional in $\mathcal{L}\left(  \mathcal{G}_{c}\left(  \Omega\right)
,\widetilde{\mathbb{C}}\right)  $.
\end{remark}

\section{The \ frequential microlocal analysis}

After a short overview of classical results in distribution theory and
propagation of singularities under linear and nonlinear operators, we will
give a characterisation of the local regularity in the two more general cases
of generalized functions which summarize the other ones. This leads to the
definition of corresponding wave front sets which gives a general frequential
microanalysis of generalized singularities.

\subsection{Microlocal analysis in distribution spaces}

\subsubsection{Wave front set and microlocal regularity of product}

As it is recalled in Introduction it follows from Lemma 8.1.1. in
\cite{HorPDOT1} that if $\Omega$ is an open set in $\mathbb{R}^{n}$ and
$u\in\mathcal{D}^{\prime}\left(  \Omega\right)  $, one can set: $\Sigma
_{x}(u)=\underset{\Phi}{\cap}\Sigma\left(  \Phi u\right)  ;$ $\Phi
\in\mathcal{D}\left(  \Omega\right)  ,$ $\Phi\left(  x\right)  \neq0$ and
define the wave front set of $u$ as
\[
WF(u)=\left\{  \left(  x,\xi\right)  \in\Omega\times\left(  \mathbb{R}%
^{n}\backslash0\right)  ;\xi\in\Sigma_{x}(u)\right\}  .
\]

Then, if $u\in\mathcal{D}^{\prime}\left(  \Omega\right)  $ and $\left(
x,\xi\right)  \in\Omega\times\left(  \mathbb{R}^{n}\backslash0\right)  $, $u$
is said to be in $H_{\left(  s\right)  }^{loc}$ at $\left(  x,\xi\right)  $ if
$\left(  x,\xi\right)  \notin WF(u-v)$ for some $v\in H_{\left(  s\right)
}\left(  \mathbb{R}^{n}\right)  $. The microlocal regularity of products is
proved in\medskip

\noindent\textbf{Theorem }(8.3.3\ in \cite{HorLNLHDE}): \textit{Let }$u_{j}\in
H_{\left(  s_{j}\right)  }\left(  \mathbb{R}^{n}\right)  $, $j=1,2.$%
\textit{Then}

$\left(  i\right)  u_{1}u_{2}\in H_{\left(  s_{2}\right)  }^{loc}$
\textit{outside} $WF(u_{1})$ \textit{if }$s_{1}>n/2$ \textit{and} $s_{1}%
+s_{2}>n/2.$

$\left(  ii\right)  u_{1}u_{2}\in H_{\left(  s\right)  }^{loc}$
\textit{outside} $WF(u_{1})$ \textit{if} $s_{1}<n/2$ \textit{and} $s_{1}%
+s_{2}-n/2>s\geq0.$

$\left(  iii\right)  u_{1}u_{2}\in H_{\left(  s_{1}+s_{2}-n/2\right)  }^{loc}$
\textit{outside} $WF(u_{1})\cup WF(u_{2})$ \textit{if} $s_{1}+s_{2}>0.$

\subsubsection{Pseudo-differential operators}

We recall that the space $\mathcal{S}_{\rho,\delta}^{m}\left(  \mathbb{R}%
^{n}\times\mathbb{R}^{n}\right)  $ of symbols of order $m$ and type $\left(
\varrho,\delta\right)  $ consist in all $a\in\mathrm{C}^{\infty}\left(
\mathbb{R}^{n}\times\mathbb{R}^{n}\right)  $ such that
\[
\left\vert \partial_{x}^{\beta}\partial_{\xi}^{\alpha}a(x,\xi)\right\vert \leq
C_{\alpha\beta}(1+\left\vert \xi\right\vert ^{m-\varrho\left\vert
\alpha\right\vert +\delta\left\vert \beta\right\vert }%
\]
with $0\leq\delta<\rho\leq1$. We consider here the simplified case
$\mathcal{S}^{m}=\mathcal{S}_{1,0}^{m}$ which defines a pseudo-differential
operator (belonging to $Op\mathcal{S}^{m}$)%
\[
a\left(  x,D\right)  u\left(  x\right)  =%
{\textstyle\int}
e^{i\left\langle x,\xi\right\rangle }a\left(  x,\xi\right)  \hat{u}\left(
\xi\right)  d\xi
\]
which maps continuously $\mathcal{S}\left(  \mathbb{R}^{n}\right)  $ into
$L^{\infty}\left(  \mathbb{R}^{n}\right)  \cap\mathrm{C}\left(  \mathbb{R}%
^{n}\right)  $. It extends for every $s\in\mathbb{R}$ to a continuous map
$H_{\left(  s\right)  }\left(  \mathbb{R}^{n}\right)  \rightarrow H_{\left(
s-m\right)  }\left(  \mathbb{R}^{n}\right)  $.

If $\Omega$ is open then a continuous linear map $A:\mathcal{D}\left(
\Omega\right)  \rightarrow\mathrm{C}^{\infty}\left(  \Omega\right)  $ is said
to be a pseudo-differential operator of order $m$ in $\Omega$ (an element of
$\Psi^{m}\left(  \Omega\right)  $) if for any $\varphi$,$\psi\in
\mathcal{D}\left(  \Omega\right)  $ the operator $\mathcal{S}\left(
\mathbb{R}^{n}\right)  \ni u\mapsto\varphi A\left(  \psi u\right)  $ is in
$Op\mathcal{S}^{m}$. If $u\in\mathcal{E}^{\prime}\left(  \Omega\right)  $ and
$A\in\Psi^{m}\left(  \Omega\right)  $, then $Au\in\mathcal{D}^{\prime}\left(
\Omega\right)  $ is well defined and leads to the inclusion%
\[
\text{sing supp}Au\subset\text{sing supp}u
\]
which is the projection of the microlocal property%
\[
WF\left(  Au\right)  \subset WF\left(  u\right)
\]
when $A$ is properly supported and $u$ belongs to $\mathcal{D}^{\prime}\left(
\Omega\right)  $.

The propagation of non characteristic regularity for semi-linear equations
studied in \cite{Rauch} is given by\medskip

\noindent\textbf{Theorem }(8.4.13\ in \cite{HorLNLHDE}): \textit{Let }$u\in
H_{\left(  s+k\right)  }^{loc}\left(  \Omega\right)  $, \textit{where }%
$\Omega$\textit{ is an open set in }$R^{n}$\textit{ and }$s>n/2$\textit{, be a
solution of the semi-linear equation}%
\[
P(x,D)u=f\left(  x,J_{k}u\right)
\]
\textit{where }$J_{k}u=\left(  \partial^{\alpha}u\right)  _{\left\vert
\alpha\leq k\right\vert }$\textit{, }$f$\textit{ and the coefficients of }%
$P$\textit{ are }$C^{\infty}$\textit{ and }$k$\textit{ is smaller than the
order }$m$\textit{ of }$P(x,D)$\textit{. If }$P$\textit{ is noncharacteristic
at }$\left(  x,\xi\right)  \in\Omega\times\left(  \mathbb{R}^{n}%
\backslash0\right)  $\textit{, it follows that }$u\in H_{\left(
2s+m-n/2\right)  }^{loc}$\textit{ at }$\left(  x,\xi\right)  $\textit{.}

\subsubsection{Application of paradifferential calculus}

The paradifferentiel calculus of Bony \cite{Bony I} is based on some
regularization of non smooth symbols. We don't intend to develop this theory
here but look at it as a powerfull tool to prove good results for nonlinear
problems such as the following\medskip

\noindent\textbf{Theorem }(10.3.6\ in \cite{HorLNLHDE}): \textit{Let }$u\in
H_{\left(  s+m-1/2\right)  }^{loc}\left(  \Omega\right)  $, $s>\max(\left(
n-1\right)  /2$ , $n/4)$, a\textit{nd assume that }$u$\textit{ verifies the
quasilinear differential equation}
\[%
{\textstyle\sum\limits_{\left\vert \alpha\right\vert =m}}
a_{\alpha}(x,J_{m-1}u(x))\partial^{\alpha}u+c(x,J_{m-1}u(x))=
\]
\textit{where }$a_{\alpha}$\textit{ and }$c$\textit{ are }$C^{\infty}%
$\textit{. Then it follows that }$u\in H_{\left(  2s+m-n/2\right)  }^{loc}%
$\textit{ at every non characteristic point }$\left(  x,\xi\right)
$\textit{.}

\subsubsection{Propagation of singularities}

Roughly speaking we know that under some conditions for linear or
pseudo-differential equations, the singularities of solutions propagate along
bicharacteristics This remains valid for nonlinear equations in the sense of

\noindent\textbf{Theorem }(11.4.1\ in \cite{HorLNLHDE}: \textit{Let }$u\in
H_{\left(  s+m\right)  }^{loc}\left(  \Omega\right)  ,s>n/2+1,$\textit{be a
real valued solution of the differential equation }%
\[
F(x,J_{m}u(x))=0
\]
\textit{where }$F\in\mathrm{C}^{\infty}$\textit{. If }$\sigma\leq
2s-n/2$\textit{, then the set of }$\left(  x,\eta\right)  \in\Omega
\times\left(  \mathbb{R}^{n}\backslash0\right)  $\textit{ where }$u\notin
H_{\left(  \sigma+m-1\right)  }^{loc}$\textit{ is contained in the
characteristic set and it is invariant under the Hamilton flow defined by the
principal symbol of the linearized equation}.

Beals \cite{Beals I, Beals II} has studied the case of second order hyperbolic
equations (extended by \cite{Bony II} to arbitrary order) for which we can
give a special version as\medskip

\noindent\textbf{Theorem }(11.5.10\ in \cite{HorLNLHDE}: \textit{Let }$u\in
H_{\left(  s\right)  }^{loc}\left(  \Omega\right)  ,s>n/2$ \textit{be a
solution of of the hyperbolic second order semi-linear equation}%
\[
P(x,\partial)=f(x,u)
\]
\textit{where} $f\in\mathrm{C}^{\infty}$. \textit{If }$u\in H_{\left(
s\right)  }^{loc}$ \textit{at a characteristic point }$\left(  x,\xi\right)  $
\textit{and if }$s\leq\sigma<3s-n+1,$\textit{it follows that }$u\in H_{\left(
s\right)  }^{loc}$\textit{ at the bicharacteristic }$\gamma$\textit{ through
}$\left(  x,\xi\right)  $.

\subsection{The \ frequential microlocal analysis in $\mathcal{G}^{r}$}

We refer the reader to \cite{ADRReg, ADJAM} for more details.

\subsubsection{Characterization of $\mathcal{G}^{r,\mathcal{R}}$-local
regularity}

We consider an open subset $\Omega$of $\mathbb{R}^{d}$ and the Schwartz space
$\mathcal{S}\left(  \Omega\right)  $ of rapidly decreasing functions defined
on $\Omega$, endowed with the family of seminorms $\mathcal{Q}\left(
\Omega\right)  =\left(  \mu_{q,\alpha}\right)  _{\left(  q,\alpha\right)
\in\mathbb{N}\times\mathbb{N}^{d}}$ defined by
\[
\mu_{q,\alpha}\left(  f\right)  =\sup_{x\in\Omega}\left(  1+\left\vert
x\right\vert \right)  ^{q}\left\vert \partial^{\alpha}f\left(  x\right)
\right\vert .
\]

The space of "rough"rapidly decreasing functions can be defined as%
\[
\mathcal{S}_{\ast}\left(  \Omega\right)  =\left\{  f\in\mathrm{C}^{\infty
}\left(  \Omega\right)  \,\left\vert \,\forall q\in\mathbb{N},\;\;\mu
_{q,0}\left(  f\right)  <+\infty\right.  \right\}  .
\]

In order to make easier the comparison between the distributional case and the
generalized case, we begin by recalling the classical theorem and complete it
by some equivalent statements given in the following result (Theorem 16\ in
\cite{ADRReg}): for $u$ in $\mathcal{E}^{\prime}\left(  \mathbb{R}^{n}\right)
$, the following equivalences hold:%
\[%
\begin{tabular}
[c]{lll}%
$\left(  i\right)  ~\,u\in\mathrm{C}^{\infty}\left(  \mathbb{R}^{n}\right)  $
& $\Leftrightarrow\left(  ii\right)  $ & $\!\!\mathcal{F}\left(  u\right)
\in\mathcal{S}\left(  \mathbb{R}^{n}\right)  $\\
& $\Leftrightarrow\left(  iii\right)  $ & $\!\!\mathcal{F}\left(  u\right)
\in\mathcal{S}_{\ast}\left(  \mathbb{R}^{n}\right)  $\\
& $\Leftrightarrow\left(  iv\right)  $ & $\!\!\mathcal{F}\left(  u\right)
\in\mathcal{O}_{M}^{\prime}\left(  \mathbb{R}^{n}\right)  $\\
& $\Leftrightarrow\left(  v\right)  $ & $\!\!\mathcal{F}\left(  u\right)
\in\mathcal{O}_{C}^{\prime}\left(  \mathbb{R}^{n}\right)  .$%
\end{tabular}
\ \ \ \ \
\]
where $\mathcal{F}$ is the classical Fourier transform defined as topological
automorphism of $\mathcal{S}\left(  \mathbb{R}^{d}\right)  $. The result is
based on the following inclusions%
\begin{align*}
\mathcal{S}\left(  \mathbb{R}^{n}\right)   &  \subset\mathcal{S}_{\ast}\left(
\mathbb{R}^{n}\right)  \subset\mathcal{O}_{M}^{\prime}\left(  \mathbb{R}%
^{n}\right)  \subset\mathcal{O}_{C}^{\prime}\left(  \mathbb{R}^{n}\right)  ,\\
\mathcal{F}\left(  \mathcal{E}^{\prime}\left(  \mathbb{R}^{n}\right)  \right)
&  \subset\mathcal{O}_{C}\left(  \mathbb{R}^{n}\right)  ;\mathcal{O}%
_{C}\left(  \mathbb{R}^{n}\right)  \cap\mathcal{O}_{C}^{\prime}\left(
\mathbb{R}^{n}\right)  =\mathcal{S}\left(  \mathbb{R}^{n}\right)  .
\end{align*}

The Fourier transform has been extended to some spaces of rapidly decreasing
generalized functions (like $\mathcal{G}_{s}\left(  \mathbb{R}^{n}\right)
=\mathcal{X}_{\mathcal{S}}\left(  \mathbb{R}^{n}\right)  /\mathcal{N}%
_{\mathcal{S}}\left(  \mathbb{R}^{n}\right)  $) and more completely described
in \cite{ADRReg} in the framework of $\mathcal{R}$-regular spaces. We can
point out that in any framework, the elements with compact support have always
a Fourier transform.

\begin{definition}
Let $\mathcal{R}$ be a regular subset of $\mathbb{R}_{+}^{\mathbb{N}}$ and
$\Omega$ an open subset or $\mathbb{R}^{n}$.\ Set%
\begin{align*}
\mathcal{X}_{\mathcal{S}_{\ast}}^{r,\mathcal{R}}\left(  \Omega\right)   &
=\left\{  \left(  f_{\varepsilon}\right)  _{\varepsilon}\in\mathcal{S}_{\ast
}\left(  \Omega\right)  ^{\Lambda}\,\left\vert \,\exists N\in\mathcal{R}%
,\ \forall q\in\mathbb{N},\;\;\mu_{q,0}\left(  f_{\varepsilon}\right)
=\mathrm{O}\left(  r_{\lambda}^{-N(q)}\right)  \;\mathrm{as}\;\lambda
\rightarrow0\right.  \right\}  ,\\
\mathcal{N}_{\mathcal{S}_{\ast}}^{r}\left(  \Omega\right)   &  =\left\{
\left(  f_{\varepsilon}\right)  _{\varepsilon}\in\mathcal{S}_{\ast}\left(
\Omega\right)  ^{\Lambda}\,\left\vert \,\forall N\in\mathbb{R}_{+}%
^{\mathbb{N}},\;\forall q\in\mathbb{N},\;\;\mu_{q,0}\left(  f_{\varepsilon
}\right)  =\mathrm{O}\left(  r_{\lambda}^{N(q)}\right)  \;\mathrm{as}%
\;\lambda\rightarrow0\right.  \right\}  .
\end{align*}

\end{definition}

One can show that $\mathcal{X}_{\mathcal{S}_{\ast}}^{r,\mathcal{R}}\left(
\Omega\right)  $ is a subalgebra of $\mathcal{S}_{\ast}\left(  \Omega\right)
^{\left(  0,1\right]  }$ and that $\mathcal{N}_{\mathcal{S}_{\ast}}^{r}\left(
\Omega\right)  $ is an ideal of $\mathcal{X}_{\mathcal{S}_{\ast}%
}^{r,\mathcal{R}}\left(  \Omega\right)  $. (The proof is similar to that of
Proposition 1 in \cite{ADRReg}).

\begin{definition}
\label{Gs}The algebra $\mathcal{G}_{s_{\ast}}^{r,\mathcal{R}}\left(
\Omega\right)  =\mathcal{X}_{\mathcal{S}_{\ast}}^{r,\mathcal{R}}\left(
\Omega\right)  /\mathcal{N}_{\mathcal{S}_{\ast}}^{r}\left(  \Omega\right)  $
is called the algebra of $\left(  r,\mathcal{R}\right)  $-regular rough
rapidly decreasing generalized functions.
\end{definition}

\begin{theorem}
Let $x_{0}\in\Omega\subset\mathbb{R}^{n}$ and $u\in\mathcal{G}^{r}\left(
\Omega\right)  $. Then, $u$ is $\mathcal{G}^{r,L}$ at $x_{0}$ (in the sense of
definition 2) iff there exist some neighbourhood $W$ of $x_{0}$, some
$\varphi\in\mathcal{D}\left(  W\right)  $, $\varphi\left(  x_{0}\right)
\neq0$, such that $\widehat{\varphi u}\in\mathcal{G}_{s_{\ast}}^{r,\mathcal{R}%
}\left(  \mathbb{R}^{n}\right)  $
\end{theorem}

\begin{proof}
(Sketch). Let $u$ be an element in $\mathcal{G}^{r}\left(  \Omega\right)  $
$\mathcal{G}^{r,\mathcal{R}}$ at $x_{0}$. There exists a neighbourhood $W$ of
$x_{0}$ such that $u\left\vert _{W}\right.  \in\mathcal{G}^{r,\mathcal{R}%
}\left(  W\right)  $. We can extend any given $\varphi\in\mathcal{D}\left(
W\right)  $, $\varphi\left(  x_{0}\right)  \neq0$, into $\widetilde{\varphi
}\in\mathcal{D}\left(  \mathbb{R}^{d}\right)  $ such that $\widetilde{\varphi
}u\in\mathcal{G}_{c}^{r,\mathcal{R}}\left(  \mathbb{R}^{d}\right)
=\mathcal{G}_{c}\left(  \mathbb{R}^{n}\right)  \cap\mathcal{G}^{r,\mathcal{R}%
}\left(  \mathbb{R}^{n}\right)  $.\ We follow the arguments of Theorem 22 in
\cite{ADRReg}, replacing the $\varepsilon$-estimates by the $r_{\lambda}%
$-ones. This leads to prove that $\widehat{\widetilde{\varphi}u}\in
\mathcal{G}_{s_{\ast}}^{r,\mathcal{R}}\left(  \mathbb{R}^{n}\right)  $.
Conversely, this last assertion with the above hypothesis permits to prove
that $\widetilde{\varphi}u\in\mathcal{G}_{c}^{r,\mathcal{R}}\left(
\mathbb{R}^{n}\right)  $, and then there exits a neighbourhood $V$ of $x_{0}$
such that $u\left\vert _{V}\right.  \in\mathcal{G}^{r,\mathcal{R}}\left(
V\right)  $. In this last part one needs to define an inverse Fourier
transform $\mathcal{F}^{-1}$ in $\mathcal{G}_{s_{\ast}}^{r,\mathcal{R}}\left(
\mathbb{R}^{n}\right)  $ for which one introduces the space $\mathcal{G}%
_{\mathcal{B}}^{r,\mathcal{R}}\left(  \mathbb{R}^{n}\right)  $ of $\left(
r,\mathcal{R}\right)  $-regular bounded generalized functions such that
$\mathcal{F}^{-1}\left(  \mathcal{G}_{s_{\ast}}^{r,\mathcal{R}}\left(
\mathbb{R}^{n}\right)  \right)  \subset\mathcal{G}_{\mathcal{B}}%
^{r,\mathcal{R}}\left(  \mathbb{R}^{n}\right)  $. The result follows from the
equality $\mathcal{G}_{\mathcal{B}}^{r,\mathcal{R}}\left(  \mathbb{R}%
^{n}\right)  \cap\mathcal{G}_{c}\left(  \mathbb{R}^{n}\right)  =\mathcal{G}%
^{r,\mathcal{R}}\left(  \mathbb{R}^{n}\right)  \cap\mathcal{G}_{c}\left(
\mathbb{R}^{n}\right)  $.
\end{proof}

\subsubsection{The $\mathcal{G}^{r,\mathcal{R}}$-generalized wave front set}

\begin{definition}
An element $u\in\mathcal{G}^{r}\left(  \Omega\right)  $ is said to be
microlocally $\left(  r,\mathcal{R}\right)  $-regular at $\left(  x_{0}%
,\xi_{0}\right)  \in\Omega\times\left(  \mathbb{R}^{d}\backslash0\right)  $(we
set: $u\in\mathcal{G}^{r,\mathcal{R}}\left(  x_{0},\xi_{0}\right)  $) if there
exist some neighbourhood $W$ of $x_{0}$, some $\varphi\in\mathcal{D}\left(
W\right)  $, $\varphi\left(  x_{0}\right)  \neq0$, some conic neighborhood
$\Gamma$ of $\xi_{0}$ such that $\widehat{\varphi u}\in\mathcal{G}_{s_{\ast}%
}^{r,\mathcal{R}}\left(  \Gamma\right)  $.
\end{definition}

\begin{definition}
The $\mathcal{G}^{r,\mathcal{R}}$-generalized wave front set of $u\in
\mathcal{G}^{r}\left(  \Omega\right)  $, denoted by $WF^{\left(
r,\mathcal{R}\right)  }\left(  u\right)  $ is the complement in $\Omega
\times\left(  \mathbb{R}^{n}\backslash0\right)  $ of the set of all pairs
$\left(  x_{0},\xi_{0}\right)  $ such that $u$ is microlocally $\left(
r,\mathcal{R}\right)  $-regular at $\left(  x_{0},\xi_{0}\right)  $.
\end{definition}

\begin{theorem}
The projection of $WF^{(r,\mathcal{R})}\left(  u\right)  $ in $\Omega$ is
equal to $sing\,supp^{^{(r,\mathcal{R})}}\left(  u\right)  .$
\end{theorem}

The proof follows from the arguments involved in \cite{HorPDOT1} using lemma 8.1.1.

\begin{example}
\medskip

taking $\lambda=\varepsilon\in\left]  0,1\right]  $, $r_{\varepsilon
}=\varepsilon$ and $\mathcal{R}=\mathcal{B}o$ (the set of bounded sequences),
we obtain the $\mathcal{G}^{\infty}$ microlocal analysis of elements in
$\mathcal{G}$ \cite{NePiSc, Scarpa1}.

taking $\lambda=\varepsilon\in\left]  0,1\right]  $ and $r_{\varepsilon
}=\varepsilon$, we obtain for any $\mathcal{R}$ the $\mathcal{G}^{\mathcal{R}%
}$ microlocal analysis of elements in $\mathcal{G}$ \cite{ADRReg}.
\end{example}

\subsubsection{Characterization of $\mathcal{G}^{r,\mathcal{R},L}$-local
regularity}

When starting from previous cases (like $\mathcal{G}$, $\mathcal{G}%
^{\mathcal{R}}$ or $\mathcal{G}^{L}$) the problem is to change simultaneously
the asymptotic scale into a new one, and the $\mathcal{G}^{\infty}$-regularity
subordinated to $L$-conditions into $\mathcal{G}^{\mathcal{R}}$-regularity
subordinated to $L$-conditions. To do that we have to mix carefully the
techniques used in \cite{ADRReg} and \cite{JAM3}. This study is done in
\cite{ADJAM}. In this subsection, we only give the definitions and results
without proofs.

\noindent\textbf{Theorem }\textit{Let }$x_{0}\in\Omega\subset R^{n}$\textit{
and }$u\in G^{r}\left(  \Omega\right)  $\textit{. Then, }$u$\textit{ is
}$G^{r,\mathcal{R},L}$\textit{ at }$x_{0}$\textit{ (in the sense of Definition
2) iff there exist some neighbourhood }$W$\textit{ of }$x_{0}$\textit{, a
compact }$K$\textit{ such }$W\subset K\Subset\Omega$\textit{, a sequence of
functions }$\chi_{k\text{,}}$\textit{each in} $D_{K}\left(  \Omega\right)
$\textit{ and valued in }$\left[  0,1\right]  $\textit{ with }$\chi_{k}%
u=u$\textit{ on }$W$\textit{, a representative }$\left(  u_{\lambda}\right)
_{\lambda}$\textit{of }$u$\textit{, a regular sequence }$N\in R$\textit{, a
positive constant }$c$\textit{, and }$\lambda_{0}\in\Lambda$\textit{ such that
for all }$\xi\in R$%
\begin{equation}
\forall k\in\mathbb{N},\forall\lambda\prec\lambda_{0},\text{ \ }\left\vert
\xi\right\vert ^{k}\left\vert \widehat{u_{k,\lambda}}\left(  \xi\right)
\right\vert \leq cr_{\lambda}^{-N\left(  k\right)  }\left(  cL_{k}\right)
^{k}. \tag{*}\label{star}%
\end{equation}

\subsubsection{The $\mathcal{G}^{r,\mathcal{R},L}$-generalized wave front set}

\begin{definition}
An element $u\in\mathcal{G}^{r}\left(  \Omega\right)  $ is said to be
microlocally $\left(  r,\mathcal{R},L\right)  $-regular at $\left(  x_{0}%
,\xi_{0}\right)  \in\Omega\times\left(  \mathbb{R}^{n}\backslash0\right)  $(we
set: $u\in\mathcal{G}^{r,\mathcal{R},L}\left(  x_{0},\xi_{0}\right)  $) if
there exist a neighborhood $W$ of $x_{0}$, a conic neighborhood $\Gamma$ of
$\xi_{0}$, a sequence $\left(  u_{k}=\chi_{k}u\right)  _{k\in\mathbb{N}}$ of
generalized functions where each $\chi_{k}$ is valued in $\left[  0,1\right]
$ and is in $\mathcal{D}_{K}\left(  \Omega\right)  $, with $W\subset
K\Subset\Omega$, $u_{k\text{ }}$being equal to $u$ in $W$, a sequence
$N\in\mathcal{R}$, a positive constant $c$, and $\lambda_{0}\in\Lambda$ such
that (\ref{star}) holds when $\xi\in\Gamma$.
\end{definition}

\begin{definition}
\label{wavefront}The $\mathcal{G}^{r,\mathcal{R},L}$-generalized wave front
set of u$\in\mathcal{G}^{r}\left(  \Omega\right)  $, denoted by $WF^{\left(
r,\mathcal{R},L\right)  }\left(  u\right)  $ is the complement in
$\Omega\times\left(  \mathbb{R}^{n}\backslash0\right)  $ of the set of all
pairs $\left(  x_{0},\xi_{0}\right)  $ such that $u$ is microlocally $\left(
r,\mathcal{R},L\right)  $-regular at $\left(  x_{0},\xi_{0}\right)  $.
\end{definition}

$WF^{r,\mathcal{R},L}\left(  u\right)  $ is a closed subset of $\Omega
\times\left(  \mathbb{R}^{n}\backslash0\right)  $, and its projection in
$\Omega$ is given by the following result:

\noindent\textbf{Theorem }The projection of $WF^{(r,\mathcal{R},L)}\left(
u\right)  $ in $\Omega$ is equal to $sing\,supp^{^{(r,\mathcal{R},L)}}\left(
u\right)  .$

\begin{example}
Taking $\lambda=\varepsilon\in\left]  0,1\right]  $, $r_{\varepsilon
}=\varepsilon$ , $\mathcal{R}=\mathcal{B}o$ (the set of bounded sequences), we
obtain for any $L$ the $\mathcal{G}^{L}$ microlocal analysis of elements in
$\mathcal{G}$ \cite{JAM3}.

Taking $\lambda=\varepsilon\in\left]  0,1\right]  $, $r_{\varepsilon
}=\varepsilon$ , $\mathcal{R}=\mathcal{B}o$ and $L_{k}=k+1$, we get the
$\mathcal{G}^{A}$ microlocal analysis of elements in $\mathcal{G}$ \cite{PSV}.

Taking $\lambda=\varepsilon\in\left]  0,1\right]  $, $r_{\varepsilon
}=e^{\varepsilon^{-\frac{1}{2s-1}}}$, $\mathcal{R}=\mathcal{B}o$ and
$L_{k}=\left(  k+1\right)  ^{s}$, we obtain the $\mathcal{G}^{s,\infty}$
micolocal analysis of elements in $\mathcal{G}^{s}$ \cite{BouBen}.
\end{example}

\subsubsection{Propagation of singularities under differential (or
pseudo-differential) operators}

a) We can summarize the first investigations in the following results proved
in \cite{ADJAM}\smallskip

\noindent\textbf{Proposition }\textit{Suppose that }$\left(  a,u\right)
$\textit{ is in }$G^{r}\left(  \Omega\right)  \times G^{r}\left(
\Omega\right)  $\textit{, we have\smallskip}

$\left(  i\right)  $\textit{ If }$a\in G^{r,\mathcal{R}}\left(  \Omega\right)
$\textit{ (resp. }$a\in G^{r,\mathcal{R},L}\left(  \Omega\right)  $\textit{),
then }$WF^{(r,\mathcal{R)}}\left(  au\right)  \subset WF^{(r,\mathcal{R}%
)}\left(  u\right)  $

\textit{(resp. }$WF^{(r,\mathcal{R},L)}\left(  au\right)  \subset
WF^{(r,\mathcal{R},L)}\left(  u\right)  $\textit{)\smallskip}

$\left(  ii\right)  $\textit{ }$WF^{(r,\mathcal{R)}}\left(  \partial^{\alpha
}u\right)  \subset WF^{(r,\mathcal{R})}\left(  u\right)  $\textit{ and
}$WF^{(r,\mathcal{R},L)}\left(  \partial^{\alpha}u\right)  \subset
WF^{(r,\mathcal{R},L)}\left(  u\right)  $\textit{.\smallskip}

\noindent\textbf{Proposition}\textit{ Let }$P\left(  \partial\right)
=\underset{\left\vert \alpha\right\vert \leq m}{%
{\textstyle\sum}
}a_{\alpha}\partial^{\alpha}$\textit{ a differential operator in }%
$G^{r}\left(  \Omega\right)  $\textit{.}

\textit{If the coefficients }$a_{\alpha}$\textit{ lie in }$G^{r,\mathcal{R}%
}\left(  \Omega\right)  $\textit{ (resp. in }$G^{r,\mathcal{R},L}\left(
\Omega\right)  $\textit{), then we have}%
\[
WF^{(r,\mathcal{R)}}\left(  P\left(  \partial\right)  u\right)  \subset
WF^{(r,\mathcal{R})}\left(  u\right)  \text{ (resp. }WF^{(r,\mathcal{R}%
,L)}\left(  P\left(  \partial\right)  u\right)  \subset WF^{(r,\mathcal{R}%
,L)}\left(  u\right)  \text{.}%
\]

\noindent b) In the special case of $\mathcal{G}^{\infty}$ singularities of
$\mathcal{G}$, we can quote the results based on pseudodifferential operators
and pseudodifferential techniques. In \cite{GarHorm} analogues of
H\"{o}rmander definition of the distributional wave front set given in
\cite{HorFIO} are obtained by characterizations of generalized wave front set
in terms of intersection over some non-ellipticity domains. This intersection
is taken over all slow scale pseudo-differential operators $a\in
\widetilde{\underline{\mathcal{S}}}_{sc}^{m}\left(  \Omega\times\mathbb{R}%
^{n}\right)  $ (def. 1.1) verifying some other regularity conditions. More
precisely, if $Ell_{sc}\left(  a\right)  $ denote the set of all $\left(
x,\xi\right)  \in\Omega\times T^{\ast}\left(  \Omega\right)  \backslash0$
where $a$ is slow scale micro-elliptic (def. 1.2), Theorem\textbf{
}2.1\ proves that \textit{for all }$u\in G\left(  \Omega\right)  $%
\[
WF_{g}\left(  u\right)  =WF_{sc}\left(  u\right)  :=\underset{a\left(
x,D\right)  u\in\mathcal{G}^{\infty}\left(  \Omega\right)  }{\underset
{a\left(  x,D\right)  \in\text{ }_{pr}\Psi^{0}\left(  \Omega\right)  }{\cap}%
}Ell_{sc}\left(  a\right)  ^{c}%
\]
where\textit{ }$_{pr}\Psi^{0}\left(  \Omega\right)  $\textit{ }denote the set
of all properly supported slow scale operators of order\textit{ }$0$\textit{.}

Another pseudo-differential characterisation of $WF_{g}\left(  u\right)  $ is
given by Theorem 2.1.1 which proves that for all\textit{ }$u\in G\left(
\Omega\right)  $%
\[
WF_{g}\left(  u\right)  =\cap Char\left(  A\right)
\]
where the intersection is taken over all classical properly supported
classical pseudo-differential operators $A$ such that $Au$ belongs to
$G^{\infty}\left(  \Omega\right)  $.\textit{ }

Following these characterizations some refined results on propagation of
singularities can be obtained. For example, Theorem 3.1 proves that \textit{if
}$A=a\left(  x,D\right)  $\textit{ }is a properly supported
pseudo-differential operator with slow scale symbol and $u\in G\left(
\Omega\right)  $%
\[
WF_{g}\left(  Au\right)  \subset WF_{g}\left(  u\right)  \subset WF_{g}\left(
Au\right)  \cup Ell_{sc}\left(  a\right)  ^{c}\text{.}%
\]

\subsection{The frequential microlocal analysis in $\mathcal{L}\left(
\mathcal{G}_{c}(\Omega),\widetilde{\mathbb{C}}\right)  $}

\subsubsection{The generalized wave front set $WF_{\mathcal{G}}\left(
T\right)  $}

Inspired by the results and definitions of \cite{GarHorm} recalled in the
previous subsection, Garetto (Def. 3.3, \cite{Gar}) defines the $\mathcal{G}%
$-wave front set of a functional $T\in\mathcal{L}\left(  \mathcal{G}%
_{c}(\Omega),\widetilde{\mathbb{C}}\right)  $ as%
\[
WF_{\mathcal{G}}\left(  T\right)  :=\underset{a\left(  x,D\right)
T\in\mathcal{G}\left(  \Omega\right)  }{\underset{a\left(  x,D\right)
\in_{pr}\Psi^{0}\left(  \Omega\right)  }{\cap}}Ell_{sc}\left(  a\right)
^{c}\text{.}%
\]

And even the $\mathcal{G}^{\infty}$-wave front set of $T$ is defined in the
same way by replacing $\mathcal{G}\left(  \Omega\right)  $ by $\mathcal{G}%
^{\infty}\left(  \Omega\right)  $. Proposition 3.5 shows that the projection
on $\Omega$ of $WF_{\mathcal{G}}\left(  T\right)  $ is exactly
$singsupp_{\mathcal{G}}T$.

When $A$ is a properly supported pseudo-differential operator with symbol
$a\in\widetilde{\mathcal{S}}_{\rho,\delta}^{m}\left(  \Omega\times
\mathbb{R}^{n}\right)  $, the inclusion%
\[
WF_{\mathcal{G}}\left(  AT\right)  \subset WF_{\mathcal{G}}\left(  T\right)
\]
can be refined by introducting the concept of $\mathcal{G}$-microsupport of
$a$, denoted $\mu supp_{\mathcal{G}}\left(  a\right)  $. It is the complement
of all $\left(  x,\xi\right)  \in\Omega\times T^{\ast}\left(  \Omega\right)
\backslash0$ where $a$ is $\mathcal{G}$-smoothing (Def.3.6). Then we have
(Corollary 3.9)%
\[
WF_{\mathcal{G}}\left(  a\left(  x,D\right)  T\right)  \subset WF_{\mathcal{G}%
}\left(  T\right)  \cap\mu supp_{\mathcal{G}}\left(  a\right)  \text{.}%
\]
This result is reformulated in terms of $\mathcal{G}$-microsupport of the
operator $A$ ($\mu supp_{\mathcal{G}}\left(  A\right)  $) in the form%
\[
WF_{\mathcal{G}}\left(  AT\right)  \subset WF_{\mathcal{G}}\left(  T\right)
\cap\mu supp_{\mathcal{G}}\left(  A\right)
\]
where the $\mathcal{G}$-microsupport of $A$ is defined (Def.3.11) by%
\[
\mu supp_{\mathcal{G}}\left(  A\right)  :=\underset{a\left(  x,D\right)
=A}{\underset{a\in\widetilde{\mathcal{S}}_{\rho,\delta}^{m}\left(
\Omega\times\mathbb{R}^{n}\right)  }{\cap}}\mu supp_{\mathcal{G}}\left(
a\right)  .
\]

\subsubsection{Fourier transform characterisation of $T$ and propagation of
singularities}

When $\varphi\in\mathcal{D}\left(  \Omega\right)  $ and $T\in\mathcal{D}%
^{\prime}\left(  \Omega\right)  $, we recall that the regularity of $\varphi
T$ can be measured by the rapid decay of its Fourier transform in some conic
region $\Gamma\subset\mathbb{R}^{n}\backslash0$. Following this idea, Garetto
introduces the subset $\mathcal{G}_{\mathcal{S},0}\left(  \Gamma\right)  $ of
$\mathcal{G}_{\tau}\left(  \mathbb{R}^{n}\right)  $ (algebra of tempered
generalized functions) such that%
\[
\mathcal{G}_{\mathcal{S},0}\left(  \Gamma\right)  =\left\{  u=\left[
u_{\varepsilon}\right]  \in\mathcal{G}_{\tau}\left(  \mathbb{R}^{n}\right)
\text{ }\forall l\in\mathbb{R}\text{ }\mathbb{\exists N\in N}\text{ }%
\underset{x\in\Gamma}{\text{sup}}\left(  1+\left\vert \xi\right\vert \right)
^{l}\left\vert u_{\varepsilon}\left(  x\right)  \right\vert =O\left(
\varepsilon^{-N}\right)  \text{ as }\varepsilon\rightarrow0\right\}
\]
which is similar to $\mathcal{G}_{s_{\ast}}^{r,\mathcal{R}}\left(
\Gamma\right)  $ introduced in Definition \ref{Gs}.

This leads to the Fourier transform characterization of $T$ given in

\noindent\textbf{Theorem} 3.15\ of \cite{GarHorm} (or 3.10 of \cite{Gar}):
\textit{Let }$T$\textit{ be a basic functional in }$\mathcal{L}\left(
\mathcal{G}_{c}\left(  \Omega\right)  ,\widetilde{\mathbb{C}}\right)
$\textit{. Then }$\left(  x,\xi\right)  \notin WF_{\mathcal{G}}\left(
T\right)  $\textit{ if and only if there exist a conic neighbourhood of }$\xi
$\textit{ and a cutoff function }$\Phi\in D\left(  \Omega\right)  $\textit{
with }$\Phi\left(  x\right)  =1$\textit{ such that}%
\[
\mathcal{F}\left(  \Phi T\right)  \subset\mathcal{G}_{\mathcal{S},0}\left(
\Gamma\right)  \text{.}%
\]

Then an extension of Theorem 4.1 in \cite{GaGrOb} follows:

\noindent\textbf{Theorem} 4.1 in \cite{Gar}: \textit{If }$A=a\left(
x,D\right)  $is a properly supported pseudo-differential operator with symbol
$a\in\widetilde{\underline{\mathcal{S}}}_{sc}^{m}\left(  \Omega\times
\mathbb{R}^{n}\right)  $ and $T$\textit{ a basic functional in }$L\left(
\mathcal{G}_{c}\left(  \Omega\right)  ,\widetilde{\mathbb{C}}\right)  $, then%
\[
WF_{\mathcal{G}}\left(  AT\right)  \subset WF_{\mathcal{G}}\left(  T\right)
\subset WF_{\mathcal{G}}\left(  AT\right)  \cup Ell_{sc}\left(  a\right)
^{c}\text{.}%
\]

\section{The asymptotic microlocal analysis}

\bigskip Let $\Omega$ be an open set in $X$. Fix $u=\left[  u_{\lambda
}\right]  \in\mathcal{A}(\Omega)$ and $x\in\Omega$. The idea of the
$(a,\mathcal{F)}$\textit{-}microlocal analysis is the following: $\left(
u_{\lambda}\right)  _{\lambda}$ may not tend to a section of $\mathcal{F}$
above a neighborhood of $x$, that is, there may not exist $V\in\mathcal{V}%
_{x}$ and $f\in\mathcal{F}\left(  V\right)  $\ such that $\lim\limits_{\Lambda
}\left.  _{\mathcal{F}(V)}\right.  u_{\lambda}=f$.\ Nevertheless, in this
case, there may exist $V\in\mathcal{V}_{x}$, $r\geq0$ and $f\in\mathcal{F}%
\left(  V\right)  $ such that $\lim\limits_{\Lambda}\left.  _{\mathcal{F}%
(V)}\right.  a_{\lambda}(r)u_{\lambda}=f$, that is $\left[  a_{\lambda
}(r)u_{\lambda}\left\vert _{V}\right.  \right]  $ is in the subspace (resp.
subalgebra) $\mathcal{F}_{\mathcal{A}}(V)$ of $\mathcal{A}(V)$ introduced in
Subsection 2.5.\ These preliminary remarks lead to the following concept and
results which we summarize from the results given in \cite{JAM0, JAM1,
ADJAMMO}.

\subsection{The $(a,\mathcal{F})$\textit{-}singular parametric spectrum}

We recall that $a$ is a map from $\mathbb{R}_{+}$ to $A_{+}$ such that
$a(0)=1$ and $\mathcal{F}$ is a presheaf of topological vector spaces (or
topological algebras). For any open subset $\Omega$ of $X$, $u=\left[
u_{\lambda}\right]  \in\mathcal{A}(\Omega)$ and $x\in\Omega,$ set
\begin{align*}
N_{\left(  a,\mathcal{F}\right)  ,x}\left(  u\right)   &  =\left\{
r\in\mathbb{R}_{+}\ \mid\ \exists V\in\mathcal{V}_{x},\ \exists f\in
\mathcal{F}(V)\ :\ \lim\limits_{\Lambda}\left.  _{\mathcal{F}(V)}\right.
(a_{\lambda}(r)\,u_{\lambda}\left\vert _{V}\right.  )=f\right\} \\
&  =\Big\{r\in\mathbb{R}_{+}\ \mid\ \exists V\in\mathcal{V}_{x}\ :\ \left[
a_{\lambda}\left(  r\right)  u_{\lambda}\left\vert _{V}\right.  \right]
\in\mathcal{F}_{\mathcal{A}}(V)\Big\}.
\end{align*}
It is easy to check that $N_{\left(  a,\mathcal{F}\right)  ,x}\left(
u\right)  $ does not depend on the representative of $u$. If no confusion may
arise, we shall simply write
\[
N_{\left(  a,\mathcal{F}\right)  ,x}\left(  u\right)  =N_{x}(u).
\]

Assume that:

$\left(  a\right)  $ For all $\lambda\in\Lambda$
\[
\forall\left(  r,s\right)  \in\mathbb{R}_{+},\ \ a_{\lambda}(r+s)\leq
a_{\lambda}(r)a_{\lambda}(s),
\]
and, for all $r\in\mathbb{R}_{+}\backslash\left\{  0\right\}  $, the net
$\left(  a_{\lambda}\left(  r\right)  \right)  _{\lambda}$ converges to $0$ in
$\mathbb{K}$

$\left(  b\right)  $ $\mathcal{F}$ is a presheaf of Hausdorff locally convex
topological vector spaces.\smallskip

\noindent Then, from Theorem 7 in \cite{ADJAMMO} we have, for $u\in
\mathcal{A}(\Omega)$:

\medskip

$\left(  i\right)  $~If $r\in N_{x}(u)$, then $\left[  r,+\infty\right)  $ is
included in $N_{x}(u)$. Moreover, for all $s>r$, there exists $V\in
\mathcal{V}_{x}$ such that: $\lim\limits_{\Lambda}\left.  _{\mathcal{F}%
(V)}\right.  (a_{\lambda}(s)\,u_{\lambda}\left\vert _{V}\right.  )=0$.
Consequently, $N_{x}(u)$ is either empty, or a sub-interval of $\mathbb{R}%
_{+}$.

$\left(  ii\right)  $~More precisely, suppose that for $x\in\Omega$, there
exist $r\in\mathbb{R}_{+}$, $V\in\mathcal{V}_{x}$ and$\ f\in\mathcal{F}(V)$,
nonzero on each neighborhood of $x$ included in $V$, such that $\lim
\limits_{\Lambda}\left.  _{\mathcal{F}(V)}\right.  (a_{\lambda}(r)\,u_{\lambda
}\left\vert _{V}\right.  )=f$. Then $N_{x}(u)=\left[  r,+\infty\right)  .$

$\left(  iii\right)  $~In the situation of $\left(  i\right)  $ and $\left(
ii\right)  $, we have that $0\in N_{x}(u)$ iff $N_{x}(u)=\mathbb{R}_{+}%
$.\ Moreover, if one of these assertions holds, the limits $\lim
\limits_{\Lambda}\left.  _{\mathcal{F}(V)}\right.  (a_{\lambda}\left(
s\right)  \,u_{\lambda}\left\vert _{V}\right.  )$ can be non null only for
$s=0$.

\medskip

Now, we set
\begin{gather*}
\Sigma_{\left(  a,\mathcal{F}\right)  ,x}(u)=\Sigma_{x}(u)=\mathbb{R}%
_{+}\backslash N_{x}(u),\\
R_{\left(  a,\mathcal{F}\right)  ,x}\left(  u\right)  =R_{x}(u)=\inf N_{x}(u).
\end{gather*}
According to the previous remarks and comments, $\Sigma_{\left(
a,\mathcal{F}\right)  ,x}(u)$ is an interval of $\mathbb{R}_{+}$ of the form
$\left[  0,R_{\left(  a,\mathcal{F}\right)  ,x}\left(  u\right)  \right)  $ or
$\left[  0,R_{\left(  a,\mathcal{F}\right)  ,x}\left(  u\right)  \right]  $,
the empty set, or $\mathbb{R}_{+}$. This leads to the following

\bigskip

\noindent\textbf{Definition }(4 in \cite{ADJAMMO})\textit{The }$\left(
a,\mathcal{F}\right)  $\textit{-singular spectrum of }$u\in A(\Omega)$\textit{
is the set}%
\[
\mathcal{S}_{\mathcal{A}}^{\left(  a,\mathcal{F}\right)  }\left(  u\right)
=\left\{  (x,r)\in\Omega\times\mathbb{R}_{+}\,\left\vert \,r\in\Sigma
_{x}(u)\right.  \right\}  .
\]

\bigskip

\noindent\textbf{Example }(4 in \cite{ADJAMMO})\textit{Set }$X=\mathbb{R}^{d}%
$\textit{, }$\mathcal{E}=\mathrm{C}^{\infty}$\textit{, }$\mathcal{F}%
=\mathrm{C}^{p}$\textit{ (}$p\in\overline{\mathbb{N}}=\mathbb{N}\cup\left\{
+\infty\right\}  $\textit{), }$f\in\mathrm{C}^{\infty}\left(  \Omega\right)
$\textit{.\ Set }$u=\left[  \left(  \varepsilon^{-1}f\right)  _{\varepsilon
}\right]  $\textit{ and }$v=\left[  \left(  \varepsilon^{-1}\left\vert
\ln\varepsilon\right\vert f\right)  _{\varepsilon}\right]  $\textit{ in
}$\mathcal{A}\left(  \Omega\right)  =\mathcal{G}\left(  \Omega\right)
$\textit{. Then, for all }$x\in\mathbb{R}$\textit{, }%
\[
N_{\left(  a,\mathrm{C}^{p}\right)  ,x}\left(  u\right)  =\left[
1,+\infty\right)  \,,\ \ \text{ }N_{\left(  a,\mathrm{C}^{p}\right)
,x}\left(  v\right)  =\left(  1,+\infty\right)  \,,\ \ \ R_{\left(
a,\mathrm{C}^{p}\right)  ,x}\left(  u\right)  =R_{\left(  a,\mathrm{C}%
^{p}\right)  ,x}\left(  v\right)  =1.
\]

\bigskip

\noindent\textbf{Remark }(5 in \cite{ADJAMMO})\textit{We have: }%
$\Sigma_{\left(  a,\mathcal{F}\right)  ,x}(u)=\varnothing$\textit{ iff
}$N_{\left(  a,\mathcal{F}\right)  ,x}(u)=\mathbb{R}_{+}$\textit{ and,
according to Theorem }7 in\textit{ }\cite{ADJAMMO}\textit{, iff }$0\in
N_{\left(  a,\mathcal{F}\right)  ,x}(u)$\textit{, that is, there exist
}$\left(  V,f\right)  \in V_{x}\times\mathcal{F}(V)$\textit{ such that }%
$\lim\limits_{\Lambda}\left.  _{\mathcal{F}(V)}\right.  (a_{\lambda
}(0)\,u_{\lambda}\left\vert _{V}\right.  )=f$\textit{.\ As }$a_{\lambda
}(0)\equiv1$\textit{, this last assertion is equivalent to }$x\in
\mathcal{O}_{\mathcal{A}}^{\mathcal{F}}\left(  u\right)  $\textit{. Thus
}$\Sigma_{\left(  a,\mathcal{F}\right)  ,x}(u)=\varnothing$\textit{ iff
}$x\notin\mathcal{S}_{\mathcal{A}}^{\mathcal{F}}(u)$\textit{.\medskip}

This remark implies directly the:

\bigskip

\noindent\textbf{Proposition }(8 in \cite{ADJAMMO})\textit{The projection of
the }$\left(  a,\mathcal{F}\right)  $\textit{-singular spectrum of }%
$u$\textit{ on }$\Omega$\textit{ is the }$\mathcal{F}$\textit{-singular
support of }$u$\textit{.}

\subsection{Some properties of the $(a,\mathcal{F})$-singular parametric
spectrum}

\textbf{Notation }For $u=\left[  u_{\lambda}\right]  \in\mathcal{A}\left(
\Omega\right)  $, $\lim\limits_{\Lambda}\left.  _{\mathcal{F}(V)}\right.
\left(  a_{\lambda}(r)\,u_{\lambda}\left\vert _{V}\right.  \right)
\in\mathcal{F}\left(  V\right)  $ means that there exists $f\in\mathcal{F}%
\left(  V\right)  $ such that $\lim\limits_{\Lambda}\left.  _{\mathcal{F}%
(V)}\right.  \left(  a_{\lambda}(r)\,u_{\lambda}\left\vert _{V}\right.
\right)  =f$.

\subsubsection{Linear and differential properties}

It is easy to prove that for any $u,v\in\mathcal{A}(\Omega)$, we have%
\[
\mathcal{S}_{\mathcal{A}}^{\left(  a,\mathcal{F}\right)  }\left(  u+v\right)
\subset\mathcal{S}_{\mathcal{A}}^{\left(  a,\mathcal{F}\right)  }\left(
u\right)  \cup\mathcal{S}_{\mathcal{A}}^{\left(  a,\mathcal{F}\right)
}\left(  v\right)  .
\]

As a corollary: for any $u$, $u_{0}$, $u_{1}$\ in $\mathcal{A}(\Omega)$ with%
\[
\left(  i\right)  \mathit{\ }u=u_{0}+u_{1}\ \ \ \ \ \ (ii)\mathit{\ }%
\mathcal{S}_{\mathcal{A}}^{\left(  a,\mathcal{F}\right)  }\left(
u_{0}\right)  =\varnothing,
\]
we have: $\mathcal{S}_{\mathcal{A}}^{\left(  a,\mathcal{F}\right)  }\left(
u\right)  =\mathcal{S}_{\mathcal{A}}^{\left(  a,\mathcal{F}\right)  }\left(
u_{1}\right)  $.

\medskip

Assume that $\mathcal{F}$ is a sheaf of topological differential vector
spaces, with continuous differentiation, admitting $\mathcal{E}$ as a subsheaf
of topological differential algebras. Then the sheaf $\mathcal{A}$ is also a
sheaf of differential algebras with, for any $\alpha\in\mathbb{N}^{d}$ and
$u\in\mathcal{A}\left(  \Omega\right)  $,
\[
\partial^{\alpha}u=\left[  \partial^{\alpha}u_{\lambda}\right]  \text{, where
}\left(  u_{\lambda}\right)  _{\lambda}\text{ is any representative of
}u\text{.}%
\]
The independence of $\partial^{\alpha}u$ on the choice of representative
follows directly from the definition of $\mathcal{J}_{(I_{A},\mathcal{E}%
,\mathcal{P})}$.) It follows that if \textit{ }$u\ $is in $\mathcal{A}%
(\Omega)$\textit{ }and $g$ in $\mathcal{E}(\Omega)$\textit{, }for all\textit{
}$\partial^{\alpha}$, $\alpha\in\mathbb{N}^{d}$, we have\textit{\ }%
\[
\mathcal{S}_{\mathcal{A}}^{\left(  a,\mathcal{F}\right)  }\left(
g\partial^{\alpha}u\right)  \subset\mathcal{S}_{\mathcal{A}}^{\left(
a,\mathcal{F}\right)  }\left(  u\right)  \text{.}%
\]

This leads to the more general statement: Let $P(\partial)=%
{\displaystyle\sum\limits_{\left\vert \alpha\right\vert \leq m}}
C_{\alpha}\partial^{\alpha}$\ be a differential polynomial with coefficients
in $E(\Omega).$ For any $u\in A(\Omega)$, we have%
\[
\mathcal{S}_{\mathcal{A}}^{\left(  a,\mathcal{F}\right)  }\left(
P(\partial)u\right)  \subset\mathcal{S}_{\mathcal{A}}^{\left(  a,\mathcal{F}%
\right)  }\left(  u\right)  \text{.}%
\]

\subsubsection{Nonlinear properties}

When $\mathcal{F}$ is a presheaf of algebras, the $\left(  a,\mathcal{F}%
\right)  $-singular spectrum inherits new properties with respect to nonlinear
operations. It is the purpose of following results.

\bigskip

\noindent\textbf{Theorem }(15 in \cite{ADJAMMO}) \textit{We suppose that }%
$F$\textit{ is a presheaf or algebras. For }$u$\textit{ and }$v\in A(\Omega
)$\textit{, let }$D_{i}$\textit{ (}$i=1,2,3$\textit{) be the following
disjoint sets:}%
\[
D_{1}=\mathcal{S}_{\mathcal{A}}^{\mathcal{F}}(u)\diagdown(\mathcal{S}%
_{\mathcal{A}}^{\mathcal{F}}(u)\cap\mathcal{S}_{\mathcal{A}}^{\mathcal{F}%
}(v))\ ;\ \ D_{2}=\mathcal{S}_{\mathcal{A}}^{\mathcal{F}}(v)\diagdown
(\mathcal{S}_{\mathcal{A}}^{\mathcal{F}}(u)\cap\mathcal{S}_{\mathcal{A}%
}^{\mathcal{F}}(v))\ ;\ \ D_{3}=\mathcal{S}_{\mathcal{A}}^{\mathcal{F}}%
(u)\cap\mathcal{S}_{\mathcal{A}}^{\mathcal{F}}(v).
\]
\textit{Then the (}$a$\textit{,}$F)$\textit{-singular asymptotic spectrum of
}$uv$\textit{ verifies}%
\[
\mathcal{S}_{\mathcal{A}}^{\left(  a,\mathcal{F}\right)  }\left(  uv\right)
\subset\left\{  (x,\Sigma_{x}(u)),x\in D_{1}\right\}  \cup\left\{
(x,\Sigma_{x}(v)),x\in D_{2}\right\}  \cup\left\{  (x,E_{x}(u,v)),x\in
D_{3}\right\}
\]
\textit{where (for any }$x\in D_{3}$\textit{)}%
\[
E_{x}(u,v)=\left\{
\begin{array}
[c]{l}%
\lbrack0,\sup\Sigma_{x}(u)+\sup\Sigma_{x}(v)]\text{ if }\Sigma_{x}%
(u)\neq\mathbb{R}_{+}\text{ and }\Sigma_{x}(v)\neq\mathbb{R}_{+}\\
\mathbb{R}_{+}\text{ if }\Sigma_{x}(u)=\mathbb{R}_{+}\text{ or }\Sigma
_{x}(v)=\mathbb{R}_{+}%
\end{array}
\right.
\]

\bigskip

\noindent\textbf{Corollary }(16 in \cite{ADJAMMO})\textit{When }$F$\textit{ is
a presheaf of topological algebras, for }$u$\textit{ }$\in A(\Omega)$\textit{
and }$p\in N^{\ast}$\textit{, we have}%
\[
\mathcal{S}_{\mathcal{A}}^{\left(  a,\mathcal{F}\right)  }\left(
u^{p}\right)  \subset\left\{  (x,H_{p,x}(u)),x\in\mathcal{S}_{\mathcal{A}%
}^{\mathcal{F}}(u)\right\}
\]
\textit{where\ }$H_{p,x}(u)=\left\{
\begin{array}
[c]{l}%
\lbrack0,p\sup\Sigma_{x}(u)]\text{ if }\Sigma_{x}(u)\neq\mathbb{R}_{+}\\
\mathbb{R}_{+}\text{ if }\Sigma_{x}(u)=\mathbb{R}_{+}\text{.}%
\end{array}
\right.  $

\subsection{Some examples and applications to partial differential equations}

In this subsection we shall give some examples of $\left(  a,\mathcal{F}%
\right)  $\emph{-}singular spectra of solutions to nonlinear partial
differential equations given in (\cite{ADJAMMO}, subsection 4.2). Throughout
we shall suppose that $\Lambda=]0,1]$, $X=\mathbb{R}^{d}$, $\mathcal{E}%
=\mathrm{C}^{\infty}$, $\mathcal{F}=\mathrm{C}^{p}$ ($1\leq p\leq\infty$) or
$\mathcal{F}={\mathcal{D}}^{\prime}$, $a_{\varepsilon}(r)=\varepsilon^{r}$.
The results will hold for any $(\mathcal{C},\mathcal{E},\mathcal{P})$%
\emph{-}algebra%
\[
\mathcal{A}=\mathcal{H}_{(A,\mathcal{E},\mathcal{P})}/\mathcal{J}%
_{(I_{A},\mathcal{E},\mathcal{P})}%
\]
such that $\left(  a_{\varepsilon}(r)\right)  _{\varepsilon}\in A_{+}$ for all
$r\in\mathbb{R}_{+}$ and the hypothesis given in 2.6.2 holds.

\subsubsection{On the singular spectrum of powers of the delta
function\label{Ex:deltam}}

We can compare the $(a,\mathrm{C}^{p})$-singular spectrum and the
$(a,{\mathcal{D}}^{\prime})$-singular spectrum of powers of the delta
function. Given a mollifier of the form
\[
\varphi_{\varepsilon}\left(  x\right)  =\dfrac{1}{\varepsilon^{d}}%
\varphi\left(  \dfrac{x}{\varepsilon}\right)  ,\ x\in\mathbb{R}^{d}\text{
\ where }\varphi\in\mathcal{D}(\mathbb{R}^{d}),\varphi\geq0\text{
and\ }{\textstyle\int}\varphi\left(  x\right)  dx=1,
\]
its class in ${\mathcal{A}}(\mathbb{R}^{d})$ defines the delta function
$\delta(x)$ as an element of ${\mathcal{A}}(\mathbb{R}^{d})$. Its powers are
given by ($m\in\mathbb{N}$)
\[
\delta^{m}=\big[\varphi_{\varepsilon}^{m}\big]=\Big[\dfrac{1}{\varepsilon
^{md}}\;\varphi^{m}\left(  \dfrac{.}{\varepsilon}\right)  \Big].
\]
Clearly, the $\mathrm{C}^{0}$-singular spectrum is given by
\[
{\mathcal{S}}_{{\mathcal{A}}}^{(a,\mathrm{C}^{0})}(\delta^{m}%
)=\big(0,[0,md]\big).
\]
Differentiating $\varphi^{m}(x)$ and observing that for each derivative there
is a point $x$ at which this function does not vanish we obtain the
$(a,\mathrm{C}^{p})$-singular spectrum of $\delta^{m}:$%
\[
{\mathcal{S}}_{{\mathcal{A}}}^{(a,\mathrm{C}^{p})}(\delta^{m}%
)=\big(0,[0,md+p]\big).
\]

Given now a test function $\psi\in{\mathcal{D}}(\mathbb{R}^{d})$, we have
\[
\int\varphi_{\varepsilon}^{m}(x)\psi(x)\,dx=\int\dfrac{1}{\varepsilon^{md-d}%
}\,\varphi^{m}(x)\psi(\varepsilon x)\,dx,
\]
thus the $(a,D^{\prime})$\textit{-}singular spectrum of $\delta^{m}$ is%
\[
{\mathcal{S}}_{{\mathcal{A}}}^{(a,{\mathcal{D}}^{\prime})}(\delta
^{m})=\varnothing\ \ \mathrm{for}\ m=1,\qquad{\mathcal{S}}_{{\mathcal{A}}%
}^{(a,{\mathcal{D}}^{\prime})}(\delta^{m}%
)=\big(0,[0,md-d[\big)\ \ \mathrm{for}\ m>1.
\]

\subsubsection{The singular spectrum of solutions to semilinear hyperbolic
equations}

The singular spectrum of solutions of a semilinear transport problem%
\[
\left(  P_{\lambda}\right)  \left\{
\begin{array}
[c]{l}%
\partial_{t}u_{\varepsilon}(x,t)+\lambda(x,t)\partial_{x}u_{\varepsilon
}(x,t)=F(u_{\varepsilon}(x,t)),\quad x\in\mathbb{R},\ t\in\mathbb{R}\\
u_{\varepsilon}(x,0)=u_{0\varepsilon}(x),\quad x\in\mathbb{R}%
\end{array}
\right.
\]
where $\lambda$ and $F$ are smooth functions of their arguments. may decrease
or increase with respect to the one of the data, depending on the function
$F$. We observe that by a change of coordinates we may assume without loss of
generality that $\lambda\equiv0$.

\medskip

\noindent a) For $F(u_{\varepsilon}(x,t))=-u_{\varepsilon}^{3}(x,t)$ (the
dissipative case: Example 8 in \cite{ADJAMMO}), the problem $\left(
P_{0}\right)  $ has the solution%
\[
u_{\varepsilon}(x,t)=\frac{u_{0\varepsilon}(x)}{\sqrt{2tu_{0\varepsilon}%
^{2}(x)+1}}=\frac{1}{\sqrt{2t+1/u_{0\varepsilon}^{2}(x)}}.
\]
When the initial data are given by a power of the delta function,
$u_{0\varepsilon}(x)=\varphi_{\varepsilon}^{m}(x)$, the solution formula shows
that $u_{\varepsilon}(x,t)$ is a bounded function (uniformly in $\varepsilon$)
supported on the line $\{x=0\}$. Thus $u_{\varepsilon}(x,t)$ converges to zero
in ${\mathcal{D}}^{\prime}(\mathbb{R}\times]0,\infty\lbrack)$, and so
\[
{\mathcal{S}}_{{\mathcal{A}}}^{(a,{\mathcal{D}}^{\prime})}(u_{0}%
)=\big(0,[0,m-1[\big),\qquad{\mathcal{S}}_{{\mathcal{A}}}^{(a,{\mathcal{D}%
}^{\prime})}(u)=\varnothing.
\]

\noindent b) For $F(u_{\varepsilon}(x,t))=\sqrt{1+u_{\varepsilon}^{2}%
(x,t)},\quad x\in\mathbb{R},\ t>0$ (\cite{ADJAMMO}, Example 9), the problem
$\left(  P_{0}\right)  $ has the solution
\[
u_{\varepsilon}(x,t)=u_{0\varepsilon}(x)\cosh t+\sqrt{1+u_{0\varepsilon}%
^{2}(x)}\,\sinh t.
\]

b$_{1}$) with a delta function as initial value, that is, $u_{0\varepsilon
}(x)=\varphi_{\varepsilon}(x)$ we obtain%
\begin{align*}
\iint u_{\varepsilon}(x,t)\psi(x,t)\,dxdt  &  =\iint\Big(\varphi(x)\cosh
t+\sqrt{\varepsilon^{2}+\varphi^{2}(x)}\,\sinh t\Big)\psi(\varepsilon
x,t)\,dxdt\\
&  \rightarrow\iint\Big(\varphi(x)\cosh t+|\varphi(x)|\sinh t\Big)\psi
(0,t)\,dxdt
\end{align*}
for $\psi\in{\mathcal{D}}(\mathbb{R}^{2})$. Thus in this case
\[
{\mathcal{S}}_{{\mathcal{A}}}^{(a,{\mathcal{D}}^{\prime})}(u_{0}%
)={\mathcal{S}}_{{\mathcal{A}}}^{(a,{\mathcal{D}}^{\prime})}(u)=\varnothing.
\]

b$_{2}$) with the derivative of a delta function as initial value,
$u_{0\varepsilon}(x)=\varphi_{\varepsilon}^{\prime}(x)$, a similar calculation
shows that
\[
\iint u_{\varepsilon}(x,t)\psi(x,t)\,dxdt=\iint\Big(\varphi(x)\cosh
t+\dfrac{1}{\varepsilon}\sqrt{\varepsilon^{4}+(\varphi^{\prime})^{2}%
(x)}\,\sinh t\Big)\psi(\varepsilon x,t)\,dxdt
\]
and so
\[
{\mathcal{S}}_{{\mathcal{A}}}^{(a,{\mathcal{D}}^{\prime})}(u_{0}%
)=\varnothing,\qquad{\mathcal{S}}_{{\mathcal{A}}}^{(a,{\mathcal{D}}^{\prime}%
)}(u)=\{(0,t,r):t>0,0\leq r<1\}.
\]

Example 10 in \cite{ADJAMMO} shows that it is quite possible for the singular
spectrum to increase with time.

\noindent c) when taking $F(u_{\varepsilon}(x,t))=\big(u_{\varepsilon
}(x,t)+1)\big)\log\big(u_{\varepsilon}(x,t)+1\big),\quad x\in\mathbb{R}%
,\ t>0$, the problem $\left(  P_{0}\right)  $ has the solution
\[
u_{\varepsilon}(x,t)=\big(u_{0\varepsilon}(x)+1\big)^{e^{t}},
\]
provided $u_{0\varepsilon}>-1$ in which case the function on the right hand
side of the differential equation is smooth in the relevant region. To
demonstrate the effect, we take a power of the delta function as initial
value, that is $u_{0\varepsilon}(x)=\varphi_{\varepsilon}^{m}(x)$. Then
\[
{\mathcal{S}}_{{\mathcal{A}}}^{(a,{\mathcal{D}}^{\prime})}(u_{0}%
)=\{(0,r):0\leq r<m-1\},\qquad{\mathcal{S}}_{{\mathcal{A}}}^{(a,{\mathcal{D}%
}^{\prime})}(u)=\{(0,t,r):t>0,0\leq r<me^{t}-1\}.
\]

\subsubsection{Blow-up in finite time}

In situations where blow-up in finite time occurs, microlocal asymptotic
methods allow to extract information beyond the point of blow-up. This can be
done by regularizing the initial data and truncating the nonlinear term. This
is shown in Example 11 of \cite{ADJAMMO} for a simple situation.

The problem to be treated is formally the initial value problem%

\[%
\begin{array}
[c]{l}%
\partial_{t}u(x,t)=u^{2}(x,t),\quad x\in\mathbb{R},\ t>0\\
u(x,0)=H(x),\quad x\in\mathbb{R}%
\end{array}
\]
where $H$ denotes the Heaviside function. Clearly, the local solution
$u(x,t)=H(x)/(1-t)$ blows up at time $t=1$ when $x>0$. Choose $\chi
_{\varepsilon}\in\mathrm{C}^{\infty}\left(  \mathbb{R}\right)  $ with%
\[
0\leq\chi_{\varepsilon}(z)\leq1\ ;\ \chi_{\varepsilon}(z)=1\text{ if }%
|z|\leq\varepsilon^{-s}\,,\ \chi_{\varepsilon}(z)=0\ \text{if }|z|\geq
1+\varepsilon^{-s}\,,\ s>0.
\]
Further, let $H_{\varepsilon}(x)=H\ast\varphi_{\varepsilon}(x)$ where
$\varphi_{\varepsilon}$ is a mollifier as in ~\ref{Ex:deltam}. One considers
the regularized problem
\[%
\begin{array}
[c]{l}%
\partial_{t}u_{\varepsilon}(x,t)=\chi_{\varepsilon}\big(u_{\varepsilon
}(x,t)\big)u_{\varepsilon}^{2}(x,t),\quad x\in\mathbb{R},\ t>0\\
u_{\varepsilon}(x,0)=H_{\varepsilon}(x),\quad x\in\mathbb{R}.
\end{array}
\]
When $x<0$ and $\varepsilon$ is sufficiently small, $u_{\varepsilon}(x,t)=0$
for all $t\geq0$. For $x>0$, $u_{\varepsilon}(x,t)=1/(1-t)$ as long as
$t\leq1-\varepsilon^{s}$. The cut-off function is chosen in such a way that
$|\chi_{\varepsilon}(z)z^{2}|\leq(1+\varepsilon^{-s})^{2}$ for all
$z\in\mathbb{R}$. Therefore,
\[
\partial_{t}u_{\varepsilon}\leq(1+\varepsilon^{-s})^{2}%
\ \mbox{always\ and}\ \partial_{t}u_{\varepsilon}%
=0\ \mbox{when}\ |u_{\varepsilon}|\geq1+\varepsilon^{-s}.
\]
Some computations and estimates permit to obtain the following $\mathrm{C}%
^{0}$-singular support and $\left(  a,\mathrm{C}^{0}\right)  $-singular
spectrum (for $a_{\varepsilon}(r)=\varepsilon^{r}$) of $u=\left[
u_{\varepsilon}\right]  $:%
\[
{\mathcal{S}}_{{\mathcal{A}}}^{\mathrm{C}^{0}}(u)={\mathcal{S}}_{1}%
(u)\cup{\mathcal{S}}_{2}(u)\text{ with }{\mathcal{S}}_{1}(u)=\{(0,t):0\leq
t<1\}~;~{\mathcal{S}}_{2}(u)=\{(x,t):x\geq0,t\geq1\},
\]%
\[
{\mathcal{S}}_{{\mathcal{A}}}^{(a,\mathrm{C}^{0})}(u)=\left(  {\mathcal{S}%
}_{1}(u)\times\left\{  0\right\}  \right)  \cup\left(  {\mathcal{S}}%
_{2}(u)\times\left[  0,s\right]  \right)  .
\]

These results give a microlocal precision on the the blow-up: The $C^{0}%
$-singularities (resp. $\left(  a,\mathrm{C}^{0}\right)  $-singularities) of
$u$ are described by means of two sets: $\mathcal{S}_{1}(u)$ and
$\mathcal{S}_{2}(u)$ (resp. $\mathcal{S}_{1}(u)\times\left\{  0\right\}  $ and
$\mathcal{S}_{2}(u)\times\left[  0,s\right]  $). The set $\emph{S}_{1}(u)$
(resp. $\mathcal{S}_{1}(u)\times\left\{  0\right\}  $) is related to the data
$\mathrm{C}^{0}$ (resp. $\left(  a,\mathrm{C}^{0}\right)  $)-singularity. The
set $\mathcal{S}_{2}(u)$ (resp. $\mathcal{S}_{2}(u)\times\left[  0,s\right]
$) is related to the singularity due to the nonlinearity of the equation
giving the blow-up at $t=1$. The blow-up locus is the edge $\left\{
x\geq0,t=1\right\}  $ of $\mathcal{S}_{2}(u)$ and the strength of the blow-up
is measured by the length $s$ of the fiber $\left[  0,s\right]  $ above each
point of the blow-up locus. This length is closely related to the diameter of
the support of the regularizing function $\chi_{\varepsilon}$ and depends
essentially on the nature of the blow-up: Changing simultaneously the scales
of the regularization and of the cut-off (i.e. replacing $\varepsilon$ by some
function $h(\varepsilon)\rightarrow0$ in the definition of $\varphi
_{\varepsilon}$ and $\chi_{\varepsilon}$) does not change the fiber and
characterizes a sort of moderateness of the strength of the blow-up.

\subsubsection{The strength of a singularity and the sum law}

We point out the following remark (\cite{ADJAMMO}, subsection 4.3): when
studying the propagation and interaction of singularities in semilinear
hyperbolic systems, Rauch and Reed \cite{RauchReed} defined the strength of a
singularity of a piecewise smooth function. This notion is recalled in the
one-dimensional case. Assume that the function $f:\mathbb{R}\rightarrow
\mathbb{R}$ is smooth on $]-\infty,x_{0}]$ and on $[x_{0},\infty\lbrack$ for
some point $x_{0}\in\mathbb{R}$. The \emph{strength of the singularity of $f$
at $x_{0}$} is the order of the highest derivative which is still continuous
across $x_{0}$. For example, if $f$ is continuous with a jump in the first
derivative at $x_{0}$, the order is $0$; if $f$ has a jump at $x_{0}$, the
order is $-1$. Travers \cite{Travers} later generalized this notion to include
delta functions. Slightly deviating from her definition, but in line with the
one of \cite{RauchReed}, it is possible to define the strength of singularity
of the $k$-th derivative of a delta function at $x_{0}$, $\partial_{x}%
^{k}\delta(x-x_{0})$, by $-k-2$.

The significance of these definitions is perceived in the description of what
Rauch and Reed termed \emph{anomalous singularities} in semilinear hyperbolic
systems. This effect is demonstrated in a paradigmatic example, also due to
\cite{RauchReed}, the $(3\times3)$-system
\begin{equation}
\left\{
\begin{array}
[c]{rclcl}%
(\partial_{t}+\partial_{x})u(x,t) & = & 0, & \quad & u(x,0)=u_{0}(x)\\
(\partial_{t}-\partial_{x})v(x,t) & = & 0, & \quad & v(x,0)=v_{0}(x)\\
\partial_{t}w(x,t) & = & u(x,t)v(x,t), & \quad & w(x,0)=0
\end{array}
\right.  \tag{**}\label{double star}%
\end{equation}
Assume that $u_{0}$ has a singularity of strength $n_{1}\geq-1$ at $x_{1}=-1$
and $v_{0}$ has a singularity of strength $n_{2}\geq-1$ at $x_{2}=+1$. The
characteristic curves emanating from $x_{1}$ and $x_{2}$ are straight lines
intersecting at the point $x=0$, $t=1$. Rauch and Reed showed that, in
general, the third component $w$ will have a singularity of strength
$n_{3}=n_{1}+n_{2}+2$ along the half-ray $\{(0,t):t\geq1\}$. This half-ray
does not connect backwards to a singularity in the initial data for $w$, hence
the term \emph{anomalous singularity}. The formula $n_{3}=n_{1}+n_{2}+2$ is
called the \emph{sum law}. Travers extended this result to the case where
$u_{0}$ and $v_{0}$ were given as derivatives of delta functions at $x_{1}$
and $x_{2}$. We are going to further generalize this result to powers of delta
functions, after establishing the relation between the strength of a
singularity of a function $f$ at $x_{0}$ and the singular spectrum of
$f\ast\varphi_{\varepsilon}$.

We consider a function $f:\mathbb{R}\rightarrow\mathbb{R}$ which is smooth on
$(-\infty,x_{0}]$ and on $[x_{0},\infty)$ for some point $x_{0}\in\mathbb{R}$;
actually only the local behavior near $x_{0}$ is relevant. A mollifier
$\varphi_{\varepsilon}(x)=\frac{1}{\varepsilon}\varphi(\frac{x}{\varepsilon})$
is fixed as in \ref{Ex:deltam} and the corresponding embedding of
${\mathcal{D}}^{\prime}(\mathbb{R})$ into the $(\mathcal{C},\mathcal{E}%
,\mathcal{P})$\emph{-}algebra ${\mathcal{A}}(\mathbb{R})$ is denoted by
$\iota$. In particular, $\iota(f)=[f\ast\varphi_{\varepsilon}]$.

If $f$ is continuous at $x_{0}$, then $\lim_{\varepsilon\rightarrow0}%
f\ast\varphi_{\varepsilon}=f$ in $\mathrm{C}^{0}$. If $f$ has a jump $x_{0}$,
this limit does not exist in $\mathrm{C}^{0}$, but $\lim_{\varepsilon
\rightarrow0}\varepsilon^{r}f\ast\varphi_{\varepsilon}=0$ in $\mathrm{C}^{0}$
for every $r>0$. The following result is\medskip

\noindent\textbf{Proposition (}16 in \cite{ADJAMMO}) \textit{Let }$x_{0}%
\in\mathbb{R}$\textit{. If }$f:\mathbb{R}\rightarrow\mathbb{R}$\textit{ is a
smooth function on }$(-\infty,x_{0}]$\textit{ and on }$[x_{0},\infty)$\textit{
or }$f(x)=\partial_{x}^{k}\delta(x-x_{0})$\textit{ for some }$k\in\mathbb{N}%
$\textit{, then the strength of the singularity of }$f$\textit{ at }$x_{0}%
$\textit{ is }$-n$\textit{ if and only if }%
\[
\Sigma_{(a,\mathrm{C}^{1}),x_{0}}\big(\iota(f)\big)=[0,n].
\]
\textit{Here }$n\in\mathbb{N}$\textit{ and }$a_{\varepsilon}(r)=\varepsilon
^{r}$\textit{.}

When returning to the model equation we find that the sum law remains valid
when the initial data are powers of delta functions. Suitable $(\mathcal{C}%
,\mathcal{E},\mathcal{P})$-algebras $\mathcal{A}(\mathbb{R})$ and
$\mathcal{A}(\mathbb{R}^{2})$ are exhibited in which the initial value problem
can be uniquely solved. When the scale is taken as $a_{\varepsilon
}(r)=\varepsilon^{r}$, the following result is obtained:\medskip

\noindent\textbf{Proposition }(17 in \cite{ADJAMMO}) \textit{Let }%
$u_{0}(x)=\delta^{m}(x+1)$\textit{, }$v_{0}(x)=\delta^{n}(x-1)$\textit{ for
some }$m,n\in\mathbb{N}^{\ast}$\textit{. Let }$w\in\mathcal{A}(\mathbb{R}%
^{2})$\textit{ be the third component of the solution to problem
(\ref{double star}). Then }$w(x,t)$\textit{ vanishes at all points }%
$(x,t)$\textit{ with }$x\neq0$\textit{ as well as }$(0,t)$\textit{ with }%
$t<1$\textit{, and }%
\[
\Sigma_{(a,\mathrm{C}^{1}),(0,t)}\big(w\big)\subset\lbrack0,m+n]
\]
\textit{for }$t\geq1$\textit{.}

\subsection{Microlocal characterisation of some regular subalgebras}

We recall that the subsheaf ${\mathcal{G}}^{\infty}$ of \emph{regular
Colombeau functions} of the sheaf ${\mathcal{G}}$ is defined as follows
\cite{Ober1}: Given an open subset $\Omega$ of $\mathbb{R}^{d}$, the algebra
${\mathcal{G}}^{\infty}(\Omega)$ comprises those elements $u$ of
${\mathcal{G}}(\Omega)$ whose representatives $(u_{\varepsilon})_{\varepsilon
}$ satisfy the condition
\[
\forall K\Subset\Omega,\ \exists m\in\mathbb{N},\ \forall l\in\mathbb{N}%
:p_{K,l}(u_{\varepsilon})=o(\varepsilon^{-m})\ \mathrm{as}\ \varepsilon
\rightarrow0.
\]

In relation with regularity theory of solutions to nonlinear partial
differential equations, a further subalgebra of ${\mathcal{G}}(\Omega)$ has
been introduced in \cite{OberBiaVolume} -- the algebra of Colombeau functions
of \emph{total slow scale type}. It consists of those elements $u$ of
${\mathcal{G}}(\Omega)$ whose representatives $(u_{\varepsilon})_{\varepsilon
}$ satisfy the condition
\[
\forall K\Subset\Omega,\ \forall r>0,\ \forall l\in\mathbb{N}:p_{K,l}%
(u_{\varepsilon})=o(\varepsilon^{-r})\ \mathrm{as}\ \varepsilon\rightarrow0.
\]
The term \emph{slow scale} refers to the fact that the growth is slower than
any negative power of $\varepsilon$ as $\varepsilon\rightarrow0$. Both
previous properties can be characterized by means of the singular spectrum.

We find in (\cite{ADJAMMO}, subsection 4.4) the proof of the corresponding
characterisations. Let $u\in{\mathcal{G}}(\Omega)$, then:

$\left(  i\right)  $ (Proposition 18) $u$ belongs to ${\mathcal{G}}^{\infty
}(\Omega)$ if and only if $\Sigma_{(a,\mathrm{C}^{\infty}),x}\big(u\big)\neq
\mathbb{R}_{+}$ for all $x\in\Omega$.

$\left(  ii\right)  $ (Proposition 19) $u$ is of total slow scale type if and
only if $\Sigma_{(a,\mathrm{C}^{\infty}),x}\big(u\big)\subset\{0\}$ for all
$x\in\Omega.$

\end{document}